\documentclass[12pt, a4paper, titlepage]{article}
\usepackage[doublespacing]{setspace}
\usepackage[{lmargin=0.75in,rmargin=0.75in, tmargin=1.25in,bmargin=1.25in}]{geometry}
\usepackage[english]{babel}
\usepackage[utf8]{inputenc}
\usepackage{amsmath,amssymb,amsfonts,amsthm}
\usepackage{tikz-cd,algorithm,algpseudocode,mathtools}
\usepackage{graphicx}
\graphicspath{D:\School\UoB Year 3\Group Project - Prime Recognition and Factorization}
\usepackage{fancyhdr,lastpage,chngpage}
\usepackage{hyperref}
\hypersetup{ 
    colorlinks=true,
    linkcolor=blue,
    filecolor=red,      
    urlcolor=red,
    citecolor=magenta
}
\DeclarePairedDelimiter{\abs}{\lvert}{\rvert}
\DeclarePairedDelimiter{\norm}{\lVert}{\rVert}
\makeatletter
\let\oldabs\abs
\def\abs{\@ifstar{\oldabs}{\oldabs*}}
\let\oldnorm\norm
\def\norm{\@ifstar{\oldnorm}{\oldnorm*}}
\makeatother

\newcommand{\reals}{\mathbb{R}}
\newcommand{\nat}{\mathbb{N}}
\newcommand{\zahle}{\mathbb{Z}}
\newcommand{\quot}{\mathbb{Q}}
\newcommand{\compl}{\mathbb{C}}
\newcommand{\nozero}{\backslash\{0\}}
\newcommand{\PriID}{\mathfrak{P}}
\newcommand{\priID}{\mathfrak{p}}
\newcommand{\loh}{\mathbf{o}}
\newcommand{\boh}{\mathbf{O}}
\newcommand{\bom}{\boldsymbol{\Omega}}
\newcommand{\lom}{\boldsymbol{\omega}}

\DeclareMathOperator{\ord}{ord}
\DeclareMathOperator{\Frob}{Frob}
\DeclareMathOperator{\Gal}{Gal}
\DeclareMathOperator{\sign}{sign}
\DeclareMathOperator{\irr}{irr}
\DeclareMathOperator{\Leg}{Leg}
\DeclareMathOperator{\Ind}{Ind}

\algnewcommand\algorithmicinput{\textbf{Input:}}
\algnewcommand\Input{\item[\algorithmicinput]}
\algnewcommand\algorithmicoutput{\textbf{Output:}}
\algnewcommand\Output{\item[\algorithmicoutput]}
\algnewcommand{\IIf}[1]{\State\algorithmicif\ #1\ \algorithmicthen}
\algnewcommand{\EndIIf}{\unskip\ \algorithmicend\ \algorithmicif}

\theoremstyle{definition}
\newtheorem{thm}{Theorem}[section] 
\newtheorem{lemma}[thm]{Lemma}
\newtheorem{prop}[thm]{Proposition}
\newtheorem{corr}[thm]{Corollary}
\newtheorem{conj}[thm]{Conjecture}
\newtheorem{defn}[thm]{Definition} 
\theoremstyle{remark}
\newtheorem*{rmk}{Remark}

\begin{document}
	\pagenumbering{gobble}
	\begin{titlepage}
	\begin{center}
	\Huge{Number Field Sieve with provable complexity}\\
	\Large{Barry van Leeuwen}\\
	\vspace*{1cm}
	\large{Supervisor: Dr. A.R. Booker}\\
	\large{\begin{align*}
	\text{Chair: }& \text{Dr. T. Dokchitser (University of Bristol)}\\	
	\text{Examiners: } & \text{Dr. J. Bober (University of Bristol)}\\
	& \text{Dr. S. Siksek (University of Warwick)}\\	
	\end{align*}}\\
	\normalsize{A dissertation submitted to the}\\
	\large{University of Bristol}\\
	\normalsize{in accordance with the requirements for award of the degree of}\\
	\large{Master of Science by Research	in Mathematics}\\
	\normalsize{ at the }\\
	\large{Faculty of Science}
	\vspace*{1cm}\\
	\Large{School of Mathematics}\\
	\large{\today}\\
	\vspace*{2cm}
	\hfill\normalsize{Word count: 28557}
	\end{center}
	\end{titlepage}
	\newpage
	\pagenumbering{arabic}
	\cfoot{Page \thepage\ of \pageref{LastPage}}
\section*{Abstract}
\addcontentsline{toc}{section}{Abstract}
In this thesis we give an in-depth introduction to the General Number Field Sieve, as it was used by Buhler, Lenstra, and Pomerance, \cite{BLP}, before looking at one of the modern developments of this algorithm: A randomized version with provable complexity. This version was posited in 2017 by Lee and Venkatesan, \cite{RNFS}, and will be preceded by ample material from both algebraic and analytic number theory, Galois theory, and probability theory.
\newpage
\section*{Dedication and Acknowledgements}
I want to thank Dr. Andrew Booker, who as my supervisor managed to find what I needed even though it may not have been what I wanted. I also want to thank Dr. James Milne, Dr. Florian Bouyer, and Dr. Lynne Walling for providing some of the material used. I also want to thank Dr. Dan Fretwell, who helped me find my footing when I just started (which feels very long ago) and my mother, Cokky van Leeuwen, who with her sorcery managed to find typos that I could not.\\
\vspace{5cm}
\begin{center}
\large{\textit{To my wife, Sarah van Leeuwen, who with her continued support and motivation made possible what I thought impossible.}}
\end{center}
\newpage
\section*{Declaration}
I declare that the work in this dissertation was carried out in accordance with the requirements of the University's Regulations and Code of Practice for Research Degree Programmes and that it has not been submitted for any other academic award. Except where indicated by specific reference in the text, the work is the candidate's own work. Work done in collaboration with, or with the assistance of, others, is indicated as such. Any views expressed in the dissertation are those of the author.\\
\\
Signed: Barry van Leeuwen\\
Date: \today
\newpage
\tableofcontents
\addcontentsline{toc}{section}{Contents}
\newpage
\section*{Notation}
\addcontentsline{toc}{section}{Notation}
A brief introduction to some of the notation used. Most of the notation listed will be formally introduced, but this can be used as a reference guide.\\
\begin{table}[h!]
\begin{tabular}{lcl}
$R[\alpha]$ & - & A ring extended by an element $\alpha$ \\
$F(\alpha)$ & - & A field extended by an element $\alpha$ \\
$\mathcal{O}_K$ & - & The ring of algebraic integers in a field extension $K$\\
$f.f.(R)$ & - & The field of fractions generated by a ring $R$\\
$\irr_F(a)$ & - & The minimal polynomial of $a$ over $F$\\
$[a]\in B$	& - & The coset of $a$ in structure $B$\\
$\langle a\rangle\subset B$ & - & The ideal generated by $a$ in structure $B$\\
$\mathcal{L}(s,\chi)$ & - & Dirichlet L-function with character $\chi$ and exponent $s$\\
$a\sim B$ & - & Draw $a$ accoridng to distribution $B$ \\
$\pi(X)$ & - & The number of primes below $X$\\
\end{tabular}
\end{table}
\\
In this paper we use Bachmann--Landau Big-O notation extended by Knuth Big-$\Omega$ notation, such that for functions $f(x)$ and $g(x)$ and $N\in\nat$:
\begin{table}[h!]
\begin{tabular}{lcl}
$f(x)=\loh\left(g(x)\right)$ & - & $\forall k \in\reals_{>0}\exists N$ such that $\forall x\geq N:\abs{f(x)}\leq k\cdot g(x)$ \\
$f(x)=\boh\left(g(x)\right)$ & - & $\exists k\in\reals_{>0}\exists N$ such that $\forall x\geq N:\abs{f(x)}\leq k\cdot g(x)$\\
$f(x)=\lom\left(g(x)\right)$ & - &$\forall k\in\reals_{>0}\exists N$ such that $\forall x\geq N:\abs{f(x)}\geq k\cdot g(x)$\\
$f(x)=\bom\left(g(x)\right)$ & - & $\exists k\in\reals_{>0}\exists N$ such that $\forall x\geq N:\abs{f(x)}\geq k\cdot g(x)$\\
\end{tabular}
\end{table}
\newpage
\section{Introduction}
In 1988 Pollard introduced a brand new factorization algorithm: The Number Field Sieve (NFS). This saw a first implementation in 1994 by Lenstra et al. and was used to factor the ninth Fermat number, $F_9$.\\
\\
The grandure of this algorithm was in the conjectured complexity of
$$L_n\bigg(\frac{1}{3},\sqrt[3]{\frac{64}{9}}+\loh(1)\bigg),$$
where $n\in\nat$ and for $a,b,x\in\reals$: 
$$L_x(a,b)=\exp(b((\log n)^a(\log\log n)^{1-a}).$$
This was a major improvement to the best algorithms at the time, such as the quadratic sieve, which were of complexity $L_n(\frac{1}{2},b)$ for some constant $b$.\\
\\
Despite abundant heuristics indicating the bound and copious research poured into a proof it was very difficult to show the complexity of the two major parts of the sieve:
\begin{enumerate}
\item Using the concept of smoothness to find a factor base.
\item Reducing this factor base to a factorization of a given $n\in\nat$.
\end{enumerate}
In chapter \ref{ch: GNFS} we will be looking at the General Number Field Sieve (GNFS) how it was explicitly described in $\cite{BLP}$ and how the sieve is structured.\\
\\
It wasn't until 2017 that Lee and Venkatesan, \cite{RNFS}. used a probabilistic and combinatorial approach to alter the algorithm in such a way that they managed to prove the following, (\cite{RNFS}, theorem 2.1 \& 2.3): 
\newpage
\begin{thm}
There is a randomised variant of the Number Field Sieve which for each $n$ finds congruences of squares $x^2\equiv y^2 \mod n$ in expected time:
$$L_n\bigg(\frac{1}{3},\sqrt[3]{\frac{64}{9}}+\mathbf{o}(1)\bigg)$$
\end{thm}
\begin{rmk}
Note that no claims are made regarding the factorization of $n$. This theorem only proves that there exists an algorithm for which the sieving process up until finding a congruence of squares, i.e. a pair $(x,y)$ for which $x^2\equiv y^2 \mod n$, with a non-zero probability in the complexity bound given in the theorem above.
\end{rmk}
This version of the algorithm will be extensively studied in chapter \ref{ch:RNFS} and finishes by proving the theorem above. \\
\\
To do this we need to use a strong basis in analytic number theory, which we will introduce in chapter \ref{ch: PREL} along with a review of some probability theory. In chapter \ref{ch: PRE1} we will recall the necessary concepts from algebraic number theory and Galois theory, but for those familiar with these subjects this chapter can be skipped.
\newpage
\section{Algebraic fundamentals}\label{ch: PRE1}
\subsection{Field extensions}\label{sec: PRE1ext}
The reason the General Number Field Sieve is so effective is that it uses strong algebraic constructions to get a factorization, which allows for far fewer computations to be made to get an effective result. The discussion starts by considering these constructions and concepts, such as field extensions, in all generality, but will swiftly collapse down to the specific cases needed for this sieve. We begin by introducing some notation that we will use throughout. 
\label{2.1}
\begin{defn}
\begin{itemize}
\item []
\item Let $K$ and $L$ be fields such that $K\subset L$, then $L$ is called a field extension of $K$ and is denoted $L/K$.
\item $K(\alpha)$ is the smallest field extension of $K$ that contains $\alpha$. If $\alpha\in K$ then $K(\alpha)=K$.
\item $K[x]$ is the polynomial extension of $K$, such that
$$K[x]=\Big\{\sum_{i=0}^\infty a_ix^i\mid \forall i\in\zahle^+:a_i\in K\Big\}.$$
\end{itemize}
\end{defn}
\begin{rmk}
To prevent confusion $\subset$ and $\supset$ are proper sub- and super-constructions, while $\subseteq$ and $\supseteq$ can have equality. 
\end{rmk}
A field extension $L/K$ can be seen as a $K$-vector space and when considered as such it becomes a natural next step to think about the dimension of $L$ as a $K$-vectorspace. This is called the degree and is denoted $\left[L:K\right]=\dim_K(L)$. If there are multiple extensions, $L\supseteq M\supseteq K$, called a tower of extensions, called a tower of extensions then the degrees of these extensions are very well behaved.
\begin{thm}
\label{2.2}\textbf{Tower Law}\\
Let $K_n\supseteq K_{n-1}\supseteq \ldots\supseteq K_0$ be fields. Then
$$\left[K_n:K\right]=\left[K_n:K_{n-1}\right]\ldots\left[K_1:K_0\right].$$
\end{thm}
If $\left[L:K\right]<\infty$ the extension is called finite, else it is infinite. To use the sieve it is important to understand when such an extension is finite and when it isn't. Let $f(x)$ be an irreducible monic monovariate polynomial with coefficients in $K$ of degree $\deg(f)=d$. As $f(x)$ is irreducible it has no roots in $K$, but it might have roots over other fields that contain $K$. 
\begin{defn}
\label{2.3}
An irreducible polynomial $f$ over a field $K$ is said to be separable over $K$ if it has no multiple zeros in a splitting field. That is there exists a field $L\supset K$ such that
$$f(x)=(x-\alpha_1)(x-\alpha_2)\ldots(x-\alpha_d)$$
in $L$, where $f(\alpha_i)=0$.
\end{defn}
This definition is not extremely helpful and in practice it will be more useful to consider [\cite{GALT}, proposition 9.14] instead:
\begin{prop}
\label{2.4}
If $K$ is a subfield of $\compl$ then every irreducible polynomial over $K$ is separable.
\end{prop}
By the fundamental theorem of algebra it is known that a polynomial $f(x)$ with $\deg(f)=d$ has exactly $d$ roots in $\compl$. These roots are not exclusive, and a root $\alpha$ is the root of many polynomials. However there is a unique monic polynomial of lowest degree that is of specific importance:
\begin{defn}
\label{2.5}
Let $f\in K[x]$ be a polynomial with coefficients in a field $K$ for which $\alpha$ is a root. Then $f$ is the minimal polynomial of $\alpha$ in $K$ if 
\begin{itemize}
\item $f$ is monic.
\item $f$ is irreducible.
\item $\forall g\in K[x]$ such that $g(\alpha)=0$, $f$ divides $g$.
\end{itemize}
This polynomial is unique up to units. 
\end{defn}
Throughout this text we adopt the convention to denote the minimal polynomial for $\alpha\in K$ over $K$ as $\irr_K(\alpha)$. The minimal polynomial gives rise to the following equivalencies, which are vital in the background throughout.
\begin{thm}
\label{2.6}
Let $f\in K[x]$ with $\deg(f)=d$ be a  monic irreducible monovariate polynomial with root $\alpha$. Let $\langle f(x)\rangle$ denote the ideal generated by $f(x)$.Then the following are equivalent:
\begin{enumerate}
\item $f(x)$ is the minimal polynomial for $\alpha$ in $K[x]$; $f(x)=\irr_K(\alpha)$.
\item $K(\alpha)\cong K[x]/\langle f(x)\rangle$.
\item $K(\alpha)$ is a field.
\item $[K(\alpha):K]=d$ and $\{1,\alpha,\alpha^2,\ldots,\alpha^{d-1}\}$ is a basis for $K(\alpha)$ over $K$.
\end{enumerate}
\end{thm}
The above theorem reveals a condition for which the extensions $K(\alpha)$ is a finite degree extension of $K$. This condition is called algebraicity and for an extension $L/K$ an element $\alpha\in L$ is called algebraic exactly if it is a root for a polynomial $f(x)\in K[x]$. 
\begin{rmk} There are two special cases to consider:
\begin{itemize}
\item A complex number, $\alpha\in\compl$, that is the root of a polynomial with coefficients in $\zahle$ is called an algebraic integer.
\item A complex number, $\alpha\in\compl$, that is the root of a polynomial with coefficients in $\quot$ is called an algebraic number
\end{itemize}
\end{rmk}
It turns out that the structure of the algebraic numbers is very simple to describe.
\begin{thm}
\label{2.7}
Every algebraic number is of the form $\frac{\alpha}{b}$ where $\alpha$ is algebraic over $\zahle$ and $b\in\zahle\nozero$.
\end{thm}
As just discussed: If $\irr_K(\alpha)$ has degree $d$, then it is only natural to ask about the other roots of $irr_K(\alpha)$. These are called the conjugates of $\alpha$ and have some interesting properties.
\begin{defn}
\label{2.8}
Let $\alpha\in\compl$ be algebraic over $K\subseteq\compl$. The conjugates of $\alpha$, $\alpha=\alpha_1,\ldots,\alpha_d$ are the roots of $\irr_K(\alpha)$.
\end{defn}
As $K\subset\compl$ proposition \ref{2.4} shows that any irreducible polynomial over $K$ is separable, hence each $\alpha_i$ is distinct:
$$\irr_K(\alpha)=(x-\alpha_1)(x-\alpha_2)\ldots(x-\alpha_d).$$
Using the conjugates [\cite{IANT},  theorem 5.6.1] proves that we only need to consider extensions by a singular algebraic number:
\begin{thm}
\label{2.9}
Let $K\subseteq\compl$ and let $\alpha,\beta\in\compl$ be algebraic over $K$. Then there exists $\gamma\in\compl$ algebraic over $K$ such that
$$K(\alpha,\beta)=K(\gamma).$$
\end{thm}
This can easily be extended to any finite extension $K(\alpha_1,\ldots,\alpha_n)$ by using the fact that $K(\alpha_1,\alpha_2)=\left(K(\alpha_1)\right)(\alpha_2)$.
\subsection{Algebraic number fields and Dedekind domains}\label{sec: PRE1alnf}
To work effectively with the general number field sieve there are some restrictions we place on our extensions. The first restriction placed upon the extensions is the restriction that they be finite. Hence, from now on, every extension $L/K\subseteq\compl$ is assumed finite. Special focus will be on algebraic number fields:
\begin{defn}\label{2.10}
Let $K\subseteq\compl$, then $K$ is an algebraic number field if $K=\quot(\alpha_1,\ldots,\alpha_n)$ such that for all $i\in\nat$, $\alpha_i$ is an algebraic number.
\end{defn}
By theorem \ref{2.9} we can equate $K=\quot(\alpha)$ for some algebraic number $\alpha\in\compl$. Moreover we know that any algebraic number can be written as $\alpha=\frac{\beta}{b}$ where $\beta$ is an algebraic integer and $b\in\zahle\nozero$, but $\quot(\frac{\beta}{b})\cong\quot(\beta)$ as $\frac{1}{b}\in\quot$ so it is safe to assume that $\alpha$ is in fact an algebraic integer.
\begin{defn}\label{2.11}
Let $A$ be the set of all algebraic integers and let $K$ be an algebraic number field. Then the set $\mathcal{O}_K=A \cap K$ is called the ring of algebraic integers over $K$.
\end{defn}
As the name suggests, $\mathcal{O}_K$ is a ring. In fact, by [\cite{IANT}, theorem 6.1.4 and 6.1.6],
\begin{thm}\label{2.12}
Let $K$ be an algebraic number field, then $\mathcal{O}_K$ is integrally closed. i.e.
\begin{enumerate}
\item $\mathcal{O}_K$ is an integral domain.
\item $\forall \alpha$: $\alpha$ algebraic over $\mathcal{O}_K \Rightarrow \alpha\in\mathcal{O}_K$.
\end{enumerate}
\end{thm}
In generality $\mathcal{O}_K$ is surprisingly well behaved allowing for some interesting properties, including that it's field of fractions is $K$:
$$f.f.(\mathcal{O}_K)=\Big\{\frac{\alpha}{\beta}\Big\vert\alpha,\beta\in\mathcal{O}_K\Big\}=K.$$
Another of the properties which is important is [\cite{IANT}, theorem 6.5.3]:
\begin{thm}\label{2.13}
Let $K$ be an algebraic number field. Then $\mathcal{O}_K$ is a Noetherian domain.
\end{thm}
\begin{rmk}
A terribly inconvenient, yet common, convention is that the ideals $\mathfrak{I}\subseteq \mathcal{O}_K$ are just called ideals of $K$. Despite the confusion this might causes it will not pose a large problem, and will be clear from context, as in practice $K$ will be a number field, hence only has trivial ideals.
\end{rmk}
It is nearly by definition that a maximal ideal of an integral domain is a prime ideal, but in general the converse is not true. However, in the case of $\mathcal{O}_K$ for an algebraic number field $K$ this converse is true.
\begin{thm} \label{2.14}
Let $K=\quot(\alpha)$ be an algebraic number field. Let $\mathfrak{P}\subseteq\mathcal{O}_K$ be a prime ideal of $\mathcal{O}_K$, then $\mathfrak{P}$ is a maximal ideal of $\mathcal{O}_K$.
\end{thm}
Now we have shown that $\mathcal{O}_K$ is integrally closed, is a Noetherian domain, and each prime ideal of $\mathcal{O}_K$ is a maximal ideal. These three define an important class of domains, which are very useful for the General Number Field Sieve.
\begin{defn}\label{2.15}
An integral domain $D$ that satisfies the following properties:
\begin{itemize}
\item D is a Noetherian domain.
\item D is integrally closed.
\item Each prime ideal of $D$ is a maximal ideal.
\end{itemize}
is called a Dedekind domain.
\end{defn}
Dedekind domains admit a very interesting property, as they generalize the idea of unique factorization into primes which we know by the fundamental theorem of arithmetic holds for all $z\in\zahle$ up to units, $\pm 1$. In Dedekind domains this holds as well, however we must use the equivalent construction for number fields.
\begin{thm}\label{2.16}
Let $D$ be a Dedekind domain. Any non-trivial integral ideal is a product of prime ideals and this factorization is unique up to order, i.e. for $\mathfrak{I}\subset D$ an ideal with factorization
$$\mathfrak{I}=\prod_{i\in I} \mathfrak{P}_i^{e_i}=\prod_{j\in J} \mathfrak{Q}_j^{f_j},$$
with prime ideals $\mathfrak{P}_i,\mathfrak{Q}_j$ then there exists a permutation $\sigma:I\rightarrow J$ such that $\mathfrak{P}_i=\mathfrak{Q}_j$ and $e_i=f_j$.
\end{thm}
To prove this we need some auxiliary measures:
\begin{prop}\label{2.17}
Let $D$ be a Dedekind domain. Then every non-zero ideal $\mathfrak{I}\in\mathcal{O}_K$ contains a product of one or more prime ideals.
\end{prop}
\begin{defn}\label{2.18}
Let $D$ be an integral domain and let $K$ be the quotient field of $D$. For each prime ideal $\mathfrak{P}$ of $D$ we define the set $\bar{\mathfrak{P}}$ by
$$\bar{\mathfrak{P}}=\{\alpha\in K:\alpha \mathfrak{P}\subseteq D\}.$$
This set is an ideal and is called a fractional ideal of $D$.
\end{defn}
\begin{prop}\label{2.19}
Let $D$ be a Dedekind domain. Let $\mathfrak{P}$ be a prime ideal of $D$. Then
$$\PriID\bar{\PriID}=D.$$
\end{prop}
With these facts recorded we proceed with the proof of theorem \ref{2.16}:
\begin{proof}
This proof will come in two parts: Existence and Uniqueness. Let $D$ be a Dedekind domain such that there exists a non-trivial ideal of $D$ that is not a product of prime ideals.\\
\\
\underline{Existance}: As $D$ is a Dedekind domain, it is Noetherian, and so by the maximality principle there is a non-trivial ideal $\mathfrak{A}$ of $D$ that is maximal with respect to the property of not being a product of prime ideals. Then by proposition \ref{2.17} there exists prime ideals $\mathfrak{P}_1,\ldots,\mathfrak{P}_k$ of $D$ such that
$$\mathfrak{P}_1\ldots\mathfrak{P}_k\subseteq \mathfrak{A}.$$
Now let $k$ be the smallest positive integer for which such a product exists. If $k=1$ then $\mathfrak{P}_1\subseteq \mathfrak{A}\subset D$ is a prime ideal, hence it is maximal. So $\mathfrak{A}=\mathfrak{P}_1$. This contradicts our assumption that $\mathfrak{A}$ was not a product of prime ideals. Hence $k\geq 2$. By proposition \ref{2.19} we have that $\bar{\mathfrak{P}_1}\mathfrak{P}_1=D$ and so
$$\bar{\mathfrak{P}_1}\mathfrak{P}_1\ldots\mathfrak{P}_k=D\mathfrak{P}_2\ldots\mathfrak{P}_k.$$
Hence we have that $\bar{\mathfrak{P}_2}\mathfrak{A}\supseteq\mathfrak{P}_2\ldots\mathfrak{P}_k$. Moreover by proposition \ref{2.19} we also have that $\mathfrak{A}\subseteq\bar{\mathfrak{P}_1}\mathfrak{A}$. Now assume that $\mathfrak{A}=\bar{\mathfrak{P}_1}\mathfrak{A}$, then
$$\mathfrak{A}\supseteq \mathfrak{P}_2\ldots\mathfrak{P}_k.$$
Which contradicts the minimality of $k$. Now assume $\mathfrak{A}\subset\bar{\mathfrak{P}_1}\mathfrak{A}$, then as the latter is an ideal of $D$, by the maximality property of $\mathfrak{A}$, we get
$$\bar{\mathfrak{P}_1}\mathfrak{A}=\mathfrak{Q}_2\ldots\mathfrak{Q}_l,$$
for some prime ideals $\mathfrak{Q}_j$, $j\in\{2,\ldots,l\}$. But then
$$\mathfrak{A}=\mathfrak{A}D=\mathfrak{A}\bar{\mathfrak{P}_1}\mathfrak{P}_1=\mathfrak{P}_1\mathfrak{Q}_2\ldots\mathfrak{Q}_l,$$
which contradicts the way $\mathfrak{A}$ was chosen. Hence every ideal of $D$ is a product of prime ideals.\\
\\
\underline{Uniqueness:} Suppose now that there exists a maximal ideal $\mathfrak{A}$ of $D$, such that 
$$\mathfrak{A}=\mathfrak{P}_1\ldots\mathfrak{P}_k=\mathfrak{Q}_1\ldots\mathfrak{Q}_l.$$
By the maximality principle this is a choice we can make. Then as $\mathfrak{P}_1\ldots\mathfrak{P}_k\subseteq\mathfrak{Q}_1$ and $\mathfrak{Q}_1$ is a prime ideal we have that $\mathfrak{P}_i=\mathfrak{Q}_i$ for some $i\in\{1,\ldots,k\}$. Relabelling $\mathfrak{Q}_i$ as $\mathfrak{Q}_1$ we get $\mathfrak{Q}_1=\mathfrak{P}_1$ as $\mathfrak{P}_1$ is a prime, hence maximal, ideal of $D$. Thus
$$\bar{\mathfrak{Q}_1}\mathfrak{Q}_1\ldots\mathfrak{Q}_l=\mathfrak{A}\mathfrak{Q}_1=\mathfrak{A}\bar{\mathfrak{P}_1}=\bar{\mathfrak{P}_1}\mathfrak{P}_1\ldots\mathfrak{P}_k=\mathfrak{P}_2\ldots\mathfrak{P}_k.$$
Now assume that $\mathfrak{A}=\mathfrak{A}\bar{\mathfrak{P}_1}$, then $\mathfrak{A}\bar{\mathfrak{P}_1}\mathfrak{P}_1=\mathfrak{A}\mathfrak{P}_1$ and so $\mathfrak{A}=\mathfrak{A}\mathfrak{P}_1$. Define the fractional ideal of $\mathfrak{A}$ as
$$\bar{\mathfrak{A}}=\bar{\mathfrak{P}_1}\ldots\bar{\mathfrak{P}_k}.$$
Then
$$\mathfrak{A}\bar{\mathfrak{A}}=\mathfrak{P}_1\ldots\mathfrak{P}_k\bar{\mathfrak{P}_1}\ldots\bar{\mathfrak{P}_k}.$$
and so
$$D=\bar{\mathfrak{A}}\mathfrak{A}=\mathfrak{A}\bar{\mathfrak{P}_1}\bar{\mathfrak{A}}=\mathfrak{P}_1.$$
but this contradicts the primality of $\mathfrak{P}_1$, so our equality assumption was faulty. Hence $\mathfrak{A}\subset\mathfrak{A}\bar{\mathfrak{P}_1}$. Since $\mathfrak{A}\bar{\mathfrak{P}_1}$ is a maximal ideal of $D$, $\mathfrak{A}\bar{\mathfrak{P}_1}$ has exactly one factorization as a product of prime ideals. Thus we deduce that $k-1=l-1$ and so $k=l$. Relabbeling this gives that for all $i\in\{1,\ldots,k\}$
$$\mathfrak{P}_i=\mathfrak{Q}_i,$$
proving that the decomposition into prime ideals is unique.
\end{proof}
\subsection{Galois theory and the Frobenius element}\label{sec: PRE1galt}
In this section we give a very concise introduction to a few ideas from Galois Theory. If you have a basic understanding of finite Galois theory then there is no loss in simply skipping this section.
\begin{defn}\label{2.20}
Let $L/K$ be number fields. Then $L$ is said to be normal over $K$ if every irreducible polynomial $f\in K[x]$ that has a zero in $L$ splits over $L$.
\end{defn}
Normality is easily checked under certain conditions:
\begin{defn}\label{2.21}
Let $L/K$ be number fields. Then $L$ is a splitting field for $K$ if for the polynomial $f\in K[x]$:
\begin{itemize}
\item $f$ splits over $L$.
\item There is no $M$, $L/M/K$, such that $M\neq L$ and $f$ splits over $M$.
\end{itemize}
\end{defn}
\begin{thm}\label{2.22}
A field extension $L/K$ is normal and finite if and only if $L$ is a splitting field for some polynomial $f\in K[x]$.
\end{thm}
We can also speak of seperable extensions
\begin{defn}\label{2.23}
Let $K$ be a number field. Then a polynomial $f\in K[x]$ is said to be separable if it splits into linear terms in $K$.
\end{defn}
\begin{defn}\label{2.24}
Let $L/K$ be number fields. Then $F$ is said to be a seperable extension of $K$ if for every $\alpha\in K$, $\irr_K(\alpha)$ is seperable over $L$.
\end{defn}
\begin{rmk}
A finite, normal, and seperable extension of a number field is also called a Galois extension.
\end{rmk}
For seperable extensions we have an equivalence allowing us to reduce to one irreducible polynomial.:
\begin{prop}\label{2.25}
For an extensions $L/K$ the following are equivalent:
\begin{enumerate}
\item $L$ is a Galois extension of $K$.
\item $L$ is the splitting field of a seperable polynomial $f\in K[x]$.
\end{enumerate}
\end{prop}
With this in place we will now start talking about automorphisms. First we define the Galois group:
\begin{defn}\label{2.26}
Let $L/K$ be number fields. Then the Galois group of $L$ over $K$ is defined as the group of $K$-automorphisms, i.e. the group of automorphisms $\pi:L\rightarrow L$ that is fixed for all elements in $K$. The group operation is composition of maps and the group is denoted $\text{Gal}(L/K).$
\end{defn}
The existance of such a group is trivial as the identity function on $L$ is a $K$-autmorphism, so will always lie in the Galois group. We make a distinction between real and complex $K$-automorphisms: A $K$-automorphism $\sigma\in\Gal(L/K)$ is real if $\sigma:L\rightarrow\reals$ and complex if $\sigma:L\rightarrow\compl$.
\begin{defn}\label{2.27}
Let $L/K$ be number fields with Galois group $\Gal(L/K)$ such that $\abs{\Gal(L/K)}=n$. Then the signature of $\Gal(L/K)$ is $\sign\left(\Gal(L/K)\right)=(r,s),$ where $r$ is the number of real $K$-automorphisms and $s$ is half the number of complex $K$-automorphisms in $\Gal(L/K)$ such that $n=r+2s$.
\end{defn}
There is one element of the Galois group that is particularly interesting, which is the Frobenius element. The idea of this comes from the Frobenius automorphism, which we recall for good measure.
\begin{defn}\label{2.28}
If $K$ is a field of characteristic $p>0$, the map $\phi:K\rightarrow K$ defined by $k\mapsto k^p$ is called the Frobenius monomorphism of $K$. When $K$ is finite, $\phi$ is called the Frobenius automorphism.
\end{defn}
To be able to define the Frobenius elements over extensions we have a little setting up to do. For this let $L/K$ be a Galois extension with Galois group $\text{Gal}(L/K)$. Let $\PriID\subseteq\mathcal{O}_L$ lying over $\priID\subseteq\mathcal{O}_K$. Then we can define the decomposition group $D(\PriID\mid\priID)$ as the set of automorphisms in $\text{Gal}(L/K)$ which fix $\PriID$:
$$D(\PriID\mid\priID)=\{\sigma\in\text{Gal}(L/K)\mid \sigma(\PriID)=\PriID\}.$$

\begin{defn}\label{2.29}
Let $L/K$ be a Galois extension with Galois group $\text{Gal}(L/K)$, let $\PriID\subseteq\mathcal{O}_L$ be a prime ideal unramified in $L/K$ lying over $\priID\in\mathcal{O}_K$. Then the Frobenius element, $\text{Frob}_\PriID$,  is the element of $D(\PriID\mid\priID)$ that acts as the Frobenius automorphism on the residue field. I.e. the unique $\text{Frob}_\PriID\in D(\PriID\mid\priID)$ such that
\begin{itemize}
\item $\text{Frob}_\PriID\in D(\PriID\mid\priID)$.
\item For all $\alpha\in\mathcal{O}_L$, $\text{Frob}_\PriID(\alpha)\equiv a^q\mod \PriID$, where $q=\abs{\mathcal{O}_K/\priID}$.
\end{itemize}
\end{defn}
This Frobenius element for $\PriID$ is uniquelly determined, but is not unrelated. For example, if $\PriID$ and $\PriID'$ both lie over $\priID\subset\mathcal{O}_K$ then the Frobenius elements are conjugate. This means that each prime ideal in $\mathcal{O}_K$ gives rise to a conjugacy class of Frobenius elements in $\text{Gal}(L/K)$.
\begin{prop}\label{2.30}
Let $L/K$ be a Galois extension with Galois group $\text{Gal}(L/K)$, let $\PriID,\PriID'\subseteq\mathcal{O}_L$ be prime ideals unramified in $L/K$, both lying over $\priID\in\mathcal{O}_K$. Then there exists $\sigma\in\text{Gal}(L/K)$ such that
$$\text{Frob}_\PriID=\sigma\text{Frob}_\PriID'\sigma^{-1}.$$
\end{prop}
\begin{proof}
Let $\alpha\in\mathcal{O}_L$. Then
$$\sigma\Frob_\PriID\sigma^{-1}(\alpha)=\sigma\left((\sigma^{-1}(\alpha))^q+a\right),$$
for some $a\in\PriID$, and
$$\sigma\left((\sigma^{-1}(\alpha))^q+a\right)=\alpha^q+\sigma(a)\equiv\alpha^q\mod \sigma\PriID.$$
\end{proof}
This is all we need for now, but in chapter \ref{ch:RNFS} we will be expand our theory a little more.
\subsection{Norms, discriminants, and ramification}\label{sec: PRE1norm}
\subsubsection{Norms and Traces}
Let $K=\quot(\alpha)$ be a number field with Galois group $\text{Gal}(K/\quot)$. Let $\mathcal{O}_K$ be the ring of algebraic integers of $K$. It will prove useful in what is to come that we can relate the elements of a number field to rational number in $\quot$. There are two natural maps to consider:
\begin{defn}\label{2.31}
Let $L/K$ be an extension of number fields with $\left[L:K\right]=n$. Let $\sigma_1,\ldots,\sigma_n$ be the $K$-automorphisms of $L$ in $\compl$. For any $\alpha\in L$ we define the following maps:
\begin{enumerate}
\item The trace of an element $\alpha$:
$$\text{Tr}_{L/K}(\alpha):L\rightarrow\compl,$$
$$\alpha\mapsto\sum_{i=1}^{n}\sigma_i(\alpha).$$
\item The field norm of an element $\alpha$:
$$N_{L/K}:L\rightarrow\compl,$$
$$\alpha\mapsto\prod_{i=1}^n\sigma_i(\alpha)^{\left[L:K(\alpha)\right]}.$$
\end{enumerate}
\end{defn}
\begin{rmk}
If it is clear from context what field we are taking the norm over we will often just write $N_{L/K}=N$.
\end{rmk}
We will mostly be concerned with Galois extensions which makes the field norm far simpler. For this consider $L/K$ with Galois group $\text{Gal}(L/K)$, then for any $\alpha\in L$ the norm becomes
$$N_{L/K}(\alpha)=\prod_{\sigma\in\text{Gal}(L/K)}\sigma(\alpha).$$
These two maps have a few important properties we will make use of throughout.
\begin{prop}\label{2.32}
\begin{enumerate}
\item []
\item $N_{L/K}$ and $\text{Tr}_{L/K}$ are homomorphisms from $L\rightarrow K$.
\item If $\alpha\in\mathcal{O}_L$ then $N_{L/K}(\alpha)\in\mathcal{O}_K$ and $\text{Tr}_{L/K}(\alpha)\in\mathcal{O}_K$.
\item Let $\alpha,\beta\in\mathcal{O}_L$ such that $\alpha\mid\beta$, then $N_{L/K}(\alpha)\mid N_{L/K}(\beta)$.
\item If $k\in K$, then $N_{L/K}(k)=k^n$ and $\text{Tr}_{L/K}(k)=nk$, where $\left[L:K\right]=n$.
\end{enumerate}
\end{prop}
\begin{rmk}
Note that (2) implies that for all $\alpha\in \mathcal{O}_L$: $\alpha\in\mathcal{O}_L^\ast \Leftrightarrow N_{L/\quot}(\alpha)=\pm1$ as $\mathcal{O}_\quot=\zahle$ and $\zahle^\ast=\pm 1$.
\end{rmk}
Despite the usefulness of the trace and field norm maps they are rather limited in that they are only defined for elements in the number field. We can also define the norm of an ideal.\footnote{ There are multiple ways to define the norm of an ideal. The definition we adopt has the benefit of being completely algebraic, which is why we use it.}
\begin{defn}\label{2.33}
Let $K$ be a number field of finite degree. Let $\mathcal{O}_K$ be the ring of algebraic integers of $K$. Let $\mathfrak{I}\subset\mathcal{O}_K$ be a non-trivial ideal of $\mathcal{O}_K$. Then
$$N(\mathfrak{I})=\abs{\mathcal{O}_K/\mathfrak{I}}.$$
\end{defn}
However these definitions are not completely independent from eachother, as they are related for principal ideals. In fact, for principal ideals the field and ideal norm overlap.
\begin{thm}\label{2.34}
Let $K$ be a number field of finite degree. Let $\mathfrak{I}=\langle\alpha\rangle\subset\mathcal{O}_K$. Then
$$\abs{\mathcal{O}_K/\mathfrak{I}}=N(\mathfrak{I})=N(\langle\alpha\rangle)=\abs{N(\alpha)}.$$
\end{thm}
It will not come as a surprise that the ideal norm is also a homomorphism, such that for every two non-zero integral ideals $\mathfrak{A,B}\in\mathcal{O}_K$, $N(\mathfrak{AB})=N(\mathfrak{A})N(\mathfrak{B})$. This is absolutely not trivial, but intuitive enough for us to understand to leave the proof to [\cite{IANT}, theorem 9.3.2]. Now we have a definition for the ideal norm and an explicit, and reasonably simple, way to compute this for principal ideals. Up until now we do not have a general way to compute the ideal norm, but we can get closer by recalling theorem \ref{2.16} as the ideal norm is a homomorphism, so we just have to consider how to compute the ideal norm for prime ideals. 
\begin{thm}\label{2.35}
Let $K$ be an algebraic number field. Let $\PriID\subset\mathcal{O}_K$. Then there exists a unique prime $p\in\zahle$, such that
$$\PriID\mid\langle p\rangle.$$
\end{thm}
As this $p\in\zahle$ is unique we say $\PriID$ lies over $p$, or equivalently $p$ lies below $\PriID$. This has an immediate effect on the ideal norm, as can be seen in the following theorem:
\begin{thm}\label{2.36}
Let $K$ be an algebraic number field with $\left[K:\quot\right]=n$. Let $\PriID\subset\mathcal{O}_K$ be a prime ideal lying over $p\in\zahle$. Then
$$N(\PriID)=p^f,$$
for some integer $f\in\{1,\ldots, n\}$.
\end{thm}
\begin{proof}
As $\PriID$ lies above $p$ we have that $P\mid\langle p\rangle$. Hence $\langle p\rangle =\PriID\mathfrak{Q}$ for some integral ideal $\mathfrak{Q}\subset\mathcal{O}_K$. As the norm is multiplicative we have that
$$N(\langle p\rangle)=N(\PriID)N(\mathfrak{Q}).$$
As the $K$-conjugates of $p$ comprise $p$ $n$ times we have
$$N(p)=p^n,$$
and so
$$N(\langle p\rangle)=\abs{N(p)}=N(\PriID)N(\mathfrak{Q})=p^n.$$
\end{proof}
The $f$ in the theorem is called the inertial degree and is used to define one of the most interesting objects in number theory, the ideal class group, which will define momentarily. First we need some preliminary definitions and properties.
 \begin{defn}\label{2.37}
Let $K$ be an algebraic number field with $\left[K:\quot\right]=n$. Let $p\in\zahle$ be prime such that
$$\langle p\rangle=\prod_{i=1}^{g}\PriID_i^{e_i}.$$
Then $g$ is called the decomposition number, $e_i$ is the ramification index for $\PriID$, and the following properties hold:
\begin{enumerate}
\item The inertia degree, $f_i$, can be computed as $f=\left[\mathcal{O}_K/\PriID_i:\zahle/\langle p\rangle\right]$.
\item $p$ is said to be ramified in $K$ if, for some $i$, $e_i>1$.
\item $p$ is said to be unramified in $K$ if for all $i$: $e_i=1$.
\item The following equation holds:
$$\sum_{i=1}^ge_if_i=n.$$
\item $P$ is said to be (totally) split if for all $i$: $e_i=f_i=1$, equivalently $g=n$.
\end{enumerate}
\end{defn}
\begin{rmk}
This can also be defined generally for extensions $L/K$, and follows naturally from the above definition. For our purposes the above will suffice however.
\end{rmk}
\subsubsection{Discriminants and Minkowski's bound}
We now define the concept of the (absolute) discriminant of a number field and immediately apply it to a theorem by Dedekind:
\begin{defn}\label{2.38}
Let $K$ be a number field of finite degree with ring of algebraic integers $\mathcal{O}_K$. Let $\{a_1,\ldots,a_n\}$ be an integral basis for $\mathcal{O}_K$ and let $\{\sigma_1,\ldots,\sigma_n\}$ be the set of $K$-automorphisms into $\compl$. Then the (absolute) discriminant of $K$ is
$$\Delta_K=\det\begin{pmatrix}
\sigma_1(a_1)&\ldots&\sigma_1(a_n)\\
\vdots & \ddots & \vdots\\
\sigma_n(a_1) &\ldots & \sigma_n(a_n)
\end{pmatrix}^2.$$
\end{defn}
\begin{thm}\label{2.39}
Let $K$ be an algebraic number field. Then the prime $p\in\zahle$ ramifies in $K$ if and only if $p\mid\Delta_K$.
\end{thm}
Now we define the ideal class group:
\begin{defn}\label{2.40}
Let $K$ be an algebraic number field and let $\mathfrak{A},\mathfrak{B}\subset\mathcal{O}_K$ be non-zero ideals. Then we say that $\mathfrak{A}\sim\mathfrak{B}$, that is they are equivalent, if there exists $\alpha,\beta\in\mathcal{O}_K\nozero$ such that $\langle b\rangle\mathfrak{A}=\langle a\rangle\mathfrak{B}$. The equivalence classes under $\sim$ are called ideal classes and the set of ideal classes of $K$, denoted $\text{Cl}_K$, is called the ideal class group. Moreover, the cardinality of the ideal class group, $\abs{Cl_K}=h_K$, is called the class number of $K$.
\end{defn}
It is well known that $h_K$ is finite and that $\text{Cl}_K$ is an abelian group for $K$ an algebraic number field. More can be said, especially about the generation of $\text{Cl}_K$, and we will record these in quick succession.
\begin{thm}\label{2.41}
Let $\mathfrak{I}\subseteq\mathcal{O}_K$ be a non-zero ideal. Then there exists a non-zero $\alpha\in\mathfrak{I}$ with
$$\abs{N_{K/\quot}(\alpha)}\leq c_KN(\mathfrak{I}),$$
for
$$c_K=\left(\frac{4}{\pi}\right)^{r_2}\frac{n!}{n^n}\sqrt{\abs{\Delta_K}},$$
where $2r_2=\abs{\{\sigma_i\mid\sigma_i \text{ a }\quot\text{-automorphism}, \sigma_i(K)\not\subseteq\reals\}}$.
\end{thm}
\begin{rmk}
$c_K$ is called the Minkowski bound.
\end{rmk}
\begin{thm}\label{2.42}
Any ideal class $c\in\text{Cl}_K$ contains an ideal $\mathfrak{I}$ such that $N(\mathfrak{I})\leq c_K.$
\end{thm}
\begin{thm}\label{2.43}
Let $K$ be a number field. Then $\text{Cl}_K$ is generated by the classes of the prime ideals $\left[\PriID\right]$, satisfying $N(\PriID)\leq c_K$. In particular, any such $\PriID$ must lie above a prime $p\leq c_K$.
\end{thm}
\subsubsection{Relativity and Transitivity}
In the last section we defined the absolute discriminant. It will not surprise the reader that there then also is a relative discriminant. To do this we need to consider the following:
\begin{defn}\label{2.44}
Let $L/K$ be number field with rings of integers $\mathcal{O}_L$ and $\mathcal{O}_K$ respectively. The fractional ideal ,
$$\mathfrak{C}_{\mathcal{O}_L\mid\mathcal{O}_K}=\{x\in L\mid\text{Tr}(x\mathcal{O}_L)\subset\mathcal{O}_K\},$$
is called Dedekind's complementary module. It's inverse,
$$\mathfrak{D}_{\mathcal{O}_L\mid\mathcal{O}_K}=\mathfrak{C}_{\mathcal{O}_L\mid\mathcal{O}_K}^{-1},$$
is called the different of $\mathcal{O}_L/\mathcal{O}_K$.
\end{defn}
\begin{rmk}
With the usual abuse of notation we will write $\mathfrak{D}_{L/K}$ for the different $\mathfrak{D}_{\mathcal{O}_L\mid\mathcal{O}_K}$
\end{rmk}
By definition Dedekind's complementary module contains $\mathcal{O}_L$ and so $\mathfrak{D}_{L/K}$ is in fact an integral ideal. The different is well-behaved in towers of number fields:
\begin{prop}\label{2.45}
For a tower of fields $L/M/K/\quot$: $\mathfrak{D}_{L/K}=\mathfrak{D}_{L/M}\mathfrak{D}_{M/K}$.
\end{prop}
This different ideal, or simply `different', is key in the definition of the relative discriminant;
\begin{defn}\label{2.46}
Let $L/K$ be number fields. The relative discriminant $\Delta_{L/K}$ is the ideal of $\mathcal{O}_K$ defined by the relative norm of the different:
$$\Delta_{L/K}=N_{L/K}(\mathfrak{D}_{L/K}).$$
\end{defn}
This allows us to define a transitivity condition for the discriminant
\begin{lemma}\label{2.47}
Let $L/M/K$ be a tower of number fields. Then
$$\Delta_{L/K}=\Delta_{M/K}^{\left[L:M\right]}N_{M/K}(\Delta_{L/M}).$$
\begin{proof}
Applying to $\mathfrak{D}_{L/K}=\mathfrak{D}_{L/M}\mathfrak{D}_{M/K}$ the norm $N_{L/K}=N_{M/K}\circ N_{L/M}$ we obtain
$$\Delta_{L/K}=N_{M/K}(\Delta_{L/M})N_{M/K}(\mathfrak{D}_{M/K}^{\left[L:M\right]})=N_{M/K}(\Delta_{L/M})\Delta_{M/K}^{\left[L:M\right]}.$$
\end{proof}
\end{lemma}
\subsection{Dirichlet's Unit Theorem}\label{sec: PRE1unit}
In this last preliminary section we want to get an idea of how $\mathcal{O}_K$ looks. Specifically, we want to get an idea of the units in $\mathcal{O}_K$. Note that the units of a ring form a group called the unit group and is denoted $\mathcal{O}_K^\times$.
\begin{defn}\label{2.48}
Let $n\in\nat$ and $\zeta_n=e^{2\pi i/n}$. Then $\zeta_n$ is a $n$th root of unity and the number field $\quot(\zeta_n)$ is called the $n$th cyclotomic field.
\end{defn}
These cyclotomic fields are in many ways the easiest number fields to work with because they behave especially nicely when $n$ is an odd prime. To emphasize this let $p\in\zahle$ be an odd prime for the remainder of the section.
\begin{thm}\label{2.49}
Let $K=\quot(\zeta_p)$, then $\mathcal{O}_K=\zahle[\zeta_p]$.
\end{thm}
With that little note we return to what is at hand. Let $\mu(K)$ be the group of units contained in a number field $K$. Note that $\mu(K)$ is a cyclic group of order $n$ under multiplication. The following proposition is then a tautology:
\begin{prop}\label{2.50}
Let $K$ be a number field. Then $\mu(K)\subset\mathcal{O}_K^\times$.
\end{prop}
The question now becomes: What are the other elements in $\mathcal{O}_K^\times$? Recall from \ref{sec: PRE1galt} that the signature of a field, $\text{sign}(K)=(r,s)$, where $r$ is the number of real $\quot$-embeddings and $2s$ the number of complex $\quot$-embeddings.
\begin{thm}\label{2.51}
Let $K$ be a number field. There exists a map $l:K\rightarrow \reals^{r+s-1}$ such that $\ker(l)$ is finite and cyclic, and $\ker(l)\cong\mu(K)$. Moreover
$$l(\mathcal{O}_K^\times)\cong \zahle^{r+s-1}.$$
\end{thm}
This immediately leads us to the concluding theorem:
\begin{thm} \textbf{Dirichlet's Unit Theorem} \label{DUT}
Let $K$ be a number field and let $\text{sign}(K)=(r,s)$. Then
$$\mathcal{O}_K^\times\cong\mu(K)\times\zahle^{r+s-1}.$$
\end{thm}
\begin{proof}
From theorem \ref{2.51} we have that there exists a map $l$ such that $\ker(l)\cong\mu(K)$ hence by the isomorphism theorems
$$K/\mu(K)=K/\ker(l)\cong \text{im}(l)=\zahle^{r+s-1},$$
and as $\mu(K)$ is an abelian group and $\zahle^{r+s-1}$ is free we get that $K\cong \mu(K)\times\zahle^{r+s-1}$.
\end{proof}
\newpage
\section{The General Number Field Sieve}\label{ch: GNFS}
For the entirety of this section let $n\in\nat$ be a given integer that we wish to factor. Without loss of generality assume that this $n$ is odd, else we could simply consider $m\in\nat$ where $n=2^sm$, and that $n$ is composite, else a primality test, which can be done in polynomial time, such as the AKS Primality Test \cite{GRAN}, will show that $n$ does not need factoring. The General Number Field Sieve (GNFS) then uses a congruence of squares modulo $n$ to perform integer factorization. Assume that $(x,y)$ is a pair that admits a congruence of squares modulo $n$ so that
$$x^2\equiv y^2 \mod n,$$
where $x\neq \pm y$. Then we can use $(x,y)$ to find a non-trivial factorization by using the well-established fact that $x^2-y^2=(x+y)(x-y)$. We then compute $\gcd(x+y,n)=d_1$ and $\gcd(x-y,n)=d_2$, in the hope that this gives a non-trivial factorization, that is $d_1,d_2\neq 1$ or $n$. Assume that we are able to produce random integers $x$ and $y$ such that the above congruence holds. Then for the hardest case, $n=pq$ semiprime, Table 1 shows an exhaustive truth table:

\begin{table}[h!]\label{tab:1}
\centering
\caption{Factorization cases for $n=pq$}
\begin{tabular}{|c|c|c|c|c|c|c|}
\hline
$p\mid x+y$ & $p\mid x-y$ & $q\mid x+y$ & $q\mid x-y$ & $\gcd(x+y,n)$ & $\gcd(x-y,n)$ & Factor?\\
\hline
\hline
T & T & T & T &$ n$ &$ n$ & F\\
T & T & T & F &$ n $&$ p$ & T\\
T & T & F & T &$ p $&$ n $& T\\
T & F & T & T &$ n $&$ q $& T\\
T & F & T & F &$ n $&$ 1 $& F\\
T & F & F & T &$ p $&$ q $& T\\
F & T & T & T &$ q $&$n $& T\\
F & T & T & F &$ q $&$ p $& T\\
F & T & F & T &$ 1 $&$ n $& F\\
\hline
\end{tabular}
\end{table}

It can be seen that we only have failure if $p$ and $q$ have the same conditions on the factorization of $x+y$ and $x-y$. This is because that means that $x\equiv \pm y\mod n$. If that is not the case then we always find a non-trivial factorization. How the NFS operates to find such congruent pairs $(x,y)$ is the goal of this section.

\subsection{Introducing the algorithm}\label{sec:GNFSintro}
The Number Field Sieve is so named because of a sieving process happening both over $\zahle$ and a ring $\zahle[\alpha]$ which are contained in the field $\quot$ and the number field $\quot(\alpha)$ respectively. In this section we will give a description on how this is done, while we go into detail, introducing formal definitions and rigor, in Section \ref{sec:GNFSobstr}. However, the following commuting diagram captures the whole process: 
$$\begin{tikzcd}[row sep=scriptsize, column sep=tiny]
&\arrow [dl, "x \mapsto m" left] \zahle[x] \arrow[dr, "\mod f"]& \\
\zahle \arrow[dr, "\mod n" left] & & \zahle[x]/(f)\cong\zahle[\alpha] \arrow[dl, "x\mapsto m \mod n"]\\
& \zahle/n\zahle &
\end{tikzcd}$$

Before we can discuss how the NFS obtains a congruence we introduce the concept of smooth numbers. 
\begin{defn}\label{3.1}
Let $z\in\zahle$ and let $B\in\nat$. Then $z$ is said to be $B$-smooth if in the prime factorization of $z$,
$$z=\pm\prod_{i=1}^n p_i^{e_i},$$
it holds that $p_i\leq B$ for all $i=1,\ldots, n$.
\end{defn}
It is clear that this definition is insufficient for $\mathcal{O}_{\quot(\alpha)}$ as this ring does not solely consist of integers, however that is easily amended.\\
\\
For the sake of exposition we will consider $\mathcal{O}_{\quot(\alpha)}=\zahle[\alpha]$. This is rarely the case for number fields in general as this only happens for monogenic fields, but in section \ref{sec:GNFSobstr} we will see that this assumption can be dealt with algebraically so we only need to concern ourselves with this case. Let $a-b\alpha\in\zahle[\alpha]$ and let $\langle a-b\alpha\rangle$ be the ideal generated by it. As $\quot(\alpha)$ is an algebraic number field we get by that $\zahle[\alpha]$ is a Dedekind domain, hence any ideal splits uniquely into a product of prime ideals:
$$\langle a-b\alpha \rangle = \prod_{i=1}^n \PriID_i^{\eta_i},\hspace{5px}$$
Moreover by theorem \ref{2.36} we know that any prime ideal of $\zahle[\alpha]$ evaluates to a  prime power under the norm map. This allows us to define an analogy to smoothness over $\zahle$:
\begin{defn}\label{3.2}
Let $a-b\alpha \in\zahle[\alpha]$ and let $B\in\nat$. Consider $\langle a-b\alpha \rangle\subseteq \zahle[\alpha]$ with factorization into prime ideals
$$\langle a-b\alpha \rangle = \prod_{i=1}^n \PriID_i^{\eta_i}.$$
Then $a-b\alpha $ is said to be $B$-smooth if for all $N(\PriID_i)=p_i^{f_i}$ it holds that
$$p_i\leq B,$$
where $N:\quot(\alpha)\rightarrow\quot$ is the field norm.
\end{defn}
Using these definitions we can, for some fixed $m$ to be defined later, define sets $S_\zahle$ and $S_{\zahle[\alpha]}$ with smoothness bounds $B$ and $B'$ respectively, as follows:
$$S_\zahle=\big\{(a,b)\in\zahle^2\mid \gcd(a,b)=1,a-bm\text{ is }B\text{-smooth as in Definition \ref{3.1}}\big\}$$
$$S_{\zahle[\alpha]}=\big\{(a,b)\in\zahle^2\mid\gcd(a,b)=1,a-b\alpha\text{ is }B'\text{-smooth as in Definition \ref{3.2}} \big\}$$
The bulk of the work done in the GNFS is done in determining these sets and we will go in much more detail in the next section. If we succeed in finding a set $S\subseteq S_\zahle\cap S_{\zahle[\alpha]}$, such that $\abs{S}\geq \pi(B)+\pi(B')+1$, then the GNFS uses the remainder of its computation time to find a set $T\subset S$ such that 
\begin{align}\label{eq1}
\prod_{(a,b)\in T} (a-bm) \text{ is a square in }\zahle,
\end{align}
\begin{align}\label{eq2}
\prod_{(a,b)\in T} (a-b\alpha) \text{ is a square in }\zahle[\alpha].
\end{align}
It is not a certainty that such a set $T$ fulfilling these equations can be found. If it is not then the GNFS can be restarted with larger smoothness bounds $B$ and $B'$. For now assume we have found such a set $T$. From here completion of the algorithm is nearly trivial:\\
\\
First we observe that $f'(m)^2\prod_{(a,b)\in T} (a-bm)$ must give a square in $\zahle$, $z^2$ say, and $f'(\alpha)^2\prod_{(a,b)\in T}(a-b\alpha)$ must give a square in $\zahle[\alpha]$, $\xi^2$ say, for which there exist elementary methods to find the square roots $z$ and $\xi$.
\begin{rmk}
We have multiplied equation \eqref{eq1} and \eqref{eq2} by the square of the derivative of a polynomial at $m$ and $\alpha$ respectively. This is to ensure that $z\in\zahle$ and $\xi\in\zahle[\alpha]$. A more thorough discussion of this follows in \ref{sec:GNFSobstr}.
\end{rmk}
From here we compute $\gcd(z-N(\xi),n)$ and $\gcd(z+N(\xi),n)$. If this is non-trivial, then the result is a prime factor of $n$. Else we have to conclude failure and start again with different parameters.
\\
This describes the algorithm in grand lines. Now it is time some mathematical rigor is introduced.
\subsection{The algorithm Explained}\label{sec:GNFSexpl}
\subsubsection{Choosing a polynomial}
The first step of the NFS is that we must choose a $d\in\nat, d\neq 1$ and define $m=\lfloor n^{\frac{1}{d}}\rfloor$ s.t. $\gcd(m,n)=1$. This choice of $m$ gives rise to a polynomial $f\in\zahle[x]$ such that $n\mid f(m)$. Over the years there have been many attempts at making the choice of polynomial as effective as possible, but the simplest and effectively cheapest\footnote{When we say ``cheapest" we mean computationally most effective. } way to do this is to use the ``base-$m$" method. For each integer $m\in\zahle$ we can write $n$ as a linear combinations of powers of $m$ such that
$$n=\sum_{i=0}^d c_im^i.$$
Then we can define a polynomial $f\in\zahle[x]$ such that
$$f(x)=\sum_{i=0}^d c_ix^i,$$
where the $c_i\in\{0,\ldots,m-1\}$ are the same as the base-$m$ expansion of $n$. This guarantees that $n\mid f(m)$, in fact $f(m)=n$, and $\deg(f)=d$. By [\cite{BLP}, prop. 3.2] we have that the leading coefficient $c_d=1$, such that $f\in\zahle[x]$ is in fact monic. 
\begin{rmk}\label{rmk1}
Note that as $\abs{c_i}\in\{0,\ldots,m-1\}$ hence $c_i<m<n^{\frac{1}{d}}$ we can see that the discriminant of $f$ satisfies that
$$\abs{\Delta(f)}<d^{2d}n^{2-\frac{3}{d}}.$$
\end{rmk}
We may also assume that $f\in\zahle[x]$ is irreducible. To see this assume, to the contrary, that $f\in\zahle[x]$ is reducible. This means that there exist non-trivial polynomials $g,h\in\zahle[x]$ such that 
$$f(x)=g(x)h(x).$$
So a simple computation gives a non-trivial factorization, as $g(x)$ and $h(x)$ are assumed to be non-trivial. As we may now assume $f(x)$ to be monic and irreducible we have shown the following:
\begin{prop}\label{3.3}
Let $f(x)$ be a monovariate, monic, irreducible polynomial, with $\deg(f)=d$, in $\zahle[x]$ with roots $\alpha_1,\ldots,\alpha_d$, not necessarily in $\zahle$. Then for all $i=1,\ldots,d$,
$$f(x)=\irr_\zahle(\alpha_i).$$
\end{prop}
Note that by assumption $\deg(f)=d>1$. As $f$ is an irreducible polynomial over $\zahle[x]$ we have that none of the roots of $f$ lie in $\zahle$, for else $f$ would split into linear factors over $\zahle$. This root, $\alpha$ say, lies in $\compl$ and as $f(\alpha)=0$ this means that $\alpha$ is an algebraic integer. Moreover, by theorem \ref{2.6} we have that $\quot(\alpha)$ is a field and has basis $\big\{1,\alpha,\ldots,\alpha^{d-1}\big\}$ over $\quot$.
\subsubsection{Sieving over $\zahle$}
Now that we have defined a polynomial $f\in\zahle[x]$ and we have defined a field extension $\quot(\alpha)$ we can start sieving for pairs $(a,b)$ such that equation (\ref{eq1}) and (\ref{eq2}) hold. For the sake of exposition we will focus on the mathematical aspects and therefore consider the pairs $(a,b)$ over $\zahle$ and $\zahle[\alpha]$ separately. 
\begin{defn}\label{3.4}
Let $u\in\zahle$, then the set of integer pairs defined by 
$$U=\big\{(a,b)\in\zahle\mid \abs{a}\leq u, 0<b<u,\gcd(a,b)=1\big\}$$
is called the universe of sieving.
\end{defn}
Now let $u$ be a large integer dependent on the to-factor $n$ and define
$$\mathcal{B}=\big\{p\in\zahle\mid p\text{ prime}, p\leq B\big\} \cup \big\{\pm1\big\}$$
as the rational factor base of the sieve\footnote{The rational factor base, normally, only contains the primes $p\leq B$. In our case, however, appending $\pm 1$ improves the effectivity as we can now deal with $a-bm<0$.}, fully dependent on the chosen smoothness bound $B$. Then there is a standard procedure to work through to fill up the set $S_\zahle$:
\begin{algorithm}[H]\label{Alg1}
\caption{Procedure to populate $S_\zahle$}
\begin{algorithmic}
\Input{Universe $U$ of $(a,b)$ pairs, smoothness bound $B$}
\Output{Set $S_\zahle=\{\left(a,b\right)\in U\mid a-bm \text{ is }B\text{-smooth},\gcd(a,b)=1\}$}
\State{Initialize array comprised of $a-mb$ for all $(a,b)\in U$}
\For{each element in the array}
\State compute the set, $P_{(a,b)}$, of primes $p$, $p\leq B$, such that $a\equiv bm \mod p$
\For{each $p\in P_{(a,b)}$}
divide $a-bm$ by the maximal power of $p$, such that $p$ does not divide the quotient, and replace $a-bm$ by this quotient
\EndFor
\If{the the element in the array is $\pm 1$ \textbf{and} $\gcd(a,b)=1$}
\State add $(a,b)$ to $S_\zahle$
\EndIf
\EndFor
\end{algorithmic}
\end{algorithm}
We will accept that this algorithm terminates by choosing $B$ large enough and returns the set $S_\zahle$ such that $\abs{S_\zahle}\geq\pi(B)+1$. As we have recorded the vectors $\{e_p(a-bm)\}_{p\leq B}$ of $p$-exponents for all $p\leq B$ prime for each $a-bm$ we now define $M\in\mathcal{M}_{\abs{S_\zahle}\times \pi(B)}(\zahle)$ given by: $\forall (a_i,b_i)\in S_\zahle, \forall p_j\in B:$
$$M=\begin{pmatrix}
\text{sign}(a_1-b_1m) & e_{p_1}(a_1-b_1m) & \cdots & e_{p_{\pi(B)}}(a_1-b_1m)\\
\cdots & \cdots& \cdots & \cdots\\
\text{sign}(a_{\abs{S_\zahle}}-b_{\abs{S_\zahle}}m) & e_{p_1}(a_{\abs{S_\zahle}}-b_{\abs{S_\zahle}}m) & \cdots & e_{p_{\pi(B)}}(a_{\abs{S_\zahle}}-b_{\abs{S_\zahle}}m)
\end{pmatrix},$$
where the sign-bit is defined as: 
$$\text{sign}(a_i-b_im)=\begin{cases}
1 & a_i-b_im<0 \\
0 & \text{otherwise}
\end{cases}.$$
To achieve the situation in which equation (\ref{eq1}) holds we now want to find an independent subset of $S_\zahle$. For this we first consider a pair $a-bm\in S_\zahle$, then it is easily observed that for this to be a square in $\zahle$ we must have that
$$ x^2 = a-bm = \prod_{p\in\mathcal{B}}p^{e(p)}=\bigg(\prod_{p\in\mathcal{B}}p^{\frac{e(p)}{2}}\bigg)^2.$$
Hence 
$$x=\pm\prod_{p\in\mathcal{B}}p^{e(p)}.$$
So it suffices for us to look for the independent subset of $S_\zahle$ modulo $\mathbb{F}_2$, the finite field of characteristic 2. Hence we define 
$$M_2= M \mod 2 = (m_{i,j} \mod 2)_{i\in \{1,\ldots,\abs{S_\zahle}\}, j\in \{1,\ldots,\pi(B)\}}.$$
Let $f(a_i+b_im)$ be the $i$th row of $M_2$. As $\abs{S_\zahle}>\pi(B)$ we have that there is a linearly dependent subset in $M_2$ and hence there is a non-trivial solution $T_\zahle\subset S_\zahle$ to the linear equation:
$$\sum_{(a,b)\in T_\zahle} f(a-bm) = 0.$$
With this set we obtain (\ref{eq1}):
$$\prod_{(a,b)\in T_\zahle} (a-bm) \text{ is a square in }\zahle.$$
\subsubsection{Sieving over $\zahle[\alpha]$}
To sieve over $\zahle[\alpha]$ we will attempt to have a similar construction as the process used for $\zahle$, but to do this we will first have to deal with a few of the obstructions that arise. First recall from theorem \ref{2.36}
that for a prime ideal $\mathfrak{P}$ we have that
$$N(\mathfrak{P})=p^f,$$
where $p\in\zahle$ is prime and $f\in\nat$. Then $p$ is called the prime lying below $\mathfrak{P}$ and $f$ is called the inertia degree of $\mathfrak{P}$. Moreover a prime ideal for which $f=1$ is called a first degree prime. For a first degree prime ideal, $\mathfrak{P}$, we have that the index $[\zahle[\alpha]:\mathfrak{P}]=p$, hence gives rise to the isomorphism
$$\zahle[\alpha]/\mathfrak{P}\cong\zahle/p\zahle,$$
hence $\zahle[\alpha]/\mathfrak{P}$ is a field. This gives us a direct link between $\zahle[\alpha]$ and $\zahle/p\zahle$:
\begin{thm}\label{3.5}
Let $f(x)$ be a monic irreducible polynomial with coefficients in $\zahle$. Let $\alpha$ be a root of $f(x)$, then there is a bijective correspondence between the set $\mathcal{P}_1$ of first degree prime ideals and the set $$\{(p,m):p \text{ prime}, m\in\zahle/p\zahle, f(m)\equiv 0 \mod p\}.$$
\end{thm}
\begin{proof}
Let $\mathfrak{p}$ be a first degree prime ideal of $\zahle[\alpha]$. then $\left[\zahle[\alpha]:\mathfrak{p}\right]=p$ for some prime integer $p$ so that $\zahle[\alpha]/\mathfrak{p}\cong\zahle/p\zahle$. There is a canonical ring epimorphism $\theta:\zahle[\alpha]\rightarrow\zahle[\alpha]/\mathfrak{p}$ such that $\ker(\theta)=\mathfrak{p}$, hence for $z\in\mathfrak{p}: p\mid\theta(z)$, moreover for $n\in\zahle$ such that $p\mid n$ there is a $z\in\mathfrak{p}$ so $\theta(z)=n$. As $\theta$ is a homomorphism $\theta(1)=1$ and so $\forall n\in\zahle:\theta(n)\equiv n\mod p$.\\
\\
Now let $m=\theta(\alpha)\in\zahle/p\zahle$. If $f(x)=\sum_{i=0}^d a_ix^i$ with $a_d=1$ and $a_i\in\zahle$ then $\theta(f(\alpha))\equiv 0 \mod p$ as $f(\alpha)=0$, and so
$$0\equiv\theta(f(\alpha))\equiv\theta\left(\sum_{i=1}^d a_ix^i\right)\equiv \sum_{i=1}^d a_i\theta(x)^i\equiv \sum_{i=1}^d a_ik^i\equiv f(m)\mod p,$$
hence $f(m)\equiv 0 \mod p$ and $\mathfrak{p}$ determines the unique pair $(p,m)$.\\
\\
Conversely, let $p$ be a prime integer and $m\in\zahle/p\zahle$ with $f(r)\equiv 0\mod p$. Then there is a natural ring epimorphism that maps polynomials in $\alpha$ to polynomials in $r$. In particular $\theta(a)\equiv a \mod p$ for all $a\in\zahle$ and $\theta(\alpha)\equiv m\mod p$. Let $\mathfrak{p}=\ker(\theta)$ so that $\mathfrak{p}$ is an ideal of $\zahle[\alpha]$. Since $\theta$ is surjective that means that $\zahle[\alpha]/\mathfrak{p}\cong\zahle/p\zahle$ and so $\left[\zahle[\alpha]:\mathfrak{p}\right]=p$. This implication is unique, hence we have the two unique implications,
$$(p,m)\hookrightarrow\mathfrak{p}\hookrightarrow(p,m),$$
proving the theorem.
\end{proof}
Hence finding these first degree prime ideals is equivalent to finding roots $\mod p$ for the minimal polynomial $\irr_K(\alpha)$. Finding roots of polynomials over finite fields has a well documented background \cite{CANT}, for example Berlekamp's algorithm or Cantor-Zassenhaus. We may therefore assume that finding first degree prime ideals is easy\footnote{And with easy we mean effective, and in that not adding to our overall complexity, as the size of $\deg(f)$ is far smaller than the size of $n$.} and therefore assume we can find sufficient first degree prime ideals.\\
\\
Now that we have restricted our prime ideals we may generalize the smoothness test for $\zahle[\alpha]$. Buhler et. al suggested the following algorithm in \cite{BLP}: For a smoothness bound $B'\in\zahle$, for now further undefined, $S_{\zahle[\alpha]}$ as defined above, and the polynomial $f$ with root $\alpha=t$.
\begin{algorithm}[H]\label{Alg2}
\caption{Procedure to populate $S_{\zahle[\alpha]}$}
\begin{algorithmic}
\Input{Universe $U$, polynomial $f$, smoothness bound $B'$}
\Output{Set $S_{\zahle[\alpha]}=\{(a,b)\mid a-b\alpha \text{ is }B'\text{-smooth}, \gcd(a,b)=1\}$}
\For{each prime $p\leq B'$}
\For{each $(a,b)\in U$ such that $b\not\equiv 0\mod p$}
\State initialize an array populated by $N(a-b\alpha)$
\EndFor
\For{each $r\in\zahle/p\zahle$}
\State compute the set $R(p)=\{r\mid f(r)\equiv 0 \mod p\}$
\EndFor
\For{each $r\in R(p)$}
\If{$a\equiv br\mod p$}
\State retrieve $N(a-b\alpha)$
\State divide $N(a-b\alpha)$ by the maximal power of $p$,  such that $p$ does not divide the quotient, and replace $N(a-b\alpha)$ by this quotient.
\EndIf
\EndFor
\EndFor
\For{ each element in the array corresponding to the pair $(a,b)$}
\If{ the element is $\pm1$}
\State add the pair $(a,b)$ to $S_{\zahle[\alpha]}$
\EndIf
\EndFor
\end{algorithmic}
\end{algorithm}
\begin{rmk}
It is clear that if $b\equiv 0 \mod p$ then there are no integers with $(a,b)\in U$ and $N(\langle a-b\alpha\rangle)=N(a-b\alpha)\equiv 0 \mod p$. 
\end{rmk}
What is left to us is to find the square in $\zahle[\alpha]$ that we wish to use. However following the same procedure as the rational sieve would leave us, not with \eqref{eq2}, but with 
$$N\left(\prod_{(a,b)\in S_{\zahle[\alpha]}} (a-b\alpha)\right)\text{ is a square in }\zahle.$$
This is clearly a necessary condition for (\ref{eq2}), but it is not sufficient. To combat this obstruction we recall the pairs $(p,R(p))$ as defined in algorithm \ref{Alg2} and define the following:
\begin{prop}\label{3.6}
Let $a,b\in\zahle$, $\gcd(a,b)=1$, and $p\in\zahle$ prime. Let $r\in R(p)$. Then the function
$$e_{(p,r)}(a-b\alpha)=\begin{cases}
\ord_p(N(a-b\alpha)) & a-br\equiv 0\mod p\\
0 & otherwise
\end{cases}.$$
, where for $z\in\zahle: \ord_p(z)=f$ if $p^f\mid\mid z$, is well defined. Moreover
$$N\left(a-b\alpha\right)=\prod_{(p,r)}\left(p^{e_{(p,r)}(a-b\alpha)}\right).$$
\end{prop}
To see how this can be used to produce the square in $\zahle[\alpha]$ we need the last theorem of this section
\begin{thm}\label{3.7}
Let $S'$ be a finite set of coprime integer pairs $(a,b)$ fulfilling equation \eqref{eq2}. Then for each prime number $p$ and each $r\in R(p)$ we have
$$\sum_{(a,b)\in {S_\zahle[\alpha]}} e_{(p,r)}(a-b\alpha)\equiv 0\mod 2$$
\end{thm}
\begin{proof}
For each prime ideal $\mathfrak{P}\in\zahle[\alpha]$ define the group homomorphism $\varphi_\mathfrak{P}:\quot[\alpha]^\ast\rightarrow\zahle$ such that
\begin{enumerate}
\item $\varphi_\mathfrak{P}(\beta)\geq 0$ for all $\beta\in\zahle[\alpha],\beta\neq 0$
\item if $\beta\in\zahle[\alpha], \beta\neq 0$, then $\varphi_\mathfrak{P}(\beta)>0$ if and only if $\beta\in \mathfrak{P}$
\item for each $\beta\in\quot[\alpha]^\ast$ one has $\varphi_\mathfrak{P}(\beta)=0$ for all but finitely many $\mathfrak{P}$, and 
$$\abs{N(\beta)}=\prod_\mathfrak{P}N(\mathfrak{P})^{\varphi_\mathfrak{P}(\beta)}.$$
\end{enumerate}
Then $\varphi_\mathfrak{P}$ is a $\mathfrak{P}$-adic valuation from a multiplicative group of units to an additive group of integers , hence, for $\gcd(a,b)=1$, if $\mathfrak{P}$ is not a first degree prime then $\varphi_\mathfrak{P}(a-b\alpha)=0$ and if $\mathfrak{P}$ corresponds to the pair $(p,m)$ as defined in theorem \ref{3.6} then $\varphi_\mathfrak{P}(a-b\alpha)=e_{p,r}(a-b\alpha)$.\\
\\
Now let $\mathfrak{P}_i\in\zahle[\alpha]$ be a first degree prime ideal and let 
$$\prod_{(a,b)\in {S_\zahle[\alpha]}} (a-b\alpha)=\xi^2.$$
Then, as $\varphi_\mathfrak{P}$ is a homomorphism, we get that
$$\sum_{(a,b)\in {S_\zahle[\alpha]}} e_{(p,r)}(a-b\alpha)=\sum_{(a,b)\in {S_\zahle[\alpha]}}\varphi_{\mathfrak{P}_i}(a-b\alpha)=\varphi_{\mathfrak{P}_i}\left(\prod_{(a,b)\in {S_\zahle[\alpha]}} (a-b\alpha)\right)$$
$$=\varphi_{\mathfrak{P}_i}(\xi^2)=2\varphi_{\mathfrak{P}_i}(\xi)\equiv 0 \mod 2.$$
\end{proof}
If our assumptions hold that means we can now use this outcome to do a similar process as we did for the rational factor base and find the final set $S_{\zahle[\alpha]}$ and use this to explicitly define the set $S=S_\zahle\cap S_{\zahle[\alpha]}$. However these have not been inconsequential and will all be explained in the following section. If we, however, suspend our disbelief for a moment longer we can finish the algorithm.\\
\\
Assume that we have found $(z^2,\xi^2)\in\zahle\times\zahle[\alpha]$ such that they fulfil equations (\ref{eq1}) and (\ref{eq2}). Then all that rests us to do is finding the square roots. First consider the rational case, so $z^2$ and $\zahle$. Then, by the sieving process, we know the prime factorization of $z^2$:
$$z^2=\prod_{p_i\in\zahle}p_i^{e_i},$$
but then it is trivial to find $z$:
$$z= \prod_{p_i\in\zahle}p_i^{\frac{e_i}{2}}.$$
\\
In the algebraic case, i.e. $\xi^2$ and $\zahle[\alpha]$, there is a bigger challenge, as there is no simple way to consider factorize $\xi^2$. One naive approach would be to compute the root of $x^2-\xi^2$, but this is usually not efficiently achievable. We instead take a theoretical approach, with an eye on computational efficiency, and let $q$ be an odd prime. If $f \mod q$ is irreducible in $\mathbb{F}_q[x]$, then $\zahle[\alpha]/q\zahle[\alpha]\cong\mathbb{F}_q[x]/(f\mod q)$. It can be shown that with significant probability there exists an odd $q$ such that $f \mod q$ is irreducible. To see this consider
\begin{thm}\label{3.8}
Let $f\in\zahle[x]$ be an irreducible polynomial of degree $d$, $d>1$. Then the density, inside the set of all prime numbers, of the set of prime numbers $q$ for which $f\mod q$ factors in $\mathbb{F}_q[x]$  into distinct irreducible non-linear factors exists and is at least $\frac{1}{d}$
\end{thm}
To prove this we need the following proposition:
\begin{prop}\label{3.9}
Let $G$ be a finite group that acts transitively on a finite set $\Omega$, with $\#\Omega=d>1$. Then there are at least $\frac{\#G}{d}$ elements of $G$ that act without fixed points on $\Omega$.
\end{prop}
Now we prove the theorem.
\begin{proof}
Let $\Gamma=\Gal(f/\quot)$, viewed as a permutation group of the set $A=\{\alpha_1,\ldots,\alpha_d\}$ with roots of $f$. For each prime number $q$ that does not divide the discriminant of $f$, there is a Frobenius element $\sigma_q\in\Gamma$ with the property that the degrees of the irreducible factors of $f\mod q$ are the same as the lengths of the cycles of the permutation $\sigma_q$. Hence, we are interested in those $q$ for which $\sigma_q$ acts without fixed points on $A$. Then by the Chebotarev Density Theorem, \ref{CDT}, 
every subset $S\subset G$ that is closed under conjugation, the set of prime numbers $\sigma_q$ belongs to has density $\frac{\#C}{\#G}$. The theorem follows from proposition \ref{3.9}.
\end{proof}
So except for the extremely small probability that we can not choose any $f$ for a specific $n$ this means we can assume that there exists a $q$ so $f\mod q$ is irreducible in $\mathbb{F}_q[x]$. The following theorem completes the square root finding process.
\begin{thm}\label{3.10}
Let $\mathfrak{q}=q\zahle[\alpha]$ be an ideal of $\zahle[\alpha]/q\zahle[\alpha]$. Then there exists a $\delta\in\zahle[\alpha]$
such that  
$$\delta^2\xi^2\equiv 1 \mod \mathfrak{q}.$$
\end{thm}
\begin{proof}
As $\zahle[\alpha]/q\zahle[\alpha]\cong\mathbb{F}_q[x]/f \mod q$ we know that $\abs{\zahle[\alpha]/q\zahle[\alpha]}=q^d$, where $d=\deg(f)$. Now consider 
$$\mathfrak{I}=q\zahle[\alpha]=\bigg\{\sum_{i=1}^{d-1} a_i\alpha^i:q\mid a_i\bigg\},$$
which is a degree $d$ prime ideal in $\zahle[\alpha]$. As $f \mod q$ is assumed irreducible it follows that $f'(\alpha)\not\in\mathfrak{I}$ and for each $(a,b)\in S_{\zahle[\alpha]}$ we have that $a-b\alpha\not\in\mathfrak{I}$ since $\gcd(a,b)=1$. Therefore $\xi^2=f'(\alpha)^2\prod_{(a,b)\in S_{\zahle[\alpha]}}(a-b\alpha)\not\in\mathfrak{I}$.\\
\\
Using Berlekamp's algorithm, \cite{CANT}, we find an element $\delta \mod \mathfrak{I}$ such that $\delta^2\xi^2\equiv 1 \mod \mathfrak{I}$, completing the proof.
\end{proof}
Hence there exists an element $\delta$ that is the inverse of a square modulo $\mathfrak{q}$. Now we can apply Newton--Rhapson iteration (\cite{CANT},\cite{BLP}) to find approximations such that
$$\delta_j\equiv\frac{\delta_{j-1}(3-\delta_{j-1}^2\xi^2)}{2}\mod (q\zahle[\alpha])^2.$$
In a finite number of steps we will find a $\gamma$ such that
$$\gamma\equiv \delta_j\xi \mod (q\zahle[\alpha])^2,$$
where $\gamma=\xi$ in $\zahle[\alpha]$. To complete the algorithm we compute $\gcd(z-N(\xi),n)$ and $\gcd(z+N(\xi),n)$ to find a factorization as described in section \ref{sec:GNFSintro}.
\subsection{Dealing with obstructions}\label{sec:GNFSobstr}
So far the technique we have been trying to describe can be captured in the following two bi-implications:
\begin{align}\label{eq3}
\text{Equation \eqref{eq1}}\Leftrightarrow \prod_{(a,b)\in T}(a-bm)\text{ has non-negative even exponents at all primes }p\leq B,
\end{align}
\begin{align}\label{eq4}
\text{Equation \eqref{eq2}}\Leftrightarrow \prod_{(a,b)\in T}(a-b\alpha)\text{ has even exponents at all prime ideals } \PriID\in\zahle[\alpha].
\end{align}
It is clear to see that the first of these bi-implications holds, as we can find the root by simply dividing all exponents by $2$. The second of these is not completely clear however. One of the first assumptions we made was that $\mathcal{O}_{\quot(\alpha)}=\zahle[\alpha]$, but this is rarely the case. In fact this is something that was believed true until in the 19'th century and even featured in (eventually disproven) proposed proofs of Fermat's Last Theorem \cite{WASH}. This assumption leads to four obstructions that have to be mitigated. For this let $\omega=\prod_{(a,b)\in T}(a-b\alpha)$.
\begin{enumerate}
\item $\zahle[\alpha]$ is not necessarily $\mathcal{O}_{\quot(\alpha)}$. At best we know that $\zahle[\alpha]\subseteq \mathcal{O}_{\quot(\alpha)}$, hence we can not assume that $\zahle[\alpha]$ is a Dedekind domain so $\omega \mathcal{O}_{\quot(\alpha)}$ might not be the square of an ideal in $\zahle[\alpha]$.
\item If $\omega\mathcal{O}_{\quot(\alpha)}$ is the square of some ideal $\mathfrak{I}$, it is not certain that $\mathfrak{I}$ is a principal ideal.
\item If $\omega\mathcal{O}_{\quot(\alpha)}$ is the square of some principal ideal $\gamma\mathcal{O}_{\quot(\alpha)}$, it is not certain that $\omega=\gamma^2$ as $\omega$ agrees with $\gamma^2$ only up to units of $\mathcal{O}_{\quot(\alpha)}$.
\item And even if $\omega=\gamma^2$, then we are not assured that $\gamma\in\zahle[\alpha]$
\end{enumerate}
Clearly these are obstructions that break the number field sieve. To deal with obstruction 4 recall the standard fact that for any $\theta\in\mathcal{O}_{\quot(\alpha)}$ and $f(x)=\irr_{\quot}(\alpha)$ we have that 
$$f'(\alpha)\cdot \theta\in\zahle[\alpha].$$
So we simply multiply equation \eqref{eq2} by $f'(\alpha)^2$ and equation \eqref{eq1} by $f'(m)^2$. The only condition we must apply is that $\gcd(f'(m),n)=1$ or else the resulting element given by equation \eqref{eq1} may not be invertible, however this is simply checked and if this is not the case then we have found a factorization of $n$ so we may simply assume that $f'(m)$ and $n$ are coprime.\\
\\
Now we attempt to subvert obstruction 1, which consequently allows us to negate both obstructions 2 and 3 as well. To do this we will consider the following chain:
$$V\supset V_1\supset V_2\supset V_3 = V\cap\left(\quot(\alpha)^\times\right)^2,$$
where $V$ is the group generated by $\quot(\alpha)^\times$ with even exponents, i.e. $e_{(p,r)}(v)\equiv 0\mod 2$ for all $v\in V$. This is a group with the following subgroups:
$$V_1=\{v\in V\mid v\mathcal{O}_{\quot(\alpha)}=\mathfrak{I}^2\text{ for some }\mathfrak{I}\subset\mathcal{O}_{\quot(\alpha)}\},$$
$$V_2=\{v\in V_1\mid v\mathcal{O}_{\quot(\alpha)}=\mathfrak\langle\vartheta\rangle^2\text{ for some }\langle\vartheta\rangle\subset\mathcal{O}_{\quot(\alpha)}\},$$
$$V_3=V\cap{(\quot(\alpha)^\times)}^2.$$
We attempt to bound the index $\left[V:V_3\right]$. Considered as a $\zahle$-module $\mathcal{O}_{\quot(\alpha)}$ is free of rank $d=\deg(f)$ and by definition $\zahle[\alpha]\subset\mathcal{O}_{\quot(\alpha)}$. It is a well-known identity in algebraic number theory that
$\Delta(f)=\left[\mathcal{O}_{\quot(\alpha)}:\zahle[\alpha]\right]^2\cdot\Delta$.
and as $[\quot(\alpha):\quot]>1$ it automatically follows that $\Delta>0$ hence $\left[\mathcal{O}_{\quot(\alpha)}:\zahle[\alpha]\right]$ is bounded by $\sqrt{\Delta(f)}$.
\begin{rmk}
Recall that $\Delta$ is the discriminant of the number field as in Definition \ref{2.38}
\end{rmk}
As we have factorization into prime ideals in $\mathcal{O}_{\quot(\alpha)}$ there is a bijection between the prime ideals $\mathfrak{Q}$ of $\mathcal{O}_\quot(\alpha)$ coprime to $\left[\mathcal{O}_{\quot(\alpha)}:\zahle[\alpha]\right]$ and the ideals $\mathcal{I}\subset\zahle[\alpha]$ coprime to this index. In fact if we only consider the prime ideals $\PriID\subset\zahle[\alpha]$ we get an isomorphism of local rings:
$$\zahle[\alpha]_\PriID\rightarrow{\left(\mathcal{O}_{\quot(\alpha)}\right)}_\mathfrak{Q},$$
$$\mathfrak{Q}\mapsto\mathfrak{P}=\mathfrak{Q}\cap\zahle[\alpha].$$
This allows us to generalize the definition of the map $e_{(p,r)}$ of proposition \ref{3.6}.
\begin{prop}\label{3.11}
Let $\mathfrak{P}\in\zahle[\alpha]$ be a prime ideal of $\zahle[\alpha]$. Then
$$e_\mathfrak{P}(a-b\alpha)=\sum_{\mathfrak{Q}\supset\mathfrak{P}}f(\mathfrak{Q}/\mathfrak{P})e_{(q,r)}(a-b\alpha),$$
where $\mathfrak{Q}$ lies over a prime $q\in\zahle$.
\end{prop}
In the case that $\mathfrak{P}$ does not divide the index $\left[\mathcal{O}_{\quot(\alpha)}:\zahle[\alpha]\right]$, $e_\mathfrak{P}(a-b\alpha)=e_{(q,r)}(a-b\alpha)$ as there will only be one $\mathfrak{Q}$ such that $\mathfrak{P}\subset\mathfrak{Q}$ and $f(\mathfrak{Q}/\mathfrak{P})=1$.\\
\\
Now consider the following map
$$V/V_1\rightarrow\bigoplus_{\mathfrak{Q}\mid\left[\mathcal{O}_{\quot(\alpha)}:\zahle[\alpha]\right]}\zahle/2\zahle,$$
$$x\mapsto\left(e_\mathfrak{Q}(x)\mod 2\right)_\mathfrak{Q}.$$
This is an injective homomorphism, hence $V/V_1$ is a $\mathbb{F}_2$-vectorspace of dimension bounded by $\abs{\big\{\mathfrak{Q}\mid\mathfrak{Q}\mid\left[\mathcal{O}_{\quot(\alpha)}:\zahle[\alpha]\right]\big\}}$. As the number of rational primes dividing the index is no more than $\frac{1}{2}\abs{\Delta(f)}$ and for each of these primes there are at most $\deg(f)$ prime ideals $\mathfrak{Q}\subset\mathcal{O}_{\quot(\alpha)}$ that divide it, we obtain
$$\dim_{\mathbb{F}_2}(V/V_1)\leq\frac{d}{2}\log\Delta(f).$$
This inequality is the first step to resolving the first obstruction. Further we will bound the dimensions $\dim_{\mathbb{F}_2}(V_1/V_2)$ and finally $\dim_{\mathbb{F}_2}(V_2/V_3)$ to obtain our final result.\\
\\
As the class group, $\text{Cl}_{\quot(\alpha)}$, is a finite abelian group and that its order $h$ is, \cite{STEV}, bounded by 
$$h<\sqrt{\abs{\Delta(f)}}\frac{d-1+\log\abs{\Delta(f)}^{d-1}}{(d-1)!}.$$
Define the map $\kappa:V_1\rightarrow\text{Cl}_{\quot(\alpha)}$ by $x\mapsto\mathfrak{I}$ where $x\mathcal{O}_{\quot(\alpha)}=\mathfrak{I}^2$. By definition $\ker(\kappa)=V_2$, hence 
$$\dim_{\mathbb{F}_2}(V_1/V_2)\leq\frac{\log h}{\log 2}.$$
Not that $V_2$ consists of elements that are squares in $\quot(\alpha)^\times$ up to units in $\mathcal{O}_{\quot(\alpha)}$, hence by [\cite{STEV}, 8.3]
$$\abs{V_2/V_3}\leq\abs{\mathcal{O}_{\quot(\alpha)}^\ast/\left(\mathcal{O}_{\quot(\alpha)}^\ast\right)^2}\leq d.$$
This leads to the following theorem:
\begin{thm}\label{3.12}
Let $V$ be as above and suppose that $n>d^{2d^2}>1$. Then the subgroup $V_3=V\cap\left(\quot(\alpha)^\times\right)^2$ of squares in $V$ satisfies
$$\dim_{\mathbb{F}_2}(V/V_3)\leq\log(n)^{\frac{3}{2}}$$
\end{thm}
This shows that the general number field sieve algorithm is certain to generate an element in $V$, however we are not sure yet that this is also a square, hence is in $V_3$. To obtain this we introduce the concept of quadratic characters. Let us start by proposing the idea over $\zahle$: Assume $x,y\in\zahle$ then $x$ is a quadratic residue modulo a prime $p$ if
$$x\equiv y^2\mod p.$$
This is easily tested by the Legendre symbol:
$$\bigg(\frac{x}{p}\bigg)=x^{\frac{p-1}{2}}\equiv(y^2)^{\frac{p-1}{2}}=y^{p-1}\equiv 1 \mod p.$$
Now assume that $x$ is a square, then it is also a square modulo $p$ for every prime $p$ such that $p\nmid x$. This means that for all $p\in\zahle$, $p$ prime:
$$\left(\frac{x}{p}\right)=1$$
So if there is at least one $p$ for which$(\frac{x}{p})=-1$ then $x$ is not a square modulo $p$ and by extension not a square element. This generalizes well to $\quot(\alpha)$.\\
\\
Let $\mathfrak{P}$ be a first degree prime ideal of $\zahle[\alpha]$ lying over $p$ and let $(p,R(p))$ be the set as defined in algorithm $\ref{Alg2}$. As there exists a ring homomorphism:
$\pi:\zahle[\alpha]\rightarrow\zahle/p\zahle$
with $\ker(\pi)=\mathfrak{P}$ by definition this gives rise to a Legendre symbol:
$$\left(\frac{\cdot}{\mathfrak{P}}\right):\zahle[\alpha]\rightarrow\zahle/p\zahle\rightarrow\{\pm 1, 0\}
,$$
such that for a non-square $x\in\zahle[\alpha]$ we have $\left(\frac{x}{\mathfrak{P}}\right)=-1$ with probability $\frac{1}{2}$. Restricting to Legendre symbols coming from $\mathfrak{P}$ over a prime $p\in\zahle$ such that $p>B'$ we avoid the character value $0$, and as a consequence of Chebotarev's Density Theorem we have that the Legendre symbols coming from $\mathfrak{P}$ are equidistributed over the space of homomorphisms from $V/V_3$ to $\{\pm 1\}$, denoted $\text{Hom}(V/V_3),\{\pm 1\})$.\\
\\
Now let $\chi_\mathfrak{P}=\left(\frac{\cdot}{\mathfrak{P}}\right)$ and consider the following lemma.
\begin{lemma}\label{3.13}
Let $k,r$ be non-negative integers, and let $E$ be a $k$-dimensional $\mathbb{F}_2$ vector space. Then the probability that $k+r$ elements that are uniformly at random drawn from $E$ form a spanning set for $E$ is at least $1-2^{-r}$
\end{lemma}
Then the equidistribution over $\text{Hom}(V/V_3,\{\pm 1\})$ with the obtained bound on the dimension of $V/V_3$ as a $\mathbb{F}_2$ vector space from theorem \ref{3.12} makes it overwhelmingly likely\footnote{See \cite{BLP} for details.} that a set of $B''=\lfloor \frac{3(\log n)}{\log 2}\rfloor$ quadratic characters span the homomorphism space. If they do, then the converse to the final theorem of this section shows that if $\beta\in\zahle[\alpha]\nozero$ satisfies $\chi_\mathfrak{Q}(\beta)=1$ for all first degree primes $\mathfrak{Q}$ with $2\beta\not\in\mathfrak{Q}$, then $\beta$ is in fact a square in $\zahle[\alpha]$. Furthermore, there can be a finite number of exceptions which is explored greater in \cite{BLP}.
\begin{thm}\label{3.14}
Let $S$ be a finite set of integer pairs $(a,b)$ such that $\gcd(a,b)=1$ fulfilling that for some $\gamma\in\quot[\alpha]$:
$$\omega=\gamma^2.$$
Moreover let $\mathfrak{q}$ be a first degree prime ideal corresponding to the pair $(s,q)$ that does not divide $\langle a-b\alpha\rangle$ for any pair $(a,b)$ and for which $f'(s)\not\equiv 0\mod q$.
Then we have 
$$\prod_{(a,b)\in S} \bigg(\frac{a-bs}{q}\bigg)=1.$$
\end{thm}
\begin{proof}
Let $\varphi:\zahle[\alpha]\rightarrow\zahle/q\zahle$ be the ring homomorphism given by $\varphi(\alpha)=r \mod q$ and let $\ker(\varphi)=\mathfrak{Q}$. Then $\mathfrak{Q}$ is a first degree prime ideal corresponding to the pair $(q,r)$. Let $\chi_\mathfrak{Q}:\zahle{\alpha}\rightarrow\{\pm 1\}$ and let $\psi_\mathfrak{Q}:\zahle[\alpha]/\mathfrak{Q}\rightarrow\zahle/q\zahle\nozero$, such that $\chi_\mathfrak{Q}=\Leg_\mathfrak{Q}\circ\psi_\mathfrak{Q}$ where $\Leg_\mathfrak{Q}:\zahle/q\zahle\rightarrow\{\pm 1\}$ is the Legendre symbol. Then
$$\chi_\mathfrak{Q}(a-b\alpha)=\left(\frac{a-br}{q}\right).$$
Letting, as discussed, 
$$\xi^2=f'(\alpha)^2\prod_{(a,b)\in S}(a-b\alpha),$$
for some $\xi\in\zahle[\alpha]$. By the hypothesis the factors on the left are not in $\mathfrak{Q}$, so we have that $\xi\not\in\mathfrak{Q}$. However
$$\chi_\mathfrak{Q}(\xi^2)=1,$$
$$\chi_\mathfrak{Q}(f'(\alpha)^2)=1,$$
so
$$\chi_\mathfrak{q}\left(\prod_{(a,b)\in S}\left(\frac{a-br}{q}\right)\right)=1.$$
This completes the proof. 
\end{proof}
So now we can be almost certain that we can find a square in $\zahle[\alpha]$ and therefore have mitigated all obstructions. The only question left is: What is the complexity?
\subsection{Computational efficiency}
Let us record the number field sieve algorithm step-by-step:
\begin{algorithm}[H]\label{GNFS}
\caption{General Number Field Sieve}\label{alg3}
\begin{algorithmic}
\Input{$n\in\nat$, $n$ odd composite, degree $d>1$, universe $U$, smoothness bounds $B$ and $B'$}
\Output{A factor of $n$ or $\neg$ for failure.}
\State Choose $m=\lfloor n^{\frac{1}{d}}\rfloor$.
\State Define $f(x)\in\zahle[x]$ as the polynomial with coefficient corresponding to the base-$m$ expansion of $n$.
\For{each $(a,b)\in U$}
\State Run algorithm \ref{Alg1} to find $S_\zahle$
\State Run algorithm \ref{Alg2} to find $S_{\zahle[\alpha]}$
\State Use the techniques from section \ref{sec:GNFSobstr} to find quadratic character base $S_Q$ such that $\abs{S_Q}=B''$.
\EndFor
\If{$\abs{S_\zahle\cap S_{\zahle[\alpha]}}>\abs{B}+\abs{B'}+1$}
\State Use a reduction algorithm over $\mathbb{F}_2$ to find an dependent subset $S\subset S_\zahle\cap S_{\zahle[\alpha]}$
\State Compute $z^2=\prod_{(a,b)\in S} (a-bm)$ and use its known prime factorization to find $z$
\State Compute $\zeta^2=\prod_{(a,b)\in S}(a-b\alpha)$ and use Newton--Rhapson iteration to find $\zeta$
\If{ $z\in\zahle$ \textbf{and} $\zeta\in\zahle[\alpha]$}
\State \Return $\gcd(z-\zeta,n)$, $\gcd(z+\zeta,n)$
\Else{ find a new dependent subset $S$. If all such sets produce no result \Return $\neg$}
\EndIf
\Else{ \Return{$\neg$}}

\EndIf
\end{algorithmic}
\end{algorithm}
It is generally considered difficult to do a complexity analysis on the number field sieve as it is stated in all its generality. Already when it was analysed in \cite{BLP} there was a conjectured complexity of 
$$L_n\left(\frac{1}{3}, \sqrt[3]{\frac{64}{9}}+\loh(1)\right),$$
and over the 30 years that followed it showed that the conjectured complexity held up heuristically, however there are some things that we are able to say. Following \cite{STEV} we are able to carefully choose $f$ and thereby $\deg(f)$ to optimize the smoothness bound $B^\ast=\max{B,B'}$ and parameter $u$ for the universe of sieving $U$. The ``basic" cost of the algorithm is computed to be $u^{2+o(1)}+y^{2+o(1)}$.\\
\\
First we will address step 6. As the binary matrix that is formed will be sparse, i.e.\ there will be significantly more zeroes than ones, it allows us to use fast algorithms to find dependencies. Two of these extremely fast algorithms are Block--Lanczos (\cite{LANC},\cite{LANC2}) and Block--Wiedemann \cite{WIED} which both have similar complexity. The size of these matrices are at most $B^\ast\times B^\ast$, and by choosing $\log B^\ast\approx\log u$ we balance the contributions of $U$ and the matrix reduction.\\
\\
Looking at the sieving procedure we can see that a pair $(a,b)$ is $B^\ast$-smooth if $(a-bm)N(a-b\alpha)$ is $B^\ast$-smooth, and thus, using that the size of $m$ is approximately $\sqrt[d]{n}$, we may bound this product by
$$u\sqrt[d]{n}(d+1)u^d\sqrt[d]{n}\approx n^{\frac{2}{d}}u^{d+1},$$
by choosing $d$ optimally, namely
$$d\approx \bigg(\frac{3\log n}{\log\log n}\bigg)^{\frac{1}{3}}.$$
The probability that a number $x$, of the right size, is smooth can be approximated by $r^{-r}$ where $r=\frac{\log x}{\log y}$. In order to maximize this probability we choose
$$r=\frac{\log x}{\log u}\approx \frac{d'^2}{d'}+d'+1,$$
for $d'=\sqrt{\frac{2\log n}{\log u}}$.\\
\\
This means we get dependency in the matrix if we have at least $y\approx u^2r^{-r}$ $(a,b)$ pairs. As we have assumed $\log y\approx \log u$ taking logarithms we get $\log u\approx r\log r$, equivalently $r\approx \frac{\log u}{\log\log u}$. Hence, by taking $\frac{2}{3}$ powers:
$$\log u(\log\log u)^\frac{1}{3}\approx 2(\log n)^\frac{1}{3}. $$
As we have, for $a,s,t\in\reals$, that  $s=t(\log t)^a$ goes to $t=(1+o(1)s(\log s)^{-a}$ as $t$ goes to infinity, we get that 
$$\log y\approx \log u\approx 2(\log n)^{\frac{1}{3}}\bigg(\frac{1}{3}\log\log n\bigg)^{\frac{2}{3}}=\bigg(\frac{8}{9}\bigg)^{\frac{1}{3}}(\log n)^{\frac{1}{3}}(\log\log n)^{\frac{2}{3}}.$$
The final GCD computation, which even with the most simple implementation: Euclid's algorithm, has a complexity at most $O(\log n)$ and as such can be completely disregarded in the complexity analysis.\\
\\
Hence, with this choice of basic parameters $u,y$ and $d\approx d'\approx\big(\frac{3\log n}{\log\log n}\big)^{\frac{1}{3}}$ we get the conjectured runtime of $L_n\left(\frac{1}{3},\sqrt[3]{\frac{64}{9}}+\loh(1)\right)$, however this is fully heuristic and very little strong can be proven rigorously for this particular form of the number field sieve even with modern techniques.
\newpage
\section{Preparing for randomness}\label{ch: PREL}
\subsection{Zeta function and the Prime Number Theorem}\label{sec: PRELzeta}
Before now we have been able to suffice with elementary and algebraic number theory and some Galois theory. To introduce randomness we will be drawing strongly on analytic number theory. It was Riemann who really gave analytic number theory a kickstart by introducing the Riemann zeta function.
\begin{rmk}
Unless otherwise noted $p,q\in\zahle$ are prime numbers.
\end{rmk}
\begin{defn}\label{4.1}
Let $s\in\compl$ and $n\in\nat$, then the Riemann zeta function is defined as
$$\zeta(s)=\sum_{n=1}^\infty\frac{1}{n^s}.$$
This converges for all $s$ such that $\text{Re}(s)>1$, has an analytic continuation to $\compl$ except for a simple pole at $s=1$, and admits a functional equation
$$\pi^{-\frac{s}{2}}\Gamma\left(\frac{s}{2}\right)\zeta(s)=\pi^{-\frac{1-s}{2}}\Gamma\left(\frac{1-s}{2}\right)\zeta(1-s),$$
where $\Gamma(s)=\int_0^\infty e^{-t}t^{s-1}\, dt$. Moreover for $\text{Re}(s)>1$, $\zeta(s)$ admits an Euler product:
$$\zeta(s)=\prod_{p}\left(\frac{1}{1-p^{-s}}\right).$$
\end{defn}
\begin{rmk}
The proofs for all the properties listed in the definition can be found in any undergraduate analytic number theory book. We suggest \cite{MUNT} or \cite{DAVE}.
\end{rmk}
One of the first results of analytic number theory is the prime number theorem.
\begin{thm}\label{4.2}
Let $\pi(x)=\sum_{p\leq x}1$ be the prime counting function. Then there exists a constant $a>0$ such that
$$\pi(x)=\text{Li}(x)+\boh\left(\frac{x}{\log x}\exp\left(-a\sqrt{\log x}\right)\right)=\frac{x}{\log x}\left(1+\loh(1)\right),$$
where
$$\text{Li}(x)=\int_2^x\frac{1}{\log t}\, dt=\frac{x}{\log x}\left(1+\boh\left(\frac{1}{\log x}\right)\right)$$
is the logarithmic integral.
\end{thm}
One of the most influential open questions in mathematics is related to the zeta function: The Riemann Hypothesis.
\begin{prop}\label{RH} \textbf{ [Riemann Hypothesis] }
Let $\zeta(s)$ be the Riemann zeta function. For $\text{Re}(s)>0$, if $\zeta(s)=0$ then $\text{Re}(s)=\frac{1}{2}$.
\end{prop}
This hypothesis has been so extensively tested that it is often accepted as true in the mathematical community. Doing so there is a whole section of number theory that is considered `conditional', which will become unconditional the moment that the Riemann Hypothesis is proven. There is a very interesting result by Hardy, proven in 1914, which we will record for good measure, [theorem 14.8, \cite{MUNT}].
\begin{thm}\label{4.4}\textbf{[Hardy's Theorem]}
There exist infinitely many $t\in\reals$ such that
$$\zeta\left(\frac{1}{2}+it\right)=0.$$
\end{thm}
However for the Riemann Hypothesis to be true this is not enough, as we want all the zeroes to be on the critical line of $\Re(s)=\frac{1}{2}$. This is still one of the major research areas in analytic number theory, but there are some results. The first such was by Selberg in 1942, who showed that there is some non-zero fraction of zeroes on the critical strip. This was later improved by Levinson who showed that at least $34.7\%$ of the zeroes lie on this line in 1974-1975, and this was improved to $40\%$ by Conrey in 1989. Since then there have been, to my knowledge, no significant improvements leaving this topic wide open for future research.
\subsection{Generalizing the Riemann zeta function}\label{sec: PRELdiri}
\subsubsection{Dirichlet characters} \label{sec: PRELchar}
The Riemann zeta function is not the only function of its kind. The idea of the Riemann zeta function was generalized in two directions each with their own implications. Dirichlet introduced a generalization through the use of multiplicative character functions and Dirichlet series to define a whole family of meromorphic analytic arithmetic functions. Dedekind on the other hand chose the algebraic route and defined the Dedekind zeta functions, which are defined over number fields.  For the former we shall introduce Dirichlet characters, which extends the idea behind quadratic characters.
\begin{defn}\label{4.5}
Let $G$ be a group. A character on $G$ is a group homomorphism $\chi:G\rightarrow\compl^\times$ and the set of characters of $G$ is denoted $\hat{G}$. 
\end{defn}
Generally characters possess a few properties:
\begin{prop}\label{4.6}
Let $\chi$ be a $G$-character, then:
\begin{enumerate}
\item As $\chi$ is a homomorphism: $\forall g,g'\in G: \chi(gg')=\chi(g)\chi(g')$
\item $\forall G, \forall \chi: \chi(e_G)=1$
\item $\forall G\exists\chi_0\text{ s.t. } \forall g\in G:\chi_0(g)=1$. This is called the trivial or principal character.
\end{enumerate}
\end{prop}
These properties reveal that $\hat{G}$ might have a group structure, and this is true for the finite case:
\begin{lemma}\label{4.7}
If $G$ is a finite group then so is $\hat{G}$.
\end{lemma}
There is one more property to consider for general characters before we can introduce Dirichlet characters and that is orthogonality:
\begin{defn}\label{4.8}
Let $G$ be a finite group with set (and as $G$ is finite: group) of characters $\hat{G}$. Then $G$ is said to have the orthogonality property if the following two equations hold:
$$\sum_{g\in G}\chi(g)=\begin{cases}
\# G, & \chi=\chi_0\\
0, & \chi\neq\chi_0
\end{cases}$$
$$\sum_{\chi\in\hat{G}}\chi(g)=\begin{cases}
\#\hat{G} &, g=e_G\\
0  & ,g\neq e_G
\end{cases}$$
\end{defn}
\begin{rmk}
For the sake of clarity: Note that the first sum requires a fixed character and sums over all elements, while the second requires a fixed element and sums over all characters of $G$.
\end{rmk}
It is an elementary result from representation theory that the orthogonality property holds for all finite groups, and in particular finite cyclic groups, which is of importance as we shall now see:
\begin{defn}\label{4.9}
\begin{enumerate}
\item []
\item Let $q\in\nat$. Then a Dirichlet character $\mod q$ is a character of the form $\chi:\left(\zahle/q\zahle\right)^\times\rightarrow\compl$.
\item Let $\chi$ be a Dirichlet character $\mod q$, we extend $\chi$ to an arithmetic function $\chi:\zahle\rightarrow\compl$ as 
$$\chi(n)=\begin{cases}
\chi(n\mod q) & ,\gcd(n,q)=1\\
0 & ,\gcd(n,q)>1
\end{cases}.$$
\end{enumerate}
Moreover a Dirichlet character mod $q$ has the orthogonality property:
$$\sum_{\substack{n=1\\ \gcd(n,q)=1}}^d\chi(n)=\begin{cases}
\phi(q) & ,\chi=\chi_0\\
0 & ,\chi\neq\chi_0
\end{cases},$$
and
$$\sum_{\chi\mod q}\chi(n)=\begin{cases}
\phi(q) & ,n\equiv 1\mod q\\
0 & ,otherwise
\end{cases}.$$
\end{defn}
As the Dirichlet characters are defined over $\zahle/q\zahle$, and we will often consider $q$ a prime, we have a vested interest in how characters of cyclic groups behave. It turns out that these are rather well-behaved:
\begin{thm}\label{4.10}
Assume that $G$ is a finite cyclic group such that $\abs{G}=n$ and $\langle a \rangle=G$. Then:
\begin{enumerate}
\item $\hat{G}$ has exactly $n$ elements:
$$\chi_k(a^m)=\zeta_n^{km}, k\in\{1,\ldots,n\}.$$
where $\zeta_n=e^{\frac{2\pi i}{n}}$ is the n'th root of unity.
\item $G$ has the orthogonality property (either as defined in Definition \ref{4.8}
or as a Dirichlet character
\item $\hat{G}$ is cyclic with generator $\chi_1$, hence $G\cong \hat{G}$.
\end{enumerate}
\end{thm}
\begin{proof}
Let $\chi\in\hat{G}$. Then $\chi(a)=\zeta^k_n$ for some $k\in\{1,\ldots,n\}$. Hence 
$$\chi(a^m)=\chi(a)^m=\zeta_n^{km},$$
proving part (1). By (1) $\hat{G}$ is cyclic and generated by $\chi_1$, so $G\equiv\hat{G}$ as required.  To see that $G$ has orthogonality of characters it is simple to see that $\chi$ follows the definition of orthogonality. To illustrate we show it according to Definition \ref{4.8}. Note that 
$$\sum_{g\in G}\chi_0(g)=\abs{G}$$
is trivial. Hence assume that $k\in\{1,\ldots,n-1\}$. Then
$$\sum_{g\in G}\chi(g)=\sum_{m=0}^{n-1}\chi_k(a^m)=\sum_{m=0}^{n-1}\zeta^{km}_n=\frac{1-\zeta^{kn}_n}{1-\zeta^k_n}=0.$$
By the duality of this identity the identity of $\hat{G}$ holds too, hence we have shown (2).
\end{proof}
If $q$ is not prime we will still be dealing with a special group, as $\left(\zahle/q\zahle\right)^\times$ will still be abelian, hence can be written as a direct product of cyclic groups. Therefore we can use theorem \ref{4.10} to say the following:
\begin{lemma}\label{4.11}
Let $G_1,G_2$ be finite cyclic groups and let $G\cong G_1\times G_2$. Let $\chi_i\in\hat{G}_i$ for $i=1,2$ and define $\chi:G\rightarrow\compl^\times$ via $\chi(g_1,g_2)=\chi_1(g_1)\chi_2(g_2)$. This is a character. Conversely, if $\chi\in\hat{G}$ then there exists a unique choice of $\chi_1\in\hat{G}_1$ and $\chi_2\in\hat{G}_2$ such that $\chi(g)=\chi_1(g_1)\chi_2(g_2)$. Furthermore $\hat{G}\cong\hat{G}_1\times\hat{G}_2$ and $G$ has orthogonality of characters and if $G$ is finite abelian then $G\cong\hat{G}$.
 \end{lemma}
It is clear that if $q$ is an odd prime, then we can define 
$$\chi(n)=\left(\frac{n}{p}\right),$$
where the right hand side is the Legendre symbol. This is a character, and it is easy to show that this is, in fact, a Dirichlet character. It is therefore that Dirichlet characters are a generalization of quadratic characters, which will be extremely useful when we randomize the number field sieve. Another very useful property is that, by definition, Dirichlet characters are (quasi-)periodic.
\begin{defn}\label{4.12}
Let $\chi$ be a Dirichlet character $\mod q$. We say that $d\in\nat$ is a quasi-period if for all $m\equiv n\mod d$ such that $\gcd(mn,q)=1$: $\chi(m)=\chi(n)$. Moreover the least such $d$ is called the conductor.
\end{defn}
It is obvious from the definition that if there are two quasi-periods $d_1,d_2$ then there is a third: $\gcd(d_1,d_2)$. As $q$ is always a quasi-period we have that the conductor of a Dirichlet character $\mod q$ is at most $q$. We can however say something much stronger about the conductor of $q$:
\begin{prop}\label{4.13}
Let $\chi$ be a Dirichlet character $\mod q$ with conductor $d$. Then $d\mid q$.
\end{prop}
\begin{proof}
Let $g=\gcd(d,q)$. Suppose $m\equiv n \mod g$ and $\gcd(mn,q)=1$. By Euclid's algorithm there exists $x,y\in\zahle$ such that $m-n=dx+qy$, thus
$$\chi(m)=\chi(m-qy)=\chi(dx+n)=\chi(n).$$
Hence $g$ is a quasiperiod of $\chi$. As $d$ is the conductor that means $d\mid g$ and since $g=\gcd(d,q)$ we have that $d\mid q$.
\end{proof}
As per usual we have now found a way to consider characters ``prime":
\begin{defn}\label{4.14}
A Dirichlet character modulo $q$ is called primitive if it has conductor $q$.
\end{defn}
\begin{rmk}
Convention dictates that the trivial character $\chi_0$ is not primitive.
\end{rmk}
This idea of primitive characters is especially useful when you are working with induced characters. Taking inspiration from the notation as it is used in representation theory consider $d\mid q$ and consider $\chi^\ast$, a character $\mod d$, and define
$$\chi(n)=\begin{cases}
\chi^\ast(n) & ,\gcd(n,q)=1\\
0 & , otherwise
\end{cases}.$$
Then $\chi(n)$ is a multiplicative character and has period $q$, so it is a Dirichlet character $\mod q$. This is called the Dirichlet character $\mod q$ induced by a Dirichlet character $\mod d$ and we will denote it $\Ind_d^q(\chi)$. We conclude with the following theorem about induced Dirichlet characters:
\begin{thm}\label{4.15}
Let $\chi$ be a Dirichlet character $\mod q$ with conductor $d$. Then there exists a unique character $\chi^\ast \mod d$ such that $\chi(n)=\Ind_d^q(\chi^\ast(n)).$
\end{thm}
\subsubsection{Dirichlet $\mathcal{L}$-functions}\label{sec: PRELlfun}
The definition of Dirichlet characters leads directly to the definition of Dirichlet $\mathcal{L}$-functions:
\begin{defn}\label{4.16}
Let $q\in\nat$ and let $\chi$ be a Dirichlet character modulo $q$. Then the Dirichlet $\mathcal{L}$-function associated to $\chi$, denoted $\mathcal{L}(s,\chi)$, is defined as
$$\mathcal{L}(s,\chi)=\sum_{n=1}^\infty\frac{\chi(n)}{n^s},$$
for $\mathfrak{R}(s)>1$.
\end{defn}
As $\chi$ is a completely multiplicative character function we can immediately conclude that $\mathcal{L}(s,\chi)$ is absolutely convergent for $\mathfrak{R}(s)>1$ and admits an Euler product in this region:
$$\mathcal{L}(s,\chi)=\prod_p\left(1-\frac{\chi(p)}{p^s}\right)^{-1},$$
where $p$ is prime. The behaviour of the trivial character and the non-trivial characters is quite different. In the case of the trivial character we can rewrite the Euler product:
$$\mathcal{L}(s,\chi_0)=\prod_p\left(1-p^{-s}\right)\zeta(s).$$
This means $\mathcal{L}(s,\chi_0)$ behaves like $\zeta(s)$ and allows for a meromorphic continuation with a single simple pole at $s=1$. When $\chi$ is not trivial the behaviour changes radically.
\begin{thm}\label{4.17}
Let $q\in\nat$ and let $\chi$ be a Dirichlet character modulo $q$ such that $\chi\neq\chi_0$. Then $\mathcal{L}(s,\chi)$ converges and is holomorphic for $\mathfrak{R}(s)>0$.
\end{thm}
\begin{rmk}
$\mathfrak{R}(s)=0$ is the abscissa of convergence for $\mathcal{L}(s,\chi)$, but despite the convergence we can not extend the definition of the Euler product  of $\mathcal{L}(s,\chi)$ to $\mathfrak{R}(s)<1$.
\end{rmk}
Knowing that we have convergence for $\mathfrak{R}(s)=1$ in this case one might be interested in what these values actually are. This has been a problem for analytic number theorists for quite some time, but we can say the following:
\begin{thm}\label{4.18}
Let $\quot\in\nat$ and let $\chi$ be a Dirichlet character mod $q$. Then
$$\prod_{\chi\neq\chi_0}\mathcal{L}(1,\chi)\neq 0.$$
\end{thm}
This is absolutely not trivial and a particularly concise proof is presented in [\cite{MUNT}, theorem 4.9]. The first of two jewels in this is the following theorem by Dirichlet.
\begin{thm}\label{4.19}
Let $q\in\nat$ and let $a\in\zahle$ such that $\gcd(a,q)=1$. Then there are infinitely many primes $p\equiv a\mod q$
\end{thm}
\begin{rmk}
Despite not denoting the proof explicitly one the important step is that we show that
$$\sum_{p\equiv a\mod q}\frac{1}{p}$$
diverges. The theorem follows immediately.
\end{rmk}
It might not surprise that with the generalization of the Riemann zeta function there follows a generalization of the Riemann Hypothesis:
\begin{prop}\label{GRH}
\textbf{ [Generalized Riemann Hypothesis] }Let $\mathcal{L}(s,\chi)$ be a Dirichlet $\mathcal{L}$-function with $\chi$ a primitive character $\mod q, q>1$. Then for $\text{Re}(s)\in\left(0,1\right)$, if $\mathcal{L}(s,\chi)=0$ then $\text{Re}(s)=\frac{1}{2}$.
\end{prop}
We will come back to Dirichlet $\mathcal{L}$-functions when we adapt the prime number theorem to arithmetic progressions in Section \ref{sec: PRELarit}.
\subsubsection{Dedekind zeta functions and Siegel zeroes}\label{sec: PRELdede}
Now we consider the other generalization of the zeta function, though it can in this case very accurately be called an extension. 
\begin{defn}\label{4.21}
Let $K$ be an algebraic number field and let $\mathfrak{I}\subset K$ be an ideal. Then the Dedekind zeta function over $K$ is defined as
$$\zeta_K(s)=\sum_{\mathfrak{I}\subseteq \mathcal{O}_K}\frac{1}{\left(N_K(\mathfrak{I})\right)^s}.$$
\end{defn}
The Dedekind zeta function is for $\mathfrak{R}(s)>1$ and even admits an Euler product for that region. To see this let $\PriID\subseteq K$ be a prime ideal, then:
$$\zeta_K(s)=\prod_{\PriID\subseteq\mathcal{O}_K}\left(1-N_K(\PriID)^{-s}\right).$$
It was proven by Hecke that the Dedekind zeta function admits an analytic continuation to the whole complex plane and also admits a functional equation, so it is clear that in many ways the Dedekind zeta function behaves like the Riemann zeta function. This correlation is so close that it culminates in the following theorem:
\begin{thm}\label{4.22}
\textbf{Landau, 1903 - Prime Ideal Theorem}\\
Let $K$ be an algebraic number field and let $\mathfrak{P}\subseteq\mathcal{O}_K$ be a prime ideal and define $\pi_K(x)=\abs{\{\PriID\in\mathcal{O}_K\mid N(\PriID)\leq x\}}$. Then for $x\rightarrow\infty$, 
$$\pi_K(x)\asymp \frac{x}{\log x}.$$
\end{thm}
For our purposes it is sufficient to realize that the Dedekind zeta function over $K$ behaves a lot like the Riemann zeta function, and that the same language can be used to talk about the Dedekind zeta function as what we use for the Riemann zeta function. In fact, there is also a Riemann Hypothesis for the Dedekind zeta functions:
\begin{prop}\label{ERH}
\textbf{ [Extended Riemann Hypothesis] }
Let $K$ be a number field. For $\text{Re}(s)\in\left(0,1\right)$, if $\zeta_K(s)=0$, then $\text{Re}(s)=\frac{1}{2}.$
\end{prop}
\begin{rmk}
As the statements are equivalent for the Dirichlet $\mathcal{L}$-functions and the Dedekind zeta functions it is not uncommon to see the Extended Riemann Hypothesis be called the Generalized Riemann Hypothesis for number fields. We will in this paper refer to both as the Generalized Riemann Hypothesis, GRH, and let context infer which version we are referring to.
\end{rmk}
\subsubsection{Siegel zeroes} \label{sec: PRELsieg}
One of the main problems that we encounter when we consider both Dirichlet $\mathcal{L}$-functions and Dedekind zeta functions is that we can not generalize the zero free region. Especially, for the Riemann zeta function we were able to extend the zero-free region just past the $\mathfrak{R}(s)=1$ line, and this may not be true for Dedekind zeta functions and Dirichlet $\mathcal{L}$-functions. In fact:
\begin{thm}\label{4.24}
There is a constant $c>0$ such that if $\mathcal{L}(\sigma+it,\chi)=0$ for some primitive complex Dirichlet character $\chi\mod q$ then
$$\sigma<1-\frac{c}{\log q(\abs{t}+2)}.$$
If $\chi$ is a real primitive character then this holds for all zeroes of $\mathcal{L}(s,\chi)$ with at most one exception. The exceptional zero, if it exists is real and simple.
\end{thm}
Later we will see a result by Landau, \ref{5.28}, which develops on this idea. It is clear to see that this theorem is formulated for Dirichlet $\mathcal{L}$-functions, but a similar result exists for Dedekind zeta function. We therefore introduce the following concept:
\begin{defn}\label{4.25}
Let $\mathcal{L}(s,\chi)$ be a Dirichlet $\mathcal{L}$-function, then the hypothetical value $1-\nu$ for $\nu\in(0,\frac{1}{2})$, such that
$$\mathcal{L}(1-\nu,\chi)=0,$$
for a Dirichlet character modulo $q$, $\chi$, is called a Siegel-zero of $\mathcal{L}(s,\chi)$.\\
\\
Equivalently, for $\zeta_K(s)$, if there exists a $\nu\in\left(0,\frac{1}{2}\right)$ such that
$$\zeta_K(1-\nu)=0,$$
for a number field $K$, is called a Siegel-zero of $\zeta_K(s)$.
\end{defn}
The culminating theorem from Siegel on `his' zeroes was
\begin{thm}\label{4.26}\textbf{Siegel's theorem}
For any $\epsilon>0$ there exists a positive number $c(\epsilon)$ such that, if $\chi$ is a real primitive character $\mod q$, then
$$L\left(1,\chi\right)>\frac{c(\epsilon)}{q^{-\epsilon}}.$$
\end{thm}
Even in the definition we elude to these zeroes as hypothetical, as the assumption of the (Generalized) Riemann Hypothesis has as an immediate consequence that such zeroes do not exist. For the sake of this paper \cite{RNFS} does not assume the Riemann Hypothesis in any way, hence we must always concern ourselves with possible Siegel zeroes.
\subsection{Arithmetic progressions}\label{sec: PRELarit}
Now we turn to the next topic, which is a main factor in the Randomized Number Field Sieve: Smooth numbers on arithmetic progressions. To be able to do this we first need to know a little bit about prime numbers on arithmetic progressions and their distribution. First we note that it would be beneficial for us to have a prime counting function for arithmetic progressions. This turns out to be a simple modification:
\begin{defn}\label{4.27}
For $\gcd(a,q)=1$ the prime counting function for arithmetic progressions is defined as
$$\pi_{(a,q)}(x)=\pi(x,a,q)=\sum_{\substack{p\leq x\\p\equiv a\mod q}}1.$$
\end{defn}
The first question that comes to mind is one that was answered by Dirichlet: Are there infinitely many primes of that form? This is the case, as we stated in theorem \ref{4.19}:
\begin{thm}\label{4.28}
Let $x\in\reals$ and $a,q\in\zahle$ such that $\gcd(a,q)=1$. Then there are infinitely many primes of the form $p\equiv a\mod q$.
\end{thm}
This immediately shows that just like for $\pi(x)$ we have that $\pi(x,a,q)$ diverges. Therefore, we infer, there must be something interesting to say about the asymptotic behaviour of the prime counting function. Not only is there something interesting to say, there is a whole lot to say as well:
\begin{thm}\label{4.29}
Let $x\in\reals$, $a,q\in\zahle$ such that $\gcd(a,q)=1$, and let $\phi(x)$ be Euler's totient function. Then 
$$\pi(x,a,q)\sim\frac{x}{\log x}\frac{1}{\phi(q)}.$$
\end{thm}
This immediately shows that the density is independent of $a$ as long as $\gcd(a,q)=1$.                 
\subsubsection{Extensions and Refinements}\label{sec: PRELexre}
Before we start to talk about smooth numbers we want to consider if anything better than Dirichlet's function can be said, and we can if we allow for an ineffective bound. For this we define 
$$\psi(x,a,q)=\sum_{\substack{n\leq x\\n\equiv a\mod q}}\Lambda(n),$$
where $\Lambda(n)$ is the von Mangoldt function:
$$\Lambda(n)=\begin{cases}
\log p& \exists k\in\nat:n=p^k\\
0 & \text{ otherwise}
\end{cases}.$$
\begin{thm}\label{SW}
\textbf{Siegel-Walfisz theorem} Let $\psi(x,a,q)$ and $\Lambda(n)$ be as defined above. Given $q\in\nat$ there exists a positive constant $c(q)$ such that
$$\psi(x,a,q)=\frac{x}{\phi(q)}+\boh\left(x\exp\left(-c(q)\sqrt{\log x}\right)\right).$$
\end{thm}
The reason this is ineffective is because it is based on Siegel's theorem, \ref{4.26}, and this theorem gives no way to compute $c(n)$. Another way to go about this is looking at average cases. This probabilistic approach will be explored later and will lead to effective constant and a far stronger bound.\\
\\
It may not come as a surprise that the idea of prime numbers along arithmetic progressions can be expanded to number fields, as we have done with the Prime Ideal Theorem. For this we first consider what an `arithmetic progression' is on a number field. This does not come intuitively, but the answer lies in Galois Theory. Let $L/K$ be a Galois extension with Galois group $G=\text{Gal}(L/K)$. As $L/K$ is a  Galois extension the Frobenius element, $\Frob_\mathfrak{P}$, defines a conjugacy class
$$C=\bigg\{\text{Frob}_\mathfrak{Q}\mid \mathfrak{Q}\subset L\text{ s.t. } \mathfrak{Q}\text{ is a prime ideal and } \mathfrak{Q}\mid\mathfrak{P}\bigg\}.$$
It can be shown that this abides the rules of modular arithmetic and therefore can be used as an extension of the idea of arithmetic progressions. Using this we get the following result:
\begin{thm}\label{CDT}\textbf{Chebotarev Density Theorem}
Let $L/K$ be Galois and  let $\PriID\subseteq K$ be a prime ideal. Moreover let $C\subseteq G$ be the conjugacy class defined above. Then
$$\{\mathfrak{P}\mid \PriID\nmid\Delta_{L/K},\text{Frob}_\PriID\in C\}$$
has density $\frac{\abs{C}}{\abs{G}}$.
\end{thm}
This has long been the best result available, but over the years it has been strengthened both under the GRH assumption and without. Later we will see a version which will be of particular interest to us as it is both unconditional and has effectively defined constants.
\subsubsection{Smooth numbers on Arithmetic Progressions}
Lastly, before we move onto a different topic, we discuss how smooth numbers behave on arithmetic progressions. First we need some definitions:
\begin{defn}\label{4.32}
For $x,y,a,q\in\nat$ and a Dirichlet character $\chi$, we define:
$$\Psi(x,y)=\abs{\{z\in\nat\mid z<x,z\text{ is }y\text{-smooth}\}},$$
$$\Psi_r(x,y)=\abs{\{z\in\nat\mid z<x,z\text{ is }y\text{-smooth},\gcd(z,r)=1\}},$$
$$\Psi(x,y,a,r)=\abs{\{z\in\nat\mid z<x,z\text{ is }y\text{-smooth},z\equiv a \mod r\}},$$
$$\Psi(x,y,\chi)=\sum_{z<x}\mathbf{1}_{\{z'\mid z'\text{ is }y\text{-smooth}\}}(z)\chi(z), \text{ and}$$
$$\rho(x,y)=\Psi(x,y)x^{-1}.$$
\end{defn}
For many the `Holy Grail' in this work is showing the following, \cite{SOUN}:
\begin{conj}\label{4.33}
Let $A$ be a given positive real number. Let $y$ and $q$ be large with $q\leq y^A$. Then as $\frac{\log x}{\log y}\rightarrow \infty$ we have
$$\Psi(x,y,a,q)\sim \frac{1}{\phi(q)}\Psi_q(x,y).$$
\end{conj}
Soundararajan proved in \cite{SOUN} that the above holds for $A\geq 4\sqrt{e}-\epsilon$ given a certain bound on $y$. The latter restriction was later removed by Harper, \cite{SOHA}. We will see more of Harper's work later as we use these results.\\
\\
We now introduce some bounds based on Hildebrand and Tenenbaum's work, which will be presented proof-less here. For those especially interested in this we suggest \cite{HILD},\cite{HIT1},\cite{ILPF}. 
We will avidly be working with these sets later on, but to do so we will first need a few results starting with [\cite{RNFS}, fact 3.11]:
\begin{prop}\label{4.34}
Let $\epsilon>0$ be arbitrary and let $3\leq u\leq\left(1-\epsilon\right)\frac{\log x}{\log\log x}$ Then:
$$\Psi\left(x,x^{\frac{1}{u}}\right)=x\exp\left(-u\left(\log u+\log \log u-1+\frac{\log\log u-1}{\log u}+\boh_{\epsilon}\left(\frac{\log\log^2 u}{\log^2 u}\right)\right)\right).$$
\end{prop}
A very rough result that we use is a direct result of this fact
\begin{corr}\label{4.35}
Fix $0<a<b\leq 1$. Then uniformly in $c,d>0$:
$$\rho(L_x(b,d),L_x(a,c))=L_x\left(b-a,\frac{d(b-a)}{c}\right)^{-1+\loh(1)}.$$
\end{corr}
\begin{proof}
Let 
$$u=\frac{\log L_x(b,d)}{\log L_x(a,c)}=\frac{d}{c}\left(\frac{\log x}{\log\log x}\right)^{b-a}.$$
Then $u\rightarrow\infty$ and $u=\log\left(\frac{\log x}{\log\log x}\right)$. Hence 
$$\rho\left(L_x(b,d),L_x(a,c)\right)=\exp(-(1+\loh(1))u\log  u)$$
$$=\exp\left(-(1+\loh(1))\frac{d(b-a)}{c}\log^{b-a}(x)(\log\log x)^{b-a}\right)$$
$$=L_x\left(b-a,\frac{d(b-a)}{c}\right)^{-1+\loh(1)}.$$
\end{proof}
If we are substantially more careful however there are tighter results, whose proofs go beyond the scope of this paper, by Hildebrand and Tenenbaum. These results allow short intervals to be sieved for smooth numbers effectively:
\begin{thm}\label{4.36}
Fix any $\epsilon>0$. For any $x\geq 3$, $\log x\geq \log y\geq \left(\log\log x\right)^{\frac{5}{3}+\epsilon}$, $z\leq y^{\frac{5}{12}}$, the following estimate holds uniformly:
$$\Psi\left(x\left(1+z^{-1}\right),y\right)-\Psi(x,y)=\frac{\Psi(x,y)}{z}\left(1+\boh\left(\frac{\log(u+1)}{\log y}\right)\right).$$
\end{thm}
\begin{thm}\label{4.37}
For any $x,y$, we set $u=\frac{\log x}{\log y}$. There exists a saddle-point, $\alpha=\alpha(x,y)$, such that for any $1\leq c\leq y$:
$$\Psi(cx,y)=\Psi(x,y)c^{\alpha(x,y)}\left(1+\boh\left(\frac{1}{u}+\frac{\log y}{y}\right)\right), \text{ with}$$
$$\alpha(x,y)=\frac{\log\left(\frac{y}{\log x}+1\right)}{\log y}\left(1+\boh_c\left(\frac{\log\log(y+1)}{\log y}\right)\right).$$
\end{thm}
\begin{thm}\label{4.38}
Let $c>0$ be a constant. Let $n\in N$ such that $\omega(n)$ is the number of prime factors of $n$ (without multiplicity). Let $n$ be a $y$-smooth number with $2\leq y\leq x$ such that $\omega(n)\leq y^{c(\log(1+u))^{-1}}$. Then:
$$\Psi_n(x,y)=\frac{\phi(n)}{n}\Psi(x,y)\left(1+\boh\left(\frac{\log(1+u)\log(1+\omega(n))}{\log y}\right)\right).$$
\end{thm} 
This final theorem induces the following corollary which will prove crucial to us:
\begin{corr}\label{4.39}
Take $c>0$ an arbitrary constant, and retain $\omega$ as in \ref{4.38}. Let $2\leq y\leq x$ and with $\omega(n)\leq y^{c(\log(1+u))^{-1}}$. Then:
$$\Psi_n(x,y)=\frac{\phi(n)}{n}\Psi(x,y)\left(1+\boh_c\left(\frac{\log(1+u)\log(1+\omega(n))}{\log(y)}\right)\right)\left(1+\boh\left(\frac{\omega(n)}{y}\right)\right).$$
\end{corr}
\begin{proof}
Let $n=sr$ for $s$ y-smooth and $r$ with no prime factor less than $y$. Then $\Psi_s(x,y)=\Psi_r(x,y),\phi(n)=\phi(r)\phi(s)$ and, for $p$ prime, $\phi(r)r^{-1}=\prod_{p\mid r}(1-p^{-1})=1+\boh\left(\omega(n)y^{-1}\right)$, which implies the given bound.
\end{proof}

\subsection{Probability measures and moments}\label{sec: PRELprob} 
The Randomized Number Field Sieve is a probabilistic algorithm and makes extensive use of the discrete uniform distribution. Define for $x\in[a,b]$ and for $y\in S$, where $S$ is a finite set, the discrete uniform distribution as
$$\mathbb{P}(x)=\frac{1}{b-a+1}, \text{ or}$$
$$\mathbb{P}(y)=\frac{1}{\abs{S}}$$
This is in many ways one of the simplest distributions to work with and overall we will only need a limited amount of probability theory. There is however a need to understand how a distribution can be used to define a measure. For this recall the following definition:
\begin{defn}\label{4.40}
A set function $\mu:\mathbb{F}\rightarrow\reals$, for a field $\mathbb{F}$, is a probability measure if it satisfies these conditions:
\begin{enumerate}
\item $0\leq \mu(A)\leq 1$ for $A\subseteq \mathbb{F}.$
\item $\mu(\emptyset)=0$, $\mu(\mathbb{F})=1.$
\item $\mu\left(A_1,A_2,\ldots,A_n\right)$, for a disjoint sequence of $\mathbb{F}$-sets, such that $\bigcup_{i=1}^\infty A_i\in\mathbb{F}$, then
$$\mu\left(\bigcup_{i=1}^\infty A_i\right)=\sum_{i=1}^\infty\mu(A_i).$$
\end{enumerate}
\end{defn}
\begin{rmk}
Note that if $\mu$ is a probability measure then the support of $\mu$ is any set $A\subset\mathbb{F}$ for which $\mu(A)=1$. It is clear that this always exists by the second condition.
\end{rmk}
We can confirm that for any finite set $S$ we can define a uniform measure. To see this let $\mu$ be the uniform distribution of some finite set $S=\{s_1,\ldots,s_n\}$ such that $\abs{S}=n$. Then for any $V\subset S$ with order $\abs{V}$ we have that
$$\mu(V)=\frac{\abs{V}}{\abs{S}}.$$
Especially we have that
$$\mu(\emptyset)=0,\mu(S)=1,$$
and 
$$0\leq \mu(V)\leq 1.$$
Lastly for any collection of disjoint subsets $V_1,\ldots,V_m\subset S$:
$$\mu\left(\bigcup_{i=1}^\infty V_i\right)=\mu\left(\bigcup_{i=1}^m V_i\right)=\frac{\sum_{i=1}^m\abs{V_i}}{\abs{S}}=\sum_{i=1}^m\frac{\abs{V_i}}{\abs{S}}=\sum_{i=1}^m\mu\left(V_i\right).$$
Similarly we can confirm that for any continuous interval $[a,b]$ we can define a uniform measure. 
\begin{defn}\label{4.41}
For any finite set $S$, we denote the uniform measure over $S$ by
$$\mathcal{U}(S).$$
\end{defn}
One of the things we are going to use our uniform measure on are the zero-centered half-open integer intervals of length $L$, which we will denote as follows:
$$\mathbb{I}(L)=\left[-\frac{1}{2}L,\frac{1}{2}L\right)\cap\zahle.$$
We note here, for good measure, that this is in fact a finite set. In the section above we already discussed the Dirichlet convolution of arithmetic functions, but there is also a convolution of measures, which is normally defined as an integral, but we are interested in the discrete uniform measures as defined above. This allows us to restate the convolution of measures as follows:
\begin{defn}\label{4.42}
For any two measures $\mu,\nu$ over an additive group $G$ we define the convolution of measures as
$$\left(\mu\star\nu\right)(x)=\sum_{y\in G}\mu(y)\nu(x-y).$$
\end{defn}
\begin{rmk}
Note that it is convention to denote the convolution of measures by $\ast$, but to be distinct we will denote the convolution of measures by $\star$ and the Dirichlet convolution by $\ast$.
\end{rmk}
One of the big concepts in the RNFS is the avoidance of the General Riemann Hypothesis, and for this we need to consider the moments of a probability distribution. We will recognize the first two moments, mean and variance, of the discrete uniform distribution:
\begin{defn}\label{4.43}
Let $U$ be a discrete uniform distribution with support $S$, then the first moment, or mean, is defined as
$$\mathbb{E}(U)=\sum_{s\in S}\mathbb{P}_U(s).$$
\end{defn}
\begin{rmk}
By convention we denote subscripted to the probability, $\mathbb{P}$, the fixed variables. In this case we fix the distribution, but this is usually dropped when the distribution is clear from context.
\end{rmk}
\begin{defn}\label{4.44}
Let $U$ be a discrete uniform distribution with support $S$, then the second moment, or variance, of $x\in S$ is defined as
$$\text{Var}(x)=\mathbb{E}(x-\mathbb{E}(x))^2.$$
\end{defn}
The idea of the RNFS is to remove the dependence of the analysis of number field sieve-type algorithms on the second moment, for which certain bounds on the complexity exist. To do this we consider a probabilistic technique which Lee and Venkatesan, in \cite{RNFS}, have dubbed \textit{stochastic deepening}. The core idea is as follows:
\begin{lemma}\label{4.45}
Let $x$ be a random variable with $\mathbb{E}(x)=\mu$. Given there exists a $K\geq 1$ such that $0\leq x\leq K\mu$ uniformly, then there exists $i\in\{0,\ldots,\lceil\log_2 K\rceil\}$ such that:
$$\mathbb{P}\left(x\geq\frac{2^i\mu}{1+\lceil\log_2 K\rceil}\right)\geq\frac{1}{2^{i+1}}.$$
\end{lemma}
This lemma states that for non-negative random variables which are not too erratic, there is a substantial set where the value is large and whose contribution to the mean is large. This allows us to provide a search algorithm whose run times are shown to be near optimal without establishing variance bounds.\\
\\
The explicit statement of such a search algorithm is not too important for us, as we are analysing the theoretical complexity, so we will not define such an algorithm explicitly. It will however become important in the proof of one of the key theorems, \ref{5.2}.
\newpage
\section{The Randomized Number Field Sieve}\label{ch:RNFS}
\subsection{Introducing the algorithm}\label{sec:RNFSintro}
Recalling the GNFS algorithm, algorithm \ref{alg3}, we will focus in this section on giving a description of the differences between the sieve from chapter \ref{ch: GNFS} and a randomized version with provable complexity. It is well known that the Number Field Sieve's runtime is dominated by two problems:
\begin{enumerate}
\item Finding a set $S$ of sufficient $(a,b)$-pairs such that $a-bm$ is $B$-smooth and $a-b\alpha$ is $B'$-smooth.
\item Collapsing these $(a,b)$-pairs to find a subset $T\subseteq S$ such that a congruence of squares modulo the to-factor integer $n$ arises.
\end{enumerate}
In this version of the Number Field Sieve these two steps are randomized cleverly to produce a provable complexity. To do this there is one more change that needs to be made, and that is the choice of polynomial. It will be shown that by randomizing the choice of polynomial a rigorous bound can be obtained for the search. To start we consider the constants $\delta,\kappa,\sigma,\beta,\beta'$ under the conditions that 
\begin{align}\label{RNFSeq1}
\kappa>\delta^{-1}, 2\sigma>\max\left({\beta,\beta'}\right)+\frac{\delta^{-1}}{3\beta}(1+\loh(1))+\frac{\sigma\delta+\kappa}{3\beta'}(1+\loh(1)),
\end{align}
\begin{align}\label{RNFSeq2}
\delta^{-1}<\frac{\sigma\delta+\kappa}{2},
\end{align}
and fix the smoothness bounds
\begin{align}\label{RNFSeq3}
B=L_n\bigg(\frac{1}{3},\beta\bigg), B'=L_n\bigg(\frac{1}{3},\beta'\bigg).
\end{align}
The culminating theorem of the randomized number field sieve is as follows:
\begin{thm}\label{5.1}
Let $\delta,\kappa,\sigma,\beta,\beta'$ be constants as defined in (1) and (2). Then for any $n$, the randomized number field sieve runs in expected time
$$L_n\bigg(\frac{1}{3},\sqrt[3]{\frac{64}{9}}+\loh(1)\bigg),$$
and produces a pair $x,y$ with $x^2\equiv y^2 \mod n$.
\end{thm}
It is important to stress that this is a NFS-style algorithm, so we can not be sure that the congruence found does not give a trivial factorization. The algorithm is started in the same way as the GNFS by finding a polynomial, however in this randomized version we will insist that the polynomial is a bivariate homogeneous monic irreducible polynomial of bounded degree, that is:
$$f(x,y)\in\zahle[x,y]: f(x,y)=\hat{f}_{m,n}(x,y)+R(x,y),$$
where $\hat{f}_{m,n}(x,y)$ is the polynomial given by $n$ and $m$, such that $\hat{f}_{m,n}(m,1)$ is the base-$m$ expansion of $n$ and 
$$R(x,y)=\sum_{i=0}^{d-1}c_i(x-ym)x^{d-i-1}y^i,$$
such that every $c_i$ is determined uniformly at random and $\deg(f)=d=\delta\sqrt[3]{\frac{\log n}{\log\log n}}$, $d$ odd. Moreover we bound the choice of $m$ by $m^d\leq n<2m^d$. 
\begin{rmk}
The harshness of these bounds are all necessary to show a provable complexity, so we will often refer back to this introduction as ``the defined parameters" without restating the exact parameters again.
\end{rmk}
\begin{rmk}
We will, like \cite{RNFS}, freely switch between considering the polynomial $f(x,y)$ as a bivariate polynomial and it's monovariate equivalent $f(x)=f(x,1)$. 
\end{rmk}
The changes made to the polynomial selection method and the consequences from using random variables lead us to conclude the following theorem:
\begin{thm}\label{5.2}
Take $\delta, \kappa,\sigma, \beta,$ and $\beta'$ with the defined parameters. For any $n$, the randomized number field sieve can almost surely find an irreducible polynomial $f$ of degree $d$ and height at most $L_n(\frac{2}{3})$, with $\alpha$ a root of $f$, $n\mid f(m)$, and 
$$L_n(\frac{1}{3},\max(\beta,\beta')+\loh(1))$$
distinct pairs $a<\abs{b}\leq L_n(\frac{1}{3},\sigma)$ such that $(a-bm)$ is $B$-smooth and $a-b\alpha$ is $B'$ smooth, in expected time at most $L_n(\frac{1}{3},\lambda)$ for any 
$$\lambda>\max(\beta,\beta')+\frac{\delta^{-1}}{3\beta}(1+\loh(1))+\frac{\sigma\delta+\kappa}{3\beta'}(1+\loh(1)).$$
In particular, the probability that the randomized number field sieve fails to produce such a set is bounded above by $L_n(\frac{2}{3},\kappa-\delta^{-1})^{-1+\loh(1)}$.
\end{thm}
The key purpose for the randomization of the polynomial selection process has to do with the sieving process over the number field, as it will cause the norm of a pair, $N(a-b\alpha)$, to become a random variable in $\zahle$. This allows the consideration of smoothness over $\zahle$ and $\zahle[\alpha]$ to be completely independent. We will explain this idea in more detail when we look at the sieving process over the number field. \\
\\
This also leads us to question what the quadratic characters will then look like, and this is the second part where the RNFS differs significantly from the GNFS. Where in the GNFS we defined the quadratic characters as the Legendre symbols $\big(\frac{r}{p}\big)$ modulo a a prime $p$, we now have to account for our randomization and instead will have to choose maps from $\zahle[\alpha]$ into $\mathbb{F}_{p^k}$ stochastically and close to uniformly across all $k\log(p)<L_n(\frac{1}{3})$. This exponentially increases the sizes of the fields, but we will see they are necessary for the unconditional equidistribution of the characters. Once the sieving process is finished we are, with the exception of a set of bad $f$, guaranteed a reduction to a congruence of squares:
\begin{thm}\label{5.3}
Let $B, B'$ be with the defined parameters. Let $f$ be irreducible of degree $d$ and height at most $L_n(\frac{2}{3},\kappa)$, and let $\alpha$ be a root of $f$. Then for all but a $L_n(\frac{2}{3},\frac{\kappa-\delta^{-1}}{2}(1+\loh(1)))^{-1}$ fraction of the set of $f$, if we are given 
$$L_n\bigg(\frac{1}{3},\max(\beta,\beta')\bigg)\bom(\log\log n)$$
pairs $a<b\leq L_n(\frac{1}{3})$ such that $a-mb$ is $B$-smooth and $a-b\alpha$ is $B'$-smooth, we can find a congruence of squares modulo $n$ in expected time at most
$$L_n\bigg(\frac{1}{3},2\max\bigg(\frac{2\delta}{3},\beta,\beta'\bigg)\bigg)^{1+\loh(1)}.$$
\end{thm}
Let's dive into the details.
\subsection{Randomizing the polynomial}\label{sec: RNFSpoly}
As we mentioned the algorithm is a GNFS based algorithm with two major differences: How the polynomial is chosen and how the quadratic characters are chosen. To start we will look at how the polynomial is chosen, and consequences related to that, before proving theorem \ref{5.2}. To do this let $\beta, \beta', \kappa, \sigma, \delta, B$, and $B'$ be as defined in \eqref{RNFSeq1}, \eqref{RNFSeq2}, and \eqref{RNFSeq3}.\\
\\
To begin we make a few restrictions:
\begin{defn}\label{5.4}
Let $\mathcal{X}$ be the set of tuples $(f,m,n,a,b)$ such that the following conditions hold:
\begin{enumerate}
\item $m\in\zahle$ and $m\in\left[2^{-\frac{1}{d}}L_n\left(\frac{2}{3},\delta^{-1}\right), L_n\left(\frac{2}{3},\delta^{-1}\right)\right].$
\item $f\in\zahle[x,y]$, $\deg(f)=d=\delta\sqrt[3]{\frac{\log n}{\log\log n}}$, $2\nmid d$, with coefficients bounded by $L\left(\frac{2}{3},\kappa\right)(1+\loh(1))$. Moreover let the coefficients $\{c_i\}_{i\in I}$ be drawn uniformly at random such that 
$$c_i\in\mathbb{I}\left(2L_n\left(\frac{2}{3},\kappa-\delta^{-1}\right)\right).$$
\item $a,b\in\zahle, 0\leq a<\abs{b}\in\left[\frac{1}{2}L_n\left(\frac{1}{3},\sigma\right)\right]$, with $a-bm$ $B$-smooth and $f(a,b)$ $B'$-smooth.
\end{enumerate}
\end{defn}
Moreover let $\mathcal{X}_{f,m,n}=\{(a,b)\in\zahle^2:(f,m,n,a,b)\in\mathcal{X}\}$. The first thing we note from this definition and the definition of our parameters is that $f(x,y)$ is not generated uniformly at random as it has a  component solely determined by $m$ and $n$, $\hat{f}_{m,n}(x,y)$, which is a problem as we try to obtain that $f(a,b)$ is as likely to be $B'$-smooth as a uniformly random integer of the same size. To mitigate this we shall use the observation that $c_i$ is much larger than $m$ to prove that the randomized part of $f(x,y)$ dominates and that $f(a,b)$ therefore can be considered as a uniformly distributed along any arithmetic progression of common difference $(a-mb)$.\\
\\
Once we have achieved this we shall show that for most $B$-smooth moduli $a-mb$, the $B'$-smooth numbers are approximately uniformly distributed modulo $(a-mb)$. It is only then that we can show that $f(a,b)$ is as likely to be $B'$-smooth as a random integer of the same size. Once we have shown this we shall venture to prove that there are sufficient $(a,b)$-pairs such that $a-mb$ is $B$-smooth and $f(a,b)$ is $B'$-smooth, which will lead to a proof of theorem \ref{5.2}. A key in this is the observation that $f(a,b)$ lies on the arithmetic progression given by
$$\bigg\{\hat{f}_{m,n}(a,b)+(a-mb)z: \abs{z}\leq dL_n\left(\frac{2}{3},\kappa-\delta^{-1}\right)b^d\bigg\}.$$

\subsubsection{Uniform distribution of the polynomial}\label{sec: RNFSunif}
Continuing our exposition using the definitions and parameters we have set up until now we recall the following:
$$\forall i\in\{1,\ldots,d-1\}: \left(c_i\right)\sim\mu= \mathcal{U}\left(\mathbb{I}\left(2L_n\left(\frac{2}{3},\kappa-\delta^{-1}\right)\right)^d\right)$$
and denote $\mathbf{c}=(c_i)$ to be the coefficient vector of $f$.\\
\\
Note that for $f(a,b)$ to be $B'$-smooth it is necessary for $a-bm$ to be $B$-smooth, however this is not something we can simply assume. Therefore we consider the definition of $f(x,y)$ and note that 
$$f(a,b)\equiv \hat{f}_{m,n}(a,b) \mod a-bm.$$
As $\gcd(a,b)^d\mid f(a,b)\hat{f}_{m,n}(a,b)$ we have that $\gcd(a,b)^d(a-bm)\mid R(a,b)$. Now let $a$ and $b$ be uniform in their ranges, then we can give an explicit description of the probability that $a-mb$ is $B$-smooth:
\begin{prop}\label{5.5}
Fix $b$ in its interval. If $a,m$ are uniformly random, then:
$$\mathbb{P}_{a,m}(a-bm\text{ is } B\text{-smooth})=L_n\left(\frac{1}{3},\frac{\delta^{-1}}{3\beta}\left(1+\loh(1)\right)\right)^{-1}.$$
\end{prop}
\begin{proof}
We fix $b$. Note that $a$ is uniformly random on an interval of length $b$, and $m$ is uniformly random over an interval of length comparable to its largest value. In particular:
$$a-bm\sim\mathcal{U}\left[-bL_n\left(\frac{2}{3},\delta^{-1}\right), -b\left(2^{-\frac{1}{2d}}L_n\left(\frac{2}{3},\delta^{-1}\right)+1\right)\right)$$
$$=\mathcal{U}\left[-x\left(1+z^{-1}\right), -x\right),$$
for $x=L_n\left(\frac{2}{3},\delta^{-1}(1+\loh(1))\right)$ and $z\approx 2^{\frac{1}{2d}}-1=\boh(d^{-1})$. Note that $d=\log\left(B^{\frac{5}{12}}\right)$, and that $\log\log B=\boh(\log\log x)$. Hence from theorem \ref{4.36}
the number of smooth values of $a-mb$ is:
$$\frac{\Psi(x,B)}{z}\left(1+\boh\left(\frac{\log(u+1)}{\log B}\right)\right).$$
Since the range of values is of length $\frac{x}{z}$,
$$\mathbb{P}_{a,m}\left(a-bm\text{ is }B\text{-smooth}\right)=\rho(x,B)\left(1+\boh\left(\frac{\log(u+1)}{\log B}\right)\right).$$
Recall that $B=L_n\left(\frac{1}{3},\beta\right)$ and $x=L_n\left(\frac{2}{3},\delta^{-1}\right)^{-1+\log(1)}$. Furthermore, note that $\log u<\log\log n=\loh(\log B)$. Hence by corollary \ref{4.35}
$$\rho\left(L_n\left(\frac{2}{3},\delta^{-1}\right)^{1+\loh(1)},B\right)=L_n\left(\frac{1}{3},\frac{\delta^{-1}}{3\beta}\right)^{-1+_\loh(1)}.$$
Absorbing the multiplicative $1+\loh(1)$ into the $\loh(1)$ error term in the exponent to obtain
$$\mathbb{P}_{a,m}\left(a-bm\text{ is }B\text{-smooth}\right)=L_n\left(\frac{1}{3},\frac{\delta^{-1}}{3\beta}\left(1+\loh(1)\right)\right)^{-1}.$$
\end{proof}
Now we can fix the residue of $f(a,b)\mod a-mb$ and consider only how the randomness of our coefficients $c_i$ affect the polynomial. Therefore we will also want to be explicit about the coefficient vector, $c=(c_i)_{i\in\{1,\ldots,d-1\}}$, for $f(x,y)$, hence we will, when when we need to be explicit, denote $f(x,y)=f_c(x,y)$.\\
\\
Now we are set up to show that $f(a,b)$ can be considered as uniformly distributed along progressions of common difference $(a-mb)$. 
\begin{lemma}\label{5.6}
Let $a<b$, with $\gcd(a,b)=1$, define $\varphi=\varphi_{a,b}:\zahle^d\rightarrow\zahle$ by the following
$$\varphi\left((v_0,\ldots,v_{d-1})\right)=\sum_{i=0}^{d-1} v_ia^{d-i-1}b^i.$$
There exists a set $S\subseteq \mathbb{I}\left(4L_n\left(\frac{1}{3},\sigma\right)\right)^d$ such that $\varphi$ bijects $S$ and $\mathbb{I}\left(b^{d-1}\right)$.
\end{lemma}
\begin{proof}
For each $i\geq 0$, we claim that for any $\abs{t}\leq b^i+a^{i+1}$ there exists a representation
$$t=a^ix_0+a^{i-1}bx_1+\ldots+b^ix_i,$$
with $\abs{x_0},\ldots,\abs{x_i}\leq a+b$. We proceed by induction on the number of terms. If $i=0$ then the conclusion is trivial as $\abs{t}\leq b^0+a^1=a+1$ and so 
$$t=a^0x_0=x_0,$$
for $\abs{x_0}\leq a+b$, and as $a+1\leq a+b$ the result follows. Now let $i>0$, then we may choose $\abs{y}<a$ such that 
$$\abs{t-ya^i}\leq b^i.$$
Letting $z\in \{0,\ldots, b-1\}$ such that $za^i\equiv t-ya^i \mod b$ we can set $x_0=y+z$ and as $\abs{y}<a$ and $\abs{z}<b$ it follows that $\abs{x_0}\leq a+b$. By construction $b\mid t-x_0a^i$ and 
$$\abs{\frac{t-x_0a^i}{b}}=\abs{\frac{(t-ya^i)-za^i}{b}}\leq b^{i-1}+a^i.$$
Hence by the induction hypotheses such a representation exists for every $t$ with $\abs{t}\leq b^i+a^{i+1}$. Now we show the existence of $S$ directly. Let $t\in\mathbb{I}(b^{d-1})$, then for any such $t$: $\abs{t}\leq b^{d-1}$ and so the conditions of the above hold with $i=d-1$. Hence there exists a sequence $(x_i)_{i\in\{0,\ldots,d-1\}}$ such that
$$\sum_{i=0}^{d-1}x_ia^{d-i-1}b^i=t,$$
with $\abs{x_i}< a+b$ for all $i$. Hence the vector given by the sequence, $\mathbf{x}_t=(x_i)$, fulfills that $\varphi(\mathbf{x_t})=t$. Moreover, $\mathbf{x}_t\in\mathbb{I}\left(2(a+b)\right)^d$ and by Defn 2.4.3 this means $\mathbf{x}_t\in\mathbb{I}\left(4L_n(\frac{1}{3},\sigma)\right)^d$. Hence, taking $S=\{\mathbf{x}_t:t\in \mathbb{I}\left(b^{d-1}\right)$ gives the construction of a single vector $\mathbf{x}_t$, hence $\varphi$ is bijective on $S$.
\end{proof}
\begin{rmk}
Note that we have presupposed that $\gcd(a,b)=1$, but it is a simple exercise to show that this is trivial. By the definition of $(a,b)\in\mathcal{X}_{f,m,n}$ we know that $a-mb$ is $B$-smooth, hence $\gcd(a,b)$ is $B$-smooth. Since we will only be interested in looking at $B$-smoothness to show the uniform randomness of $f(x,y)$ we can assume without loss of generality that $\gcd(a,b)=1$ as division of two $B$-smooth numbers preserves $B$-smoothness.
\end{rmk}
We can now reformulate $f(x,y)$ using the map $\varphi$. For this consider $m,n$ fixed and observe that:
$$f(a,b)=f_c(a,b)=f_{c'}(a,b)+(a-mb)\varphi_{(a,b)}(c-c').$$
This motivates the following:
\begin{defn}\label{5.7}
Let $S$ be the set given by the lemma \ref{5.6}
. For any $l\leq b^{d-1}$, we define a set $S_l$ and a measure $\nu_l$ as follows:
$$S_l=\{v\in S\mid\varphi(v)=\mathbb{I}(l)\},$$
$$\nu_l=\mathcal{U}(S_l).$$
Particularly, $\nu_l$ gives a uniformly random element of $S$ whose image under $\varphi$ is in $\mathbb{I}(l)$. 
\end{defn}
It is clear to see that if $v\sim \nu_l$, then we can write 
$$f_v(a,b)\sim f_\mathbf{0}(a,b)+(a-mb)\mathcal{U}(\mathbb{I}(l)).$$
This allows the consideration that for the measures $\nu_l$ with support $S_l$ give us additive alterations to $c$ which change $f(a,b)$ additively by $a-mb$ times a uniformly random value. This means that if $f$ is indeed uniformly random, then we can replace the randomness of $c$ over cosets of $S_l$ with randomness of $f_c(a,b)$ over arithmetic progressions.
\begin{defn}\label{5.8}
For $\bar{\mu}:X\rightarrow\reals^+$ a measure and $F:X\rightarrow Y$ a function, we define the measure $F^{\bar{\mu}}:Y\rightarrow\reals^+$ given by
$$F^{\bar{\mu}}(y)=\bar{\mu}\left(\{F^{-1}(y)\}\right)=\sum_{x:F(x)=y}\bar{\mu}(x).$$
\end{defn}
\begin{rmk}
It will be left to the reader to prove that this is in fact a measure. Moreover note that if $F$ is bijective then $F^{\bar{\mu}}(y)=\bar{\mu}(x)$, for some unique $x\in X$.
\end{rmk}
Now consider $c\sim\nu_l$ such that we can write
$$f_c(a,b)\sim f_{0}(a,b)+(a-mb)\mathcal{U}(\mathbb{I}(l)).$$
As $\mu$ is a uniform distribution on a cube of side $2L_n\left(\frac{2}{3},\kappa-\delta^{-1}\right)$ and we have a bijective function $\varphi:S\rightarrow \mathbb{I}(b^{d-1})$ we can now show that $\varphi^\mu$ is actually close to a convolution of uniform distributions on intervals. In fact, by convolving with several distributions, $\nu_{l_i}$, we can treat each coefficient in $f$ as if it were independent and uniformly random. This is made rigorous in the following proposition:
\begin{prop}\label{5.9}
Fix $a,b$. There is a distribution $\vartheta$ such that $\vartheta$ is the convolution of uniform distribution on intervals of lengths $L_n\left(\frac{2}{3},\kappa-\delta^{-1}\right)a^ib^{d-1-i}$ for $i\in\{0,\ldots,d-1\}$ with
$$\norm{\varphi^\mu-\vartheta}_1=\boh\left(L_n\left(\frac{2}{3}, (\kappa-\delta^{-1})(1+\loh(1))\right)^{-1}\right),$$
and 
$$\abs{\mathbb{E}(\vartheta)}\leq\sum_{i=0}^{d-1}a^ib^{d-i-1}.$$
\end{prop}
\begin{proof}
Recall that the convolution of distributions, given by $\star$, is defined by Definition \ref{4.42}
and define
$$\nu=\mu\star\left[\bigstar_{i=0}^{d-1}\nu_{a^ib^{d-i-1}}\right].$$
From lemma \ref{5.6}, the support of $\nu_{a^ib^{d-i-1}}$ is contained in a cube of side $4L_n\left(\frac{1}{3},\sigma\right)$. Hence the support of $P$ of $\bigstar_{i=0}^{d-1}\nu_{a^ib^{d-i-1}}$ is contained in a cube of side $4dL_n\left(\frac{1}{3},\sigma\right)$. When $\norm{x}_\infty < L_n\left(\frac{2}{3},\kappa-\delta^{-1}\right)-4dL_n\left(\frac{1}{3},\sigma\right)$ and $p\in P$:
$$\mu(x-p)=\mu(x)=\abs{\sup(\mu)}^{-1},$$
so $\nu(x)$ is a convex combination of values in $\{\mu(x-p):p\in P\}=\{\mu(x)\}$. Hence $\nu(x)=\mu(x)$ on the $l_\infty$ ball of radius $L_n\left(\frac{2}{3},\loh(\kappa-\delta^{-1}\right)-4dL_n\left(\frac{1}{3},\sigma\right)$. Then since $L_n\left(\frac{1}{3},\sigma\right)$ is $L_n\left(\frac{2}{3},\loh(\kappa-\delta^{-1})\right)$:
$$\mathbb{P}_{x\sim\mu}(\nu(x)-\mu(x))\geq\left(1-4dL_n\left(\frac{2}{3},\kappa-\delta{-1}\right)^{-1+\loh(1)}\right)^d$$
$$\geq 1-4d^2L_n\left(\frac{2}{3},\kappa-\delta^{-1}\right)^{-1+\loh(1)}=1-L_n\left(\frac{2}{3},\kappa-\delta^{-1}\right)^{-1+\loh(1)}.$$
In particular, we have a bound on the $l_1$ distance between $\mu$ and $\nu$:
$$\norm{\mu-\nu}_1=\sum_{x\in\zahle^d}\abs{\nu(x)-\mu(x)}$$
$$\leq\mathbb{P}_{x\sim\mu}\left(\nu(x)\neq\mu(x)\right)\left(\norm{\mu}_\infty+\norm{\nu}_\infty\right)$$
$$=\boh\left(L_n\left(\frac{2}{3},\left(\kappa-\delta^{-1}\right)\left(1+\loh(1)\right)\right)^{-1}\right).$$
Now for fixed $a,b$ we apply the map $\varphi$ to $\mu$ and $\nu$ to obtain
$$\norm{\varphi^\mu-\varphi^\nu}_1\leq\norm{\mu-\nu}_1\leq\boh\left(L_n\left(\frac{2}{3},(\kappa-\delta^{-1})(1+\loh(1))\right)^{-1}\right),$$
and so the $l_1$ difference of the distributions $\varphi(\mu)$ and $\varphi(\nu)$ on $\zahle$ is small. From definition \ref{5.7} $\varphi^{\nu_l}=\mathcal{U}\left(\mathbb{I}(l)\right)$. Since applying our map $\varphi$ to a measure commutes with convolution of measures:
$$\varphi^\nu=\varphi^\mu\star\left[\bigstar_{i=0}^{d-1}\mathcal{U}\left(\mathbb{I}\left(a^ib^{d-i-1}\right)\right)\right].$$
Since $c_i\sim\mathcal{U}\left(\mathbb{I}\left(L_n\left(\frac{2}{3},\kappa-\delta^{-1}\right)\right)\right)$ are independent random variables,
$$c_ia^ib^{d-i-1}+\mathcal{U}\left(\mathbb{I}\left(a^ib^{d-i-1}\right)\right)$$
is uniformly distributed along an interval of length $L_n\left(\frac{2}{3},\kappa-\delta^{-1}\right)a^ib^{d-i-1}$. Note that $\mathbf{c}\sim\mu$. Hence there is a constant $C$ such that for all $x\in\zahle$:
$$\varphi^\nu(x)\sim\bigstar_{i=0}^{d-1}\left[\mathcal{U}\left(\mathbb{I}\left(L_n\left(\frac{2}{3},\kappa-\delta^{-1}\right)a^ib^{d-i-1}\right)\right)\right](x-C).$$
The shift $C$ accounts for the difference in expectation caused by the fact that the intervals associated with the independent random variables  are not centred. However the centre of each of these intervals has modulus at most $\frac{1}{2}a^ib^{d-i-1}+\frac{1}{2}$, and so $\abs{C}\leq\sum_{i=0}^{d-1}a^ib^{d-i-1}$. We take $\vartheta=\varphi^\nu$ to complete the proof.
\end{proof}
Hence, this convolution allows us to replace $R(a,b)$ by $a-mb$ times some convolution of uniform measures on intervals. This proves that $f(a,b)$ is close to uniformly distributed along arithmetic progressions of common difference $(a-mb)$, so it rests to show that for most $B$-smooth moduli $a-mb$, the $B'$-smooth numbers are almost uniformly distributed modulo $a-mb$. To do this we introduce an estimation of smoothness, which is called goodness:
\begin{defn}\label{5.10}
Fix a smoothness bound $B=L_n\left(\frac{1}{3}\right)$ and $F=L_n(\frac{2}{3})$. Define $\epsilon\left(F,B,r,a\right)=\loh_n(1)$. We then say that a modulus $r$ is $B$-good for $F$ if uniformly over all $(a,r)=1$:
$$\Psi(F,B,r,a)=\left(\frac{\Psi_r(F,B)}{\phi(r)}\right)^{1+\epsilon},$$
where $\phi$ is the Euler phi-function.
\end{defn}
\begin{rmk}
Recall from Definition \ref{4.32}
$$\Psi(x,y,r,a)=\abs{\{z\in\nat:z<x,z\text{ is } y\text{-smooth}, z\equiv a\mod r\}}$$
and
$$\Psi_r(x,y)=\abs{\{z\in\nat:z<x,z\text{ is } y\text{-smooth}, \gcd(z,r)=1\}}.$$
\end{rmk}
We will often supress the $\epsilon$ as we are only interested in an error exponent up to $\loh(1)$. Moreover, we will say $r$ is $B$-bad for $F$ if it is not $B$-good for $F$ and that $r$ is $B$-good (dropping the $F$) if $r$ is $B$-good for all $F\in\left[L_n\left(\frac{2}{3}\right)\omega^{-1},L_n\left(\frac{2}{3}\right)\right]$ for  $\omega=L_n(\frac{1}{3})$.\\
\\
With this notion we can start constraining the behaviour of $f(a,b) \mod a-mb$ for this we will need to characterise the moduli for which the smooth numbers are uniformly distributed across their residue classes. Once we have that it doesn't matter which residue class $f_c(a,b)$ lies in as it does not affect its probability of being smooth as $c$ varies.

\subsubsection{Uniform Distribution of Smooth Numbers}\label{sec: RNFSsmth}
From Section \ref{sec: PRELarit} we have a result of \cite{SOHA} proving that, under some conditions, $y$-smooth numbers $\leq x$ are uniformly distributed over the $\phi(q)$ residue classes $a\mod q$ with $\gcd(a,q)=1$. However the conditions prevent us from concluding uniform distribution directly, hence if we can show that the distribution of $B'$-good numbers is close enough to uniform assuming $a-mb$ is $B$-smooth, to be considered as such then we are done. To do this we consider a Bombieri--Vinogradov style theorem proposed by Harper, \cite{HARP}:
\begin{thm}\label{5.11}
Let $c$ and $K$ be fixed and effective constants. Then or any $log^KF<B<F$ with $u=\frac{\log F}{\log B}$, and $Q\leq\sqrt{\Psi(F,B)}$:
$$\sum_{r\leq Q}\max_{\left(s,t\right)=1}\abs{\Psi(F,B,r,s)=\frac{\Psi_r\left(F,B\right)}{\phi(r)}}$$
$$\ll\Psi\left(F,B\right)\left(e^{\frac{-cu}{\log^2(u)}}+B^{-c}\right)+Q\sqrt{\Psi(F,B)}\log^{7/2}F,$$
with an implied effective constant $C=C(c,k).$
\end{thm}
\begin{rmk}
We remark that the definition of ``effective" differs per application. For our purpose  we may simply claim that $c<K$ are chosen fittingly for the to-factor $n$.
\end{rmk}
If we consider 
$$Q_{max}=\max_{b,m}\abs{a+bm}=L_n\left(\frac{2}{3},\delta^{-1}(1+\loh(1))\right),$$
we can reconsider our question to the equivalent notion of bounding the probability that a $B$-smooth modulus less than $Q_{max}$ is $B'$-bad. This would immediately follow from showing that the number of $B'$-bad moduli below the $Q$-bound is much smaller than $\Psi(Q_{max}, B)$. We will, for our specific needs, bound the sum in theorem \ref{5.11} as we know that the common difference in the arithmetic progression is known to be $y$-smooth. For this let $q=a-mb$, then:
\begin{prop}\label{5.12}
Let $\epsilon>0$ be fixed. Then there exists effective constants $K,c$ such that for any $\log^Kx<y<x^{\frac{1}{\log\log x}}$ with $u=\frac{\log x}{\log y}$, $x^\epsilon\leq Q\leq \sqrt{\Psi(x,y)}$ and $\omega=\lom(1)$ with $\omega=y^{\boh(1)}$:
$$\sum_{\substack{
r\in\left[Q\omega^{-1},Q\right] \\
r \text{ is } y\text{-smooth}}}
\max_{\gcd(a,r)=1}\abs{\Psi(x,y,r,a)-\frac{\Psi_r(x,y)}{\phi(r)}}$$
$$\ll\Psi(x,y)\rho(Q,y)\left(e^{-\frac{cu}{\log^2 u}}+y^{-c}\right)+Q\sqrt{\Psi(x,y)}\log^{7/2}x,$$
for some effective implied constant.
\end{prop}
\begin{rmk}
The proof for this is mostly an exercise in restating Harper's \cite{HARP} and is therefore omitted. The full proof can be found in \cite{RNFS}.
\end{rmk}
Now we return to what we set out to prove: For most $B$-smooth moduli $a-mb$, the $B'$-smooth numbers are uniformly distributed. This follows immediately from the following two propositions:
\begin{prop}\label{5.13}
Fix $a,b,m,n$ in their intervals and let $f$ be uniformly random as before. Then:
$$\mathbb{P}_f\left(f(a,b)\text{ is }B'\text{-smooth}\mid (a-mb)\text{ is }B'\text{-good}\right)=L_n\left(\frac{1}{3},\frac{\kappa+\sigma\delta}{3\beta'}(1+\loh(1))\right)^{-1}.$$
\end{prop}
\begin{proof}
Let $a-bm=r$. By proposition \ref{5.9}:
$$\mathbb{P}_{n,f}\left(f_{n,m}\left(a,b\right)\text{ is }B'\text{-smooth}\right)=\mathbb{P}_{n,c}\left(\hat{f}_{m,n}(a,b)+r\varphi_{(a,b)}(c)\text{ is }B'\text{-smooth}\right)$$
$$=\mathbb{P}_{n,\vartheta}\left(\hat{f}_{m,n}(a,b)+r\vartheta\text{ is }B'\text{-smooth}\right)+\boh\left(L_n\left(\frac{2}{3},(\kappa-\delta^{-1})(1+\loh(1))\right)^{-1}\right).$$
By the definition of the parameters we have $\kappa>\delta^{-1}$ and therefore $\kappa-\delta^{-1}>0$. Now recall, from proposition \ref{5.9}, $\vartheta$ has $\abs{\mathbb{E}(\vartheta)}\leq\sum_{i=0}^{d-1}a^ib^{d-i-1}$ and is sampled according to the convolution of uniform measures on intervals of length $L_n\left(\frac{2}{3},\kappa-\delta^-1\right)a^ib^{d-i-1}$ for $i=0,\ldots,d-1$. Hence $\vartheta$ is unimodal with mode at some $M$ satisfying
$$\abs{M}\leq\sum_{i=0}^{d-1}a^ib^{d-i-1}<db^{d-i-1}=L_n\left(\frac{2}{3},\sigma\delta(1+\loh(1))\right).$$
Now define 
$$F_{\max}=L_n\left(\frac{2}{3},\kappa-\delta^{-1}\right)(a-mb)\sum_{i=0}^{d-1}a^ib^{d-i-1},$$
then the support of $\vartheta$ is contained in $\left[M-F_{\max}\abs{r}^{-1},M+F_{\max}\abs{r}^{-1}\right]$. We choose an $\omega=L_n\left(\frac{2}{3},\loh(1)\right)$, such that $\omega\rightarrow\infty$, and set
$$Y=L_n\left(\frac{e}{3},\kappa-\delta^{-1}\right)b^{d-1}\omega^{-1}=L_n\left(\frac{2}{3},\kappa-\delta^{-1}+\sigma\delta-\loh(1)\right).$$
Now we define a measure $\vartheta'$ to be 
$$\vartheta'(x)=\begin{cases}
\vartheta(\max(x,Y)) & x\geq 0\\
\vartheta(\min(x,-Y)) & x<0
\end{cases}.$$
Then
$$\norm{\vartheta'-\vartheta}_1\leq\mathbb{O}_{z\sim\vartheta}\left(\abs{z}<Y\right)\leq 2Y\left(L_n\left(\frac{2}{3},\kappa-\delta^{-1}\right)b^{d-1}\right)^{-1}=2\omega^{-1},$$
from the definition of $Y$. Note that $Y$ is much larger than $M$ and so $\vartheta'$ is monotone decreasing away from $0$; hence there are non-negative weights $W_y$ for $y\in\zahle$, with $W_y=0$ for $\abs{y}>F_{\max}\abs{r}^{-1}$ such that:
$$\vartheta'=\sum_{y\geq Y}W_y\mathcal{U}(\left[0,y\right])+W_{-y}\mathcal{U}(\left[-y,0\right)),$$
and $\abs{1-\sum_y W_y}\leq 2\omega^{-1}$. Hence we have:
$$\mathbb{P}_f\left(f_{m,n}(a,b)\text{ is }B'\text{-smooth}\right)=\boh\left(\omega\right)^{-1}+$$
$$\sum_{y=Y}^{F_{\max}\abs{r}^{-1}}W_y\mathbb{P}\left(\hat{f}_{m,n}(a,b)+r\mathcal{U}\left(\left[0,y\right]\right)\text{ is }B'\text{-smooth}\right)$$
$$+W_{-y}\mathbb{P}\left(\hat{f}_{m,n}(a,b)+r\mathcal{U}\left(\left[-y,0\right)\right)\text{ is }B'\text{-smooth}\right).$$
We note that $\boh\left(\omega^{-1}\right)=L_n\left(\frac{2}{3},\loh(1)\right)^{-1}$ terms can be absorbed into our $\loh(1)$ terms, and so it suffices to show that for any fixed, $B'$-good $r$ and any $y\in\left[Y,F_{\max}\abs{r}^{-1}\right]$:
$$ \mathbb{P}\left(\hat{f}_{m,n}(a,b)+r\mathcal{U}\left(\left[0,y\right]\right)\text{ is }B'\text{-smooth}\right)=L_n\left(\frac{1}{3},\frac{\kappa+\sigma\delta}{3\beta'}\right)^{-1+\loh(1)},$$
$$ \mathbb{P}\left(\hat{f}_{m,n}(a,b)+r\mathcal{U}\left(\left[-y,0\right)\right)\text{ is }B'\text{-smooth}\right)=L_n\left(\frac{1}{3},\frac{\kappa+\sigma\delta}{3\beta'}\right)^{-1+\loh(1)}.$$
Since $\abs{\hat{f}_{m,n}(a,b)}\leq\hat{F}_{max}:=YL_n\left(\frac{2}{3}\right)^{-1}$, we can absorb the probability that the value on the left is negative or positive (respectively) in the above two equations. From the definition of $B'$-good and corollary \ref{4.39}
, for any $x\in\left[\abs{r}Y-\hat{F}_{\max},F_{\max}+\hat{F}_{\max}\right]$:
$$\Psi(x,B',r,s)=\frac{\Phi_r(x,B')}{\phi(r)}L_n\left(\frac{1}{3},\loh(1)\right)=\frac{\Psi(x,B')}{r}L_n\left(\frac{1}{3},\loh(1)\right),$$
and so to finish the estimate we observe that for any $x\in\left[\abs{r}Y-\hat{F}_{\max},F_{\max}+\hat{F}_{\max}\right]$:
$$\rho(x,B')=\rho\left(L_n\left(\frac{2}{3},\kappa+\sigma\delta\right),B'\right)=L_n\left(\frac{1}{3},\frac{\kappa+\sigma\delta}{3\beta'}\right)^{-1+\loh(1)}.$$
\end{proof}
Now there is only one thing left to do: Show that $a-bm$ is $B'$-good significantly often when $a-bm$ is $B$-smooth, which we have assumed from the start and have a significant probability for from proposition \ref{5.5}.
\begin{prop}\label{5.14}
Fix any $b$. Then
$$\mathbb{P}_{a,m}(a-mb\text{ is } B'\text{-good }\mid a-mb\text{ is } B\text{-smooth})=1-\loh(1).$$
\end{prop}
\begin{proof}
We begin by bounding the number of moduli which are $F$-bad for some $$F\in\left[F_{\max}L_n\left(\frac{1}{3}\right)^{-1},F_{\max}\right].$$ We fix $\omega=B'$ for concreteness. Observe that $\Psi(F,B')=FL_n\left(\frac{1}{3}\right)^{-1}$. Since $L_n\left(\frac{2}{3}\right)=\lom\left(L_n\left(\frac{1}{3}\right)\right).$
$$Q\leq\sqrt{\Psi\left(F,B'\right)}L_n\left(\frac{2}{3},\frac{\epsilon}{4}\right)^{-1}.$$
Furthermore, for any $K$ fixed, $\lom\left(\log^KF\right)=B'=\loh\left(F^{\frac{1}{\log\log F}}\right)$. Hence we can apply proposition \ref{5.12}. Suppose that a modulus $r$ is $B$-smooth and also $B'$-bad for $F$. Then for some residue $a$ with $\gcd\left(a,r\right)=1$, the contribution to the LHS of proposition \ref{5.12} for this $r$ is at least
$$\frac{\Psi(F,B')}{\phi(r)}(1+\loh(1))=\frac{\Psi(F,B')}{r}\geq\frac{\Psi(F,B')}{Q},$$
where for the first equality we use corollary \ref{4.39}, noting that $B\leq B'$ so $r$ is $B'$-smooth, $u<\log\log n$ and the number of divisors of $r$ sis bounded by $\log r$ so that the multiplicative error is $1+\loh(1)$. Now
$$\sum_{\substack{
r\in\left[Q_{\max}\omega^{-1},Q_{max}\right] \\
r \text{ is } y\text{-smooth}}}
\max_{\gcd(a,r)=1}\abs{\Psi(F,B',r,a)-\frac{\Psi_r(F,B')}{\phi(r)}}$$
$$\leq C\Psi(F,B')\rho(Q_{\max},B')\left(e^{-\frac{cu'}{\log^2 u}}+B'^{-c}\right)+Q_{\max}\sqrt{\Psi(F,B')}\log^{7/2}F$$
$$=CF\rho(F,B')\rho(Q_{\max},B')\left(e^{-\frac{cu'}{\log^2 u}}+B'^{-c}\right)+Q_{\max}F^{\frac{1}{2}}\rho(Fc,B')^{\frac{1}{2}}\log^{7/2}F.$$
First, we observe that $F$ and $Q_{max}$ are $L_n\left({2}{3}\right)$, whilst $B'$ is $L_n\left(\frac{1}{3}\right)$. Hence both densities $\rho(Q_{\max},B')$ and $\rho(F,B')$ are $L_n\left(\frac{1}{3}\right)^{-1}$. From the definition of $Q_{\max}$ and $F_{\max}$ we have
$$Q_{\max}=L_n\left(\frac{2}{3},\delta^{-1}(1+\loh(1)\right), F_{\max}=L_n\left(\frac{2}{3},\left(\kappa+\sigma\delta\right)(1+\loh(1)\right),$$
and from the parameters defined in equation (2) we have $2\delta^{-1}<\kappa+\delta\sigma$. Hence $FQ_{\max}^{-2}=L_n\left(\frac{2}{3}\right)$. Since up to order $L_n\left(\frac{2}{3}\right)^{\loh(1)}$, the first term is $F$ and the second is $Q_{max}F^{\frac{1}{2}}$, we deduce that the first term dominates the second. If $r$ is $B'$-bad for $F$ it contributes at least $\Psi(F,B')Q_{\max}^{-1}\left(1+\loh(1)\right)$ to the sum on the left hand side. Hence the number of moduli which are in $\left[Q_{\max}\omega^{-1},Q_{\max}\right]$, are $B$-smooth and $B'$-bad for $F$ is at most
$$\frac{Q_{\max}}{\Psi(F,B')}\Psi(F,B')\rho(Q_{\max},B')\left(C+\loh(1)\right)\left(e^{-\frac{cu'}{\log^2 u'}}+B'^{-c}\right),$$
$$=\Psi\left(Q_{\max},B'\right)(C+\loh(1))\left(e^{-\frac{cu'}{\log^2 u'}}+B'^{-c}\right).$$
If a modulus is $B'$-bad near $F_{\max}$, it must be $B'$-bad for some 
$$F\in\bigg\{F_{\max}L_n\left(\frac{1}{3}\right)^{-1},F_{\max}\bigg\}\cup\bigg\{2^i:2^i\in\left[F_{\max}L_n\left(\frac{1}{3}\right),F_{\max}\right]\bigg\},$$
which is a set of logarithmic size. We can absorb a logarithmic factor into the constants $c,C$, so the number of moduli which are in $\left[Q_{\max}\omega^{-1},Q_{\max}\right]$, are $B$-smooth and $B'$-bad is at most:
$$\Psi(Q_{\max},B')\left(C+\loh(1)\right)\left(e^{-\frac{cu'}{\log^2 u'}}+B'^{-c}\right)=\loh\left(\Psi\left(Q_{\max}B'\right)\right).$$
Hence even assuming that every $B$-smooth number below $Q_{\max}\omega^{-1}$ is $B'$-bad gives:
$$\mathbb{P}_{a,m}\left(a-mb\text{ is }B'\text{-good}\mid a-mb\text{ is }B\text{-smooth}\right)$$
$$\geq 1-\frac{\Psi\left(Q_{\max}\omega^{-1},B'\right)+\Psi\left(Q_{\max},B'\right))}{\Psi\left(Q_{\max},B'\right)}=1-\loh(1).$$
\end{proof}
Hence the probability of a $B$-smooth modulus to be $B'$-good is nearing 1 and the probability for $f(a,b)$ to be $B'$-smooth given that $a-mb$ is a $B$-smooth modulus is comparable to a random integer of similar size to be $B'$-smooth under the same conditions. Hence we can now prove theorem \ref{5.2}. 
\subsubsection{Proof of theorem \ref{5.2}}\label{sec: RNFSth5.2}
Now that we have shown that the changes that randomization causes are controllable we are able to conclude that they are also sufficient to give a provable complexity. For this let us quickly reflect on our choice compared to the general number field sieve. Consider $(f,m,n,a,b)\in\mathcal{X}$ and let $\alpha\in\compl$ be a root of the first coordinate of $f$, $f(\alpha,1)=0$. Then there is a ring homomorphism
$$\zahle[\alpha]\rightarrow\zahle/n\zahle, (a-b\alpha)\mapsto(a-bm)\mod n.$$
From the field norm we also have a multiplicative map from $\mathcal{O}_{\quot(\alpha)}\rightarrow\zahle$, as we would expect from the general number field sieve. Moreover, denoting $f_d$ for the leading coefficient of $f$, the field norm is guaranteed to be in $\frac{1}{f_d}\zahle$ for any element of $\zahle+\alpha\zahle$. Hence the only thing that we can't be sure about is that $f$ is irreducible as $f$ is considered uniformly at random. If $f$ was irreducible then we can apply the strategy from the general number field sieve to obtain a congruence of squares modulo $n$. 
\begin{lemma}\label{5.15}
$$\mathbb{P}(\text{ f is reducible})\leq L_N\left(\frac{2}{3},\frac{\kappa-\delta^{-1}+\loh(1)}{3}\right)^{-1}.$$
\end{lemma}
\begin{proof}
Fix $m,n$ and let $H=2L_n\left(\frac{2}{3},\kappa-\delta^{-1}\right)$ be the range of each coefficient of the random part of our polynomial $f$. Note that if a polynomial over $\zahle$ is reducible if it is reducible modulo every prime. Hence if we bound the number of reducible polynomials modulo $\mathbb{F}_p$ for each prime $p$ and bound how often a polynomial is reducible modulo several primes $p$, we can get good bounds on the number of irreducible polynomials $f$. We count the reducible polynomials $f$  with the Turán sieve, \cite{SINT}. Let $\mathcal{A}= \{\left(c_{d-1},\ldots,c_0\right)\in\zahle^d,\abs{c_i}<H\}$ which we equate with the set of polynomials $f(x,y)=\hat{f}_{m,n}(x,y)+R(x,y)$ with $f$ and $\hat{f}_{m,n}$ are both homogeneous of degree $d$ with $\left(c_i\right)$ the coefficients of $R$. For any prime $r$, let $\mathcal{A}_r$ correspond to the subset of $\mathcal{A}$ corresponding to irreducible polynomials $\mod r$. Note that for any $f$ to correspond to $g \mod r$ we must have $(x-my)\mid \hat{f}_{m,n}-g\in\mathbb{F}_r\left[x,y\right]$ or equivalently $g(m,q)\equiv n\mod r$. We do not insist that $G$ is monic, although any irreducible $g$ must be a scalar multiple of a monic irreducible. To estimate the number of irreducibles we accept the following fact, [\cite{NTFF}, Chapter 2]:\\
\\
For any $0<i<r$, the number of monic irreducible $g$ of degree $D$ such that $g(m)\equiv i\mod r$ is $\frac{r^{D-1}}{D(D-1)}+\boh(r^{D/2})$.\\
\\
Note that for any $g$ over $\mathbb{F}_r$, such that $g(m)\equiv n$ with $r\ll\sqrt{H}$, there are
$$\left(\frac{H}{r}+\boh(1)\right)^d=\left(\frac{H}{r}\right)^d+\boh\left(\left(\frac{H}{r}\right)^{d-1}\right)$$
polynomials lying over $g$ in $\mathcal{A}$, and none if $g(m)\not\equiv n \mod r$. Hence by a union bound over the irreducibles $\mod r$:
$$\abs{\mathcal{A}_r}\leq \frac{H^d}{d(d-1)}+\boh\left(\frac{H^d}{r^{\frac{d}{2}-1}}\right)+\boh\left(H^{d-1}r\right),$$
$$\abs{\mathcal{A}_r\cap\mathcal{A}_{r'}}\leq\frac{H^d}{d^2(d^2-1)}+\boh\left(\frac{H^d}{r^{\frac{d}{2}-1}}\right)+\boh\left(\frac{H^d}{r'^{\frac{d}{2}-1}}\right)+\boh(H^{d-1}rr').$$
From the Turán sieve, considering all primes $r\leq p$ for any $p$ much smaller than $\sqrt{H}$, the number of reducible polynomials $f$ is much smaller than $H^dp^{-1}\log p+H^{d-1}p^2$, which for $p\sim H^{1/3}\log^{1/3}H$ is $H^{d-\frac{1}{3}}\log^{\frac{2}{3}}H$. Then the number of potential $f$ for fixed $n,m$ is $H^d$, and so the probability that $f$ is reducible is at most
$$H^{-\frac{1}{3}}\log^{\frac{2}{3}}H=L_n\left(\frac{2}{3},\frac{\kappa-\delta^{-1}+\loh(1)}{3}\right)^{-1}.$$
\end{proof}
\begin{rmk}
If $f$ is reducible we immediately assume that the algorithm fails, hence by sampling at most $L_n\left(\frac{1}{3}\right)$ polynomials the probability that any of them are reducible is $\loh(1)$.
\end{rmk}
Now assume that $f$ is irreducible. Then the following theorem gives us the last ingredient we need to prove theorem \ref{5.2}:
\begin{thm}\label{5.16}
With $\beta=\beta', \delta,\sigma,\kappa$ chosen subject to \eqref{RNFSeq1}, \eqref{RNFSeq2}, \eqref{RNFSeq3}, and Definition \ref{5.4}:
$$\mathbb{E}_{m,f}\left(\abs{\mathcal{X}_{f,m,n}}\right)\geq L_n\left(\frac{1}{3},\tau\right),$$
with
$$\tau=2\sigma-\frac{\delta^{-1}}{3\beta'}(1+\loh(1))+\frac{\sigma\delta+\kappa}{3\beta}(1+\loh(1)).$$
\end{thm}
\begin{proof}
As proposition \ref{5.5} and proposition \ref{5.12} randomise over $a,m$ for any fixed $b$, and uniformly over $n,f$, for any $b,n,f$:
$$\mathbb{P}_{a,m}\left(a-bm\text{ is }B\text{-smooth and }B'\text{-good}\right)=L_n\left(\frac{1}{3},\frac{\delta^{-1}}{3\beta}(1+\loh(1))\right)^{-1}.$$
Since proposition \ref{5.11} randomizes over $f$ for any fixed $a,b,m$, we have for each fixed $b$ that
$$\mathbb{P}_{a,m,f}(a-bm\text{ is }B\text{-smooth and }B'\text{-good}, f(a,b)\text{ is }B'\text{-smooth})$$
$$=L_n\left(\frac{1}{3},\frac{\delta^{-1}}{3\beta}+\frac{\kappa+\sigma\delta}{3\beta'}\right)^{-1+\loh(1)},$$
as multiplicative factors of $1+\loh(1)$ may be absorbed into the $\loh(1)$ in the exponent of the $L_n\left(\frac{1}{3}\right)$ terms. Summing over the $L_n\left(\frac{1}{3},\sigma\right)^2$ choices for a fixed pair $(a,b)$:
$$\mathbb{E}_{m,f}\left(\abs{\mathcal{X}_{f,m,n}}\right)=\sum_{a,b}\mathbb{P}_{f,m,n}((f,n,m,a,b)\in\mathcal{X})$$
$$=L_n\left(\frac{1}{3},\sigma\right)\sum_b\mathbb{P}_{f,m,n,a}\left(\begin{split}
(a-bm)\text{ is }B\text{-smooth}\\
\wedge f(a,b)\text{ is }B'\text{-smooth}\end{split}\right)$$
$$\geq L_n\left(\frac{1}{3},\sigma\right)\sum_b\mathbb{P}_{f,m,n,a}\left(
\begin{split}
(a-bm)\text{ is }B\text{-smooth}\\
\wedge(a-bm)\text{ is }B'\text{-good}\\
\wedge f(a,b)\text{ is }B'\text{-smooth}\end{split}\right)$$
\begin{align*}
\hspace{1.5cm} & \geq L_n\left(\frac{1}{3},2\sigma-\left(\frac{\delta^-1}{3\beta}\right)(1+\loh(1))+\left(\frac{\sigma\delta+\kappa}{3\beta'}\right)(1+\loh(1))\right).
\end{align*}
\end{proof}
\begin{proof}[\textbf{Proof of theorem \ref{5.2}}]
Let $\tau=2\sigma-\frac{\delta^{-1}}{3\beta'}-\frac{\sigma\delta+\kappa}{3\beta}$, and note that:
$$\lambda\geq\max(\beta,\beta')+\frac{\delta^{-1}(1+\loh(1)}{3\beta}+\frac{(\sigma\delta+\kappa)(1+\loh(1))}{3\beta'}=\max(\beta,\beta')+2\sigma-\tau+\loh(1).$$
For any fixed pair $(m,f)$, we can use the hyper-elliptic curve method to examine any pair $(a,b)$ for a suitable smoothness of $a-mb$ and $f(a,b)$ in $\max(B,B')^{\loh(1)}$ time. Hence we can determine whether a pair $(a,b)$ is in $\mathcal{X}_{f,m,n}$ in time $L_n\left(\frac{1}{3},\loh(1)\right)$.\\
\\
By lemma \ref{5.15} the probability that $f$ is reducible is $L_n\left(\frac{2}{3}\right)$, and we have an unconditional uniform bound $\abs{\mathcal{X}_{f,m,n}}\leq L_n\left(\frac{1}{3},2\sigma\right)$. Hence from theorem \ref{5.16} we deduce that
$$\mathbb{E}_{f,m}\left(\abs{\mathcal{X}_{f,m,n}}\mid f \text{ irreducible}\right)\geq L_n\left(\frac{1}{3},\tau+\loh(1)\right).$$
By using lemma \ref{4.45} we can use stochastic deepening to consider $\abs{\mathcal{X}_{f,m,n}}$ as a random variable of $(f,m)$, with $K\leq L_n\left(\frac{1}{3},2\sigma-\tau\right)$. Hence for some $j\leq 1+\lceil\log_2 K\rceil=\boh\left(\log^{\frac{1}{3}}n(\log\log n)^{\frac{2}{3}}\right)$, we have
$$\mathbb{P}_{f,m}\left(\abs{\mathcal{X}_{f,m,n}}\geq 2^jL_n\left(\frac{1}{3},\tau+\log(1)\right)\right)>2^{-j},$$
absorbing logarithmic terms where appropriate.\\
\\
Now consider the following: To find a collection of pairs the algorithm iterates through each $i\in\{0,\ldots, 1+\lceil\log_2 K\rceil\}$, and for each $i$ generates $2^i$ pairs $(f,m)$, and for each pair $(f,m)$ generates $2^{-i}L_n\left(\frac{1}{3},\max\left(\beta,\beta'\right)+2\sigma-\tau+\loh(1)\right)$ pairs $(a,b)$ and tests for smoothness of $a-mb$ and $f(a,b)$. So if $\abs{\mathcal{X}_{f,m,n}}>2^i L_n\left(\frac{1}{3},\tau+\loh(1)\right)$ then the algorithm finds $L_n\left(\frac{1}{3},\max\left(\beta,\beta'\right)+\loh(1)\right)$ pairs as required. Furthermore if 
$$\mathbb{P}_{f,m}\left(\abs{\mathcal{X}_{f,m,n}}\geq 2^i L_n\left(\frac{1}{3},\tau+\log(1)\right)\right)>2^{-i},$$
then with constant probability at least one of the pairs $(m,f)$ satisfies this condition.\\
\\
Since the time taken for a single $i$ is $L_n\left(\frac{1}{3},\max(\beta,\beta')\right)+2\sigma=\tau+\loh(1)$ we can absorb the logarithmic number of iterations into the $\loh(1)$ term. Hence, iterating it at most a logarithmic number of time reduces the probability of failure to $L_n\left(\frac{2}{3},\kappa-\delta^{-1}\right)$. Hence the expected time taken to complete the algorithm is
$$L_n\left(\frac{1}{3},\max(\beta,\beta')+2\sigma-\tau+\loh(1)\right).$$
\end{proof}
\subsection{Algebraic obstructions Vol. 2: The congruence} \label{sec: RNFScong}
Let $S$ be the set of $(a,b)$ pairs such that $a-mb$ is $B$-smooth and $f(a,b)$ is $B'$-smooth. The previous section provides a way to find a sufficient amount of $(a,b)$ pairs in a sufficiently short time, which we now will assume to be fulfilled. Identically to the General Number Field Sieve algorithm we now start on the next procedure: sieving for the congruence. In Chapter \ref{ch: GNFS} we were interested in finding a subset $S_i$ such that 
$$f'(m)^2\prod_{\left(a,b\right)\in S_i}a-bm\text{ is a square in }\zahle,$$
$$f'(\alpha)^2\prod_{\left(a,b\right)\in S_i} a-b\alpha\text{ is a square in }\zahle[\alpha].$$
In the RNFS there will be a similar approach, but it can be recalled that we used particularly broad observations to avoid having to work with the ring of algebraic integers $\mathcal{O}_{\quot[\alpha]}$ for which the structure might be unknown or especially difficult. In this section we will not only tackle the procedure to find these squares, but in particular how we can avoid similar algebraic obstructions as we saw in section \ref{sec:GNFSobstr}.\\
\\
One of the first things to observe is that, given a set $S$ such that for every $(a,b)$-pair $a-bm$ is $B$-smooth is that there is no reason for us not to use the standard method from the general number field sieve to find a square. Given an element $z\in\zahle$ with prime factorization $\prod_{i=1}^d p_i^{\eta_i}$ it is easy to check if it is a square by ensuring that $\forall i\in\{1,\ldots,d\}:2\mid\eta_i$. So given the factorizations of $a-mb$ for $B+1$ $(a,b)$-pairs we can find a dependent subset in the same way we did with the general number field sieve in Step 3.\\
\\
Once again it would be optimal if we could approach the $\zahle[\alpha]$ problem in the same way, but as we saw in Chapter 3 this requires some justification. This justification is only made more complex as we are now dealing with a uniformly random polynomial. We start in much the same fashion as we did with the general number field sieve by defining a set of pairs $(p,r)$ such that $f(r,1)\equiv 0 \mod p$ for prime $p$ and $r\in\{1,\ldots,p-1\}$. Recall that the pairs $(p,r)$ such that $p$ is prime, $0<r<p$ coprime to $p$, $p\mid f(r,1)$ is in direct correspondence with the first degree prime ideals $\PriID\subseteq\mathcal{O}_{\quot(\alpha)}$ such that $\PriID\mid \langle p\rangle, N(\PriID)=p$. Then from proposition \ref{3.8} we get that the following function is well defined:
$$e_{(p,r)}(a-b\alpha)=\begin{cases}
\ord_r\left(N(a-b\alpha)\right) & a-br\equiv 0\mod p\\
0 & otherwise
\end{cases}.$$
From the proof of theorem \ref{3.7} we know that this map is well defined from $\quot(\alpha)^\times\rightarrow\zahle$. Hence, as we assumed that we have found sufficient $(a,b)$-pairs we can use the same sieving techniques as in Section \ref{ch: GNFS} to find a subset $S$ such that $P(\alpha)=\prod_{(a,b)\in S}(a-b\alpha)\in\zahle[\alpha]$. Moreover we can again be sure that $e_{r,s}(P(\alpha))\equiv 0\mod 2$. Now it rests to us to show that we can, again by the use of characters, guarantee with high probability that the element found by $P(\alpha)$ is a square in $\zahle[\alpha]$. To give a three line summary of our technique: we will apply the pigeonhole principle to a set $H$, to show that for a randomized field with a stochastic collection of characters with large conductor the number of ways in which an element might appear square and not be square is limited.
\subsubsection{Quadratic characters}\label{sec: RNFSquad}
Now define the set we eluded to above:
$$H=\{z\in\quot(\alpha)^\times:\forall s<r, e_{(p,r)}(z)\equiv 0 \mod 2\}/\{z^2:z\in \quot(\alpha)^\times\}.$$
Because we attempt to use the pigeon hole principle we start by considering the size of this set, in fact a near-immediate result is the bounded dimension of $H$ as a vector space:
\begin{prop}\label{5.17}
$H$ is an $\mathbb{F}_2$ vector space of dimension at most
$$\left(\delta\kappa+\loh(1)\right)\log_2(n)+\frac{\delta^2\kappa}{2\log(2)}\frac{\left(\log(n)\right)^{\frac{4}{3}}}{\left(\log\log(n)\right)^{\frac{1}{3}}}.$$
\end{prop}
The proof of this follows, bar a few modifications, from the argument in [\cite{BLP}, theorem 6.7]. The full proof can be found in [\cite{RNFS}, p. 129-131], but we will accept this as a fact. To see that $H$ is a vector space we only need to observe that $\quot(\alpha)^\times$ is commutative, hence every element of $h\in H$ can be represented as a coset $[h_0]=h_0\cdot\{z^2:z\in\quot(\alpha)^\times\}$. This means that for any $h\in H$, $h^2$ is the identity since it is equivalent to $h_0^2\cdot\{z^2:z\in (\quot(\alpha)^\times\}$ and $h_0\in\quot(\alpha)^\times$.\\
\\
The next step is to find a significant amount of quadratic characters. For this let $\PriID\subset\mathcal{O}_{\quot(\alpha)}$ be a prime ideal, as $N(\PriID)=p^k$ and $\mathcal{O}_{\quot(\alpha)}$ is a Dedekind Domain we have that 
$$\mathcal{O}_{\quot(\alpha)}/\PriID=\mathbb{F}_{p^k}.$$
Hence we can identify $\PriID$ with a degree $k$ monic irreducible polynomial $p_\PriID$ over $\mathbb{F}_r$. We can say more: Given an irreducible polynomial $g$ of degree $k$ over $\mathbb{F}_p$, if $\gcd(f,g)=1$ then the quotient $\quot(\alpha)/\langle g(\alpha)\rangle$ sends everything to $0$, hence $\langle g(\alpha)\rangle$ is not prime. So we may assume $\gcd(f,g)>1$, and since $g$ is irreducible over $\mathbb{F}_p$ that means that $g$ is one of the irreducible factors of $f\mod \langle p\rangle$. In what follows we will regularly shift between following three equivalent concepts:
\begin{enumerate}
\item Prime ideal $\PriID\subseteq \mathcal{O}_{\quot(\alpha)}$.
\item The irreducible polynomial divisor $p_\PriID$ of $f(x,1) \mod p$.
\item Pair $(p,r)$, $p$ prime and $r$ s.t. $p_\PriID(r,1)\equiv 0 \mod p$ for $p_\PriID$ defined over $\mathbb{F}_{p^k}$.
\end{enumerate}
The following definition extends the Legendre symbol. 
\begin{defn}\label{5.18}
Let $K$ be a field and let $f(x)\in K[x]$ be a polynomial. Let $p(x)\in K[x]$ be an irreducible polynomial, $p(x)\nmid f(x)$. Then the polynomial Legendre symbol is defined as follows:
$$\left(\frac{f(x)}{p(x)}\right)_{L}=\begin{cases}
1, & \exists g(x)\in K[x] \text{ s.t. } f(x)\equiv g(x)^2 \mod p(x)\\
-1,& \nexists g(x)\in K[x] \text{ s.t. } f(x)\equiv g(x)^2 \mod p(x)
\end{cases}.$$
\end{defn}
This gives a very natural extension to a Dirichlet character, $\chi_\PriID$, defined by the following:
\begin{align}\label{char}
\chi_\PriID:\mathcal{O}_{\quot(\alpha)}\rightarrow \mathbb{F}_{p^k}, \quad \langle a-b\alpha\rangle \mapsto \left(\frac{a-bx}{p_\PriID}\right)_L
\end{align}
All that remains now is finding $\PriID$, for which we will seek to factorize $f(x,1)\mod p$ and look at the irreducible divisors. Given a set $\mathcal{F}$ of these characters $\chi_\PriID=\chi_{(p,r)}$, we can define the map $\Psi_{\mathcal{F}}:H\rightarrow\mathbb{F}_2^{\abs{\mathcal{F}}}$ \footnote{Implicitly we are assuming $\mathbb{F}_2\cong C_2$ by the map $b\mapsto \left(-1\right)^b$}
 sending $x$ to the tuple $\left(\chi_{(p,r)}\mid\chi_{(p,r)}\in\mathcal{F}\right)$. We will rely on the random production of such a set $\mathcal{F}$ to show that this almost surely makes $\ker(\Psi_\mathcal{F})$ small.
\begin{lemma}\label{5.19}
There is a sampleable distribution $\Upsilon$ for pairs $(p,r)$ such that $\chi_{(p,r)}$ is a character as in \ref{char}, such that for all but $\log\log n$ of the $h\in H$, considering map as a map from $H$ to $\mathbb{F}_2$:
$$\mathbb{P}_\Upsilon\left(\chi_{(p,r)}=-1\right)\geq\frac{1+\loh(1)}{2}.$$
Sampling according to $\Upsilon$ takes at most $L_n\left(\frac{1}{3},c\right)$ time for a to-define $c$. Furthermore, each character, $\chi_{(p,r)}$, can be evaluated in time at most $L_n\left(\frac{1}{3},\frac{c}{2}\right)$.
\end{lemma}
\subsubsection{Proof of lemma \ref{5.19}}\label{sec: RNFSle5.19}
To start assume that we have the finite-degree tower of number fields $L\supset K\supset\quot$, where $L/K$ is a Galois extension with Galois group $\text{Gal}(L/K)$. Recall that $\Delta_L$, $\Delta_K$ are the discriminants of $L$ and $K$ respectively, and let $d_{L/K}=\left[L:K\right]$, $d_L=\left[L:\quot\right]$ and $d_K=\left[K:\quot\right]$. Let $\PriID\subseteq K$ be a prime ideal and let $Q=\{\mathfrak{Q}\subseteq L\mid \mathfrak{Q}\text{ lies over }\PriID\}$, then from proposition \ref{2.30} we know that the Frobenius elements of all $\mathfrak{Q}\in Q$ are conjugate, which allows us to give the following definition:
\begin{defn}\label{5.20}
Let $L\supset K\supset\quot$ be a finite degree tower of number fields with $L/K$ Galois. Let $\PriID\subset K$ be a prime ideal which is unramified in $L$ and let $\text{Frob}_\PriID$ be the Frobenius element of $\PriID$ in $\text{Gal}(L/K)$. Then the Artin symbol, $\left[\frac{L/K}{\PriID}\right]$, is defined as the conjugacy class of the Frobenius automorphisms of $L/K$ corresponding to the prime ideals $\mathfrak{Q}\subset L$ dividing $\PriID$:
$$\left[\frac{L/K}{\PriID}\right]=\{\text{Frob}_\mathfrak{Q}\in\text{Gal}(L/K)\mid\hspace{5px} \mathfrak{Q}\mid\PriID\}.$$
\end{defn}
\begin{rmk}
It is a common abuse of notation to write $\text{Frob}_\PriID$ for any of the elements in $\left[\frac{L/K}{\PriID}\right]$ and just consider the element well defined up to conjugacy.
\end{rmk}
This allows us to define the set
$$\pi_C(x)=\abs{\bigg\{\PriID\mid \PriID\subset K\text{ prime}, N_{K/\quot}(\PriID)<x, \left[\frac{L/K}{\PriID}\right]\in C\bigg\}},$$ 
where $C\subseteq\text{Gal}(L/K)$ s.t. for all $g\in\text{Gal}(L/K)$: $gCg^{-1}=C$. Now let $K=\quot(\alpha)$ and let $h\in\mathcal{O}_K$ be an element of minimal norm representing a non-trivial element of $H$. Let $L=K(\sqrt{h})$, then by definition $[L:K]=2$, $d_K=d$, $d_L=2d$ and $\text{Gal}(L/K)=C_2$. Since  $\text{Gal}(L/K)$ is abelian, we have that $\left[\frac{L/K}{\PriID}\right]$ only contains one element. Hence the value $\left[\frac{L/K}{\PriID}\right]$ corresponds exactly to the action of $\chi_\PriID$ on $h$.\\
\\
As $h$ has minimal norm it fulfills the Minkowski bound, \ref{2.41}, hence:
$$N_{K/\quot}(h)\leq c_{K/\quot}=\sqrt{\Delta_K}\left(\frac{4}{n}\right)^{\frac{d}{2}}\frac{d!}{d^d}=n^{\delta\kappa(1+\loh(1)}.$$
By construction the different is generated by $2h$, and so is an integral ideal. As the relative discriminant is the norm of the different we get that:
$$\Delta_{L/\quot}\leq N_{K/\quot}(2h)\Delta_{K/\quot}^2\leq n^{\delta\kappa(5+\loh(1))}.$$
We now need an improved version of Chebotarev's Density Theorem which we alluded to before. 
\begin{thm}\textbf{Unconditional Effective Chebotarev Density Th.}
Let $L/K/\quot$ be a tower of number fields with $L/K$ Galois with $G=\text{Gal}(L/K)$. Let $C\subseteq G$ be a union of conjugacy classes such that $gCg^{-1}=C$ for all $g\in G$. Let $\abs{\bar{C}}$ be the total number of conjugacy classes in $G$. Lastly let $1-\nu$ be a Siegel zero of $\zeta_L$ if it exists and $0$ otherwise. Then there exists $c>0$ such that if $\log x\geq 10d_L\log^2( \Delta_L)$:
$$\abs{\pi_{G'}(x)-\frac{\abs{C}}{\abs{G}}\text{Li}(x)}\leq\frac{\abs{C}}{\abs{G}}\text{Li}(x^{1-\nu})+\boh\left(x\abs{\bar{C}}\exp(-{c_1}\sqrt{\frac{\log x}{d_L}}\right).$$
\end{thm}
Using the unconditional effective Chebotarev Density Theorem on the degree 2 extension $L/K$, then for $\PriID$ chosen uniformly at random with $N(\PriID)\leq x$:
$$\abs{\mathbb{P}\left(\chi_\PriID(h)=1\right)-\frac{1}{2}}<x^{-\nu(1+\loh(1))}+\boh\left(2\log x\exp\left(-c_1\sqrt{\frac{\log x}{d_L}}\right)\right),$$
where $1-\nu$ is a possible Siegel zero of the Dedekind zeta function over $L$, $\zeta_L$.  Ensuring that $\mathbb{P}(\chi_\PriID(h)=1)=\frac{1}{2}+\log(1)$ requires us to insist that $\log(x)=\lom(d_L(\log\log x)^2)$ and, if $\zeta_L$ has a Siegel zero, $\log(x)=\lom(\nu^{-1})$. \\
\\
Note that, for each non-trivial coset, we can choose a representative $h\in\mathcal{O}_K$ of minimal norm and then run through the previous paragraph for each of these. Despite the element $h$ being different every time we can get a similar result each time, dependent on $h$. Hence we make the following definition:
\begin{defn}\label{5.21}
For a field $K$ and $h$ a minimal norm representative of an element of $H$, we define $L_h=K(\sqrt{h})$. For $\epsilon>0$ we define the exceptional set:
$$E_{K,\epsilon}=\bigg\{h\cdot\{z^2\mid z\in K^\times\}\in H\text{ s.t. } \exists\nu:\zeta_{L_h}(1-\nu)=0, \nu^{-1}>L_n\left(\frac{1}{3},\epsilon\right)\bigg\}.$$
\end{defn}
It is clear that if $h\in K^\times$ and $k\in K^\times$ represent the same coset in $H$, then $L_h=L_k$. To justify this we only have to note that representing the same coset means that $h$ and $k$ differ by a square:
$$h=kz^2,$$
hence $\sqrt{h}=\sqrt{kz^2}=\sqrt{k}z$, since $z\in K^\times$ we get that $\text{irr}_K(\sqrt{h})=\text{irr}_K(\sqrt{k})$ and so $L_h=L_k$.\\
\\
The exceptional set $E_{K,\epsilon}$ is the subset of $H$ which cannot be reliably distinguished from $0$ by characters induced by primes of size $e^{L_n\left(\frac{1}{3},\epsilon\right)}$. This means that if a Siegel zero of the form specified exists, then almost every prime ideal of this size induces a character which vanishes for some element of $H$. However we can limit the size of the exceptional set.
\begin{lemma}\label{5.22}
Suppose that $K=\quot(\alpha)$ is a number field where $\alpha$ is a root of an irreducible $f=\hat{f}+(x-m)R$ where $R$ is uniformly random. Then for $\epsilon=\left(\frac{1}{3}+\loh(1)\right)\delta$,
$$\mathbb{P}_f\left(\abs{E_{K,\epsilon}}>\frac{4}{3}\log\log n\right)\leq L_n\left(\frac{2}{3},\frac{\kappa-\delta^-1}{3}(1+\loh(1)\right)^{-1} .$$
\end{lemma}
For now we accept this lemma as fact, as this will allow us to prove lemma \ref{5.19}. The proof for this lemma is based on the sparseness of Siegel zeros for Dedekind zeta functions defined over $L_h$ and will be shown in Section \ref{sec: RNFSle5.22}. Accepting the lemma, we know there exists an $x$ such that $\log(x)=\lom(d_L(\log\log x)^2)$ and, if there is a Siegel zero for $\zeta_{L_h}$, $\log(x)=\lom(\nu^{-1})$ for all but $\frac{4}{3}\log\log n$ of the $h\in H$. Moreover for all but at most $L_n\left(\frac{2}{3}\right)^{-1}$ of our maximum $L_n\left(\frac{1}{3}\right)$ polynomials $f$ we have $\log x<L_n\left(\frac{1}{3},\epsilon\right)$. Choosing to acccept failure on all these polynomials still guarantees that the algorithm fails with probability $\loh(1)$.\\
\\
Now we turn to the sampleable distribution $\Upsilon$. Recall that any prime for which $N(\PriID)<x$ is guaranteed to divide a prime $p$ with $p<x$. Moreover, if $\PriID$ is of degree $k$, then $p<\sqrt[k]{x}$. Equivalently, a $k$th degree prime ideal dividing $p$ corresponds to a simple $k$'th degree divisor of $f \mod p$. The following algorithm samples $\Upsilon$ in such a way that it outputs ideals of which all but a small fraction are prime.
\begin{algorithm}[H]\label{RNFSalg1}
\caption{Sampling $\Upsilon$ for (prime) ideals.}
\begin{algorithmic}
\Input{A polynomial $f$ of degree $d$}
\Output{A pair $(p,r)$}
\State Uniformly at random sample $k\in\zahle/d\zahle$
\State Uniformly at random sample $p\in\left(x^{(k+1)^{-1}},x^{k^{-1}}\right)$
\State Test for primality of $p$ with the Miller--Rabin test
\If{r is prime}
\State factor $f\mod p$ \textbf{and} find the irreducible and unrepeated factors $r_i$ of degree $<k$
\If{$\abs{\{r_i\}}\geq j$}
\State Uniformly at random sample one $r_i$ \textbf{and} define $r_i=r$
\Else { Sample a new $k$ and repeat the algorithm}
\EndIf
\EndIf
\end{algorithmic}
\end{algorithm}
\begin{rmk}
We will now make a few remarks regarding the runtime of the algorithm. First note that to obtain the required runtime of the sieve we need this algorithm to have runtime at most $\boh\left(\log^4(n)\right)$. Let's see how we obtain this:
\begin{itemize}
\item First note that the runtime-bound immediately means $\log x=L\left(\frac{1}{3}\right)$. This means we can not use deterministic primality tests, such as the AKS primality test we used for the GNFS. Especially we note that discarding a composite $p$ using the Miller--Rabin primality test takes $\boh(\log^3(x)\log\log x)$. 
\item The Miller--Rabin primality test discards composite $p$ with probability $1-\boh(\log^{-1}(x))$ and so $p$ is prime with probability $\bom\left((\log x)^{-1}\right)$ and so any $p$ produced here is prime with probability $1-\boh\left(d(\log x)^{-1}\right)$.
\item Factoring $f\mod p$ is possible in time $\boh\left(\left(d\log x\right)^3\right)$ probabilistically.
\item If $p$ is composite, but slips through the cracks of the Miller--Rabin primality test, then the factorisation of $f\mod p$ may fail. If it does not fail, then this $p$ still induced a quadratic character and therefore vanishes on the squares. Since we obtain at most one character from each $p$ and are guaranteed to find a character if $k=d$ and $p$ is prime, the fraction of characters which are not induced by primes is $\loh(1)$.
\end{itemize}
\end{rmk}
To finish the proof we need to show that this algorithm is sufficiently fast and that characters can be evaluated quickly. Moreover we, once again, have to show that these characters are sufficiently uniformly distributed that the bounds on $\mathbb{P}(\chi_\PriID=1)$ hold. 
\begin{prop}\label{5.23}
The expected time taken to sample $(p,r)\sim\Upsilon$ as above is at most $$L_n\left(\frac{1}{3},\left(4+\loh(1)\right)\epsilon\right).$$
\end{prop}
\begin{proof}
We noted in our remark that each attempt from the start of the algorithm takes time $\boh(\left(d\log x)^3\right)$ as the factorization of $f\mod r$ is slowest. We are guaranteed to find a factor if our degree bound $k$ is $d$, which happens with probability $\frac{1}{d}$, the integer $r$ is prime with probability $\bom(\frac{1}{\log x})$, and we successfully take an ideal in the final step with probability $\bom\left(\frac{1}{d}\right)$ if $k=d$. Hence the number of attempts needed to output a prime is bounded in expectation by $\boh\left(d^2\log x\right)$. Hence the time taken to find an ideal is bounded in expectation by $\boh(d^5\log^4(x))=L_n\left(\frac{1}{3},(4+\loh(1))\epsilon\right)$.
\end{proof}
\begin{prop}\label{5.24}
For any fixed $h$, $\mathbb{P}_\Upsilon(\chi_{(p,r)}(h)=-1)\geq\frac{1}{2d}(1+\loh(1))$.
\end{prop}
\begin{proof}
The distribution of primes $\PriID$ generated is uniform over $\PriID\mod\langle p\rangle$ for $p\in\left(x^{(k+1)^{-1}},x^{k^{-1}}\right)$ of degree at most $k$ by the Prime Ideal Theorem, \ref{4.22}. This property also holds for a uniform distribution over primes of norm $\leq x$ by an adaptation of Dirichlet's theorem of primes in arithmetic progressions, \ref{4.28}. Thus the difference between $\Upsilon$ and a uniform distribution over primes of norm $\leq x$ is the distribution of the degree of these primes. The probability that $\Upsilon$ samples $\PriID$ with $N(\PriID)\leq x$ and $\PriID\mid \langle p \rangle$ for $p$ is each of these intervals is $\frac{1}{d}$. Hence $\Upsilon$ pointwise dominates $d^{-1}$ times the uniform distribution over all primes of norm below $x$. Therefore
$$\mathbb{P}_\Upsilon(\chi_{(p,r)}(h)=-1)\geq\frac{1}{d}\mathbb{P}_{N(\PriID)\leq x}(\chi_\PriID(h)=-1)=\frac{1}{2d}(1+\loh(1)).$$
\end{proof}
The final proposition immediately completes the proof of lemma \ref{5.19}, however it is quite technical as the logarithms we are interested in are quite large. We, therefore, have to track the arithmetic carefully.
\begin{prop}\label{5.25}
Evaluating the characters $\chi_{(p,r)}$ associated with the ideal $\PriID$ sampled as above on a term $a-b\alpha$ takes time at most $L_n\left(\frac{1}{3},\left(2+\loh(1)\right)\epsilon\right)$.
\end{prop}
\begin{proof}
Let $p=2$, then $\chi_{(p,r)}=1$. Hence we may assume $p>2$. For any polynomial $P\in\mathbb{F}_p[x]$ let $\abs{P}=p^{\deg(P)}$. Recall
$$\chi_{(p,r)}(a-b\alpha)=\left(\frac{a-bx}{r(x)}\right)_L.$$
Any constant $c$ can reduce the computation to finding a Legendre symbol mod $p$:
$$\left(\frac{c}{P}\right)_L=c^{\frac{\abs{P}-1}{2}}=\left(\frac{c}{p}\right)^{\frac{p^k-1}{r-1}}=\left(\frac{c}{p}\right)^k.$$
From \cite{NTFF} we call attention to the law of quadratic reciprocity for function fields by noting that for any two relatively prime monic irreducible polynomials over $\mathbb{F}_p$:
$$\left(\frac{P}{Q}\right)\left(\frac{Q}{P}\right)=\left(-1\right)^{\frac{\abs{P}-1}{2}\frac{\abs{Q}-1}{2}}.$$
Hence
$$\chi_{(p,r)}(a-b\alpha)=\left(\frac{-b}{p}\right)^k\left(\frac{x-ab^{-1}}{r(x)}\right)_L=\left(\frac{-b}{p}\right)^k\left(-1\right)^{\frac{p-1}{2}\frac{p^{\deg(r)}-1}{2}}\left(\frac{r(x)}{x-ab^{-1}}\right)$$
$$=\left(\frac{-b}{p}\right)^k\left(-1\right)^{\frac{p-1}{2}\frac{p^{\deg(r)}-1}{2}}\left(\frac{r(ab^{-1})}{x-ab^{-1}}\right)=\left(\frac{-b}{p}\right)^k\left(-1\right)^{\frac{p-1}{2}\frac{p^{\deg(r)}-1}{2}}\left(\frac{r(ab^{-1})}{p}\right).$$
The parities of $\frac{p-1}{2}$ and $\frac{p^k-1}{2}$ are easily computed. Hence to compute $\chi_{(p,r)}(a-b\alpha)$ it suffices to compute $r(ab^{-1})$ and two legendre symbols modulo $p$ in $\boh(\log p)$ additions or subtractions of numbers of size at most $p$.\\
To compute $b^{-1} \mod p$ requires the extended Euclidean algorithm to be run, which requires $\boh(\log p)$  additions of numbers of size at most $p$.\\
To compute $ab^{-1}\mod p$ requires one multiplication.\\
To compute $s\left(ab^{-1}\right) \mod p$ requires at most $\boh(d)$ addition and multiplications modulo $p$.\\
Addition or subtraction of numbers of size $p$ or modulo $p$ takes $\boh(\log p)$ steps, while multiplication modulo $p$ takes $\boh\left(\log^2(p)\right)$ steps by iterative addition and doubling. Hence the computation time in total requires time
$$\boh(d\log^2(p))=\boh(d^{-1}\log^2(x))=L_n\left(\frac{1}{3},\left(2+\loh(1)\right)\epsilon\right).$$
\end{proof}
To complete the proof of the lemma we now simply take $c=\left(4+\loh(1)\right)\epsilon.$
\subsubsection{Proof of lemma \ref{5.22}}\label{sec: RNFSle5.22}
As we remarked earlier this proof is dependent on the sparseness of Siegel zeroes over Dedekind zeta functions. To start this we need the following proposition, [\cite{RNFS}, fact 6.22 and corollary 6.23], which we will state without proof.
\begin{prop}\label{5.26}
Let $K$ be a number field and let $\overline{K}$ be the normal closure of $K$. Let $c(K)=4$ if $K/\quot$ is normal and $c(K)=4[\overline{K}:\quot]$ otherwise. Assume that there is a real Siegel zero of $\zeta_K$, denoted $1-\nu$, such that
$$1-\left(c(K)\log\abs{\Delta_K}\right)^{-1}\leq 1-\nu\leq 1.$$
Then there is a quadratic field $F\subset K$ such that $\zeta_F(1-\nu)=0$.
\end{prop}
Now we start the proof of lemma \ref{5.22}. To do this we may assume that $f(x)\in \quot(\alpha)[x]$ is an irreducible polynomial, as the probability that $f$ is reducible can be absorbed in the error term by \ref{5.15}
\begin{prop}\label{5.27}
If $\overline{L}_h\supset L_h$ is the normal closure of $L_h=K\left(\sqrt{h}\right)$, then $\left[\overline{L}_h:\quot\right]\leq 2^dd!$.
\end{prop}
\begin{proof}
Let $\overline{K}$ be the splitting field of $K$. By construction, $\overline{K}/\quot$ is normal, especially that means it is Galois, and $[\overline{K}:\quot]\leq d!$. Let $\overline{G}=\Gal(\overline{K}/\quot)$. Given $h\in\mathcal{O}_K$, let $\mathcal{O}_h$ be the orbit of $h$ under the action of $\overline{G}$. Then $\abs{\mathcal{O}_h}\leq d$. We adjoin square roots of each element of $\mathcal{O}_h$ to $\overline{K}$ to obtain a field $\overline{L}$. Then $[\overline{L}:\quot]\leq 2^dd!$.\\
\\
Since degree 2 extensions are normal, the compositum of normal extensions is normal, and $\overline{L}/\overline{K}$ is a compositum of at most $d$ degree $2$ extensions, hence the extension $\overline{L}/\overline{K}$ and $\overline{K}/\quot$ are normal. As any $\sigma\in \text{Aut}_\quot(\overline{K})$ acts on $\mathcal{O}_h$ as a permutation we can extend $\sigma$ to an element of $\text{Aut}_\quot(\overline{L})$. So $\overline{L}/\quot$ is normal and $\overline{L}_h\subseteq\overline{L}$, concluding our hypothesis.
\end{proof}
Hence from proposition \ref{5.26}, if $\nu^{-1}>2^{d+2}d!\log\Delta_{L_h/\quot}$, then $1-\nu$ must be a zero of some quadratic field $F_h=\quot(s_h)\subseteq L_h$. By assumption $2\nmid [K:\quot]$, hence $K$ has no quadratic subfields, which means that $F_h\cap K=\quot$ and so $F_h$ is the only quadratic subfield of $L_h$. Moreover, $L_h$ is the smallest field containing both $F_h$ and $K$ and since the $h\in H$ are not related by squares of $K$ it holds that $L_h$ doesn't contain any $h'\in H$ from another class. This assures that the $L_h$ produced as $h$ varies are all distinct fields and so all $s_h$ must be distinct. \\
\\
By the transitivity of the discriminant, lemma \ref{2.47}, we have, for the towers of number field $L_h/F_h/\quot$ and $L_h/K/\quot$ that
$$\Delta_{F_h/\quot}^dN_{F_h/\quot}(\Delta_{L_h/F_h})=\Delta_{K/\quot}^2N_{K/\quot}(\Delta_{L_h/K}).$$
As $\Delta_{L_h/K}$ is the relative discriminant, hence is an ideal, we can use the Minkowski Bound:
$$N_{K/\quot}(\Delta_{L_h/K})\leq\sqrt{\Delta_{K/\quot}}\left(\frac{4}{\pi}\right)^d\left(\frac{d!}{d^d}\right)\leq\sqrt{\Delta_{K/\quot}}.$$
Since $\Delta_{K/\quot}\leq L_n\left(\frac{4}{3},\frac{1}{2}\delta^2\kappa\right)$ and $\Delta_{L_h/\quot}\leq L_n\left(\frac{4}{3},\frac{5}{4}\delta^2\kappa\right)$, we can conclude that $$\Delta_{F_h/\quot}=\boh\left(L_n\left(\frac{4}{3},\frac{5}{4}\delta\kappa\right)\right).$$
By the contribution of the Euler product of the Dedekind zeta function:
$$\zeta_{F_h/\quot}(s)=\zeta(s)\mathcal{L}\left(s,\left(\frac{\Delta_{F_h/\quot}}{\cdot}\right)\right),$$
and by reciprocity $j\mapsto\left(\frac{\Delta_{F_h/\quot}}{j}\right)$ is a Dirichlet character modulo $\Delta_{F_h/\quot}$.\\
\\
Now we invoke [\cite{RNFS}, fact 6.24]:
\begin{prop}\label{5.28}
There is an effective constant $c>0$ such that given two characters $\chi_r,\chi_{r'}$ of moduli $r,r'$ respectively, with $\chi_r\chi_{r'}$ non-principal, then the product of Dirichlet $\mathcal{L}$-functions $\mathcal{L}(s,\chi_r)\mathcal{L}(s,\chi_{r'})$ has at most one real zero in $\left(1-\frac{c}{\log rr'},1\right)$.
\end{prop}
So if there are two characters modulo $q$ and $q'$ respectively there is at most one $\mathcal{L}$-function with a zero $1-\nu$ and $\nu^{-1}>c\log qq'$ for some effective $c$. It immediately follows that there is at most one character with modulus in $[q,q^e]$ with a zero at $1-\nu$ and $\nu^{-1}>(e+1)c\log q$.\\
\\
Since $\Delta_{F_h/\quot}<L_n\left(\frac{4}{3},\left(\frac{5}{4}+\loh(1)\right)\delta\kappa\right)$ and $\Delta_{F_h/\quot}\in\zahle$ it follows that all possible discriminants can be covered by $\frac{4}{3}\log\log n$ ranges of form $[x,x^e]$. Hence there are at most $\frac{4}{3}\log\log n$ exceptional characters with exceptional zeroes such that
$$\nu^{-1}>(e+1)c\log\left(L_n\left(\frac{4}{3},\left(\frac{5}{4}+\loh(1)\right)\delta\kappa\right)\right)=\boh\left(\delta\kappa\log^{\frac{4}{3}}(n)\log\log(n)^{-\frac{1}{3}}\right),$$
as required. This bound is far weaker than the required bound of $\nu^{-1}>2^{d+2}d!\log\Delta_{L_h/\quot}$, and so there are at most $\frac{4}{3}\log\log n$ extensions $L_h/\quot$ with exceptional zeroes and $\nu^{-1}>2^{d+2}d!\log\Delta_{L_h/\quot}$. Finally:
$$2^{d+2}d!\log\Delta_{L_h/\quot}\leq d^{d(1+\loh(1)}\log^{\boh(1)}(n)=L_n\left(\frac{1}{3},\frac{\delta}{3}\left(1+\loh(1)\right)\right).$$
This completes the proof of the lemma.
\subsubsection{Proof of theorem \ref{5.3}}\label{sec: RNFSth5.3}
We have now removed all obstructions that were needed. It is therefore that we can now prove theorem \ref{5.3}. Recall that for a set of characters $\mathcal{F}$ we define the map
$$\Psi_{\mathcal{F}}:H\rightarrow\mathbb{F}_2^{\abs{F}},$$
$$x\mapsto\left(\chi_{(p,r)}(x):\chi_{(p,r)}\in\mathcal{F}\right).$$
With the work we have done up to now we are ready to produce our linear $\Psi_{\mathcal{F}}$ with small kernel, and this to produce a congruence of squares. Due to the size of some of the numbers involved, having $L_n\left(\frac{1}{3}\right)$ bits, we track the arithmetic very closely. First, using lemma \ref{5.19}, we sample $4d\left(\delta\kappa\log n+\frac{\delta^2\kappa}{2\log 2}\frac{\log^{\frac{4}{3}}(n)}{\log \log n^{\frac{1}{3}}}\right)$ pairs $(p_i,r_i)$ from $\Upsilon$ independently. \\
\\
Note that this sample is of size $\loh(\log^2(n))=L_n\left(\frac{1}{3},\loh(1)\right)$ and each individual sample takes at most $L_n\left(\frac{1}{3},c\right)$, so we can produce the complete sample in $L_n\left(\frac{1}{3},c+\loh(1)\right)$. After this process we have $M=1+B+dB'+4d\left(\delta\kappa\log n+\frac{\delta^2\kappa}{2\log 2}\frac{\log^{\frac{4}{3}}(n)}{\log \log n^{\frac{1}{3}}}\right)$ pairs from Section \ref{sec: RNFSunif} and the sample from $\Upsilon$. Note that $M=L_n\left(\frac{1}{3},\max\left(\beta,\beta'\right)\right)^{1+\loh(1)}$. For each of these, we need to evaluate each of our characters and as each character evaluates in $L_n\left(\frac{1}{3},\frac{c}{2}\right)$ we get that this process takes $L_n\left(\frac{1}{3},\frac{c}{2}+\max(\beta,\beta')\right)^{1+\loh(1)}$.\\
\\
Fix some $h\in H\backslash\{1\}$ such that $h$ is not in the exceptional set of size $\log\log n$. Each map $\chi_{(p_i,r_i)}$ is independent and induces a map in $\text{Hom}(H,\mathbb{F}_2)$ such that 
$$\mathbb{P}(h\not\in\text{ker}\left(\chi_{(p_i,r_i)}\right))\geq\frac{1+\loh(1)}{2d}.$$
As a corollary it follows that:
$$\mathbb{P}(h\in\text{ker}(\Psi_{\mathcal{F}}))\leq\left(1-\frac{1+\loh(1)}{2d}\right)^{M-(1+B+dB')}\leq\abs{H}^{-2+\loh(1)}.$$
Hence by the union bound over the non-trivial elements of $H$ the probability that any of these non-exceptional and non-zero elements is in the kernel is $\loh(1)$. Hence with high probability the kernel of $\Psi_{\mathcal{F}}$ has size at most $\frac{4}{3}\log\log n$.\\
\\
From here we proceed as normal. With the $M$ pairs representing linear polynomials we use a fast kernel finding algorithm for sparse matrices, such as Block--Wiedemann, \cite{WIED}, to find a suitable subset $S_i$ to construct a polynomial $P_i$ in time 
$$\boh(M^2)=L_n\left(\frac{1}{3},2\max(\beta,\beta')(1+\loh(1)\right),$$
such that $P_i(m)$ is a square in $\zahle$ and $P_i(\alpha)$ is a square in $\mathcal{O}_{\quot(\alpha)}$ multiplied by one of at most $\frac{4}{3}\log\log n$ elements of $H$. Repeating the algorithm $l=\frac{4}{3}\log\log n$ times to generate polynomials $P_1,\ldots,P_l$, we are able to guarantee that for some $i<j$, $P_i$ and $P_j$ lie over the same element $h$, hence $P_iP_j$ is a square in $\mathcal{O}_{\quot(\alpha)}$. In what follows we consider the $\binom{l}{2}\sim \frac{8}{9}\left(\log\log n\right)^2$ polynomials separately.\\
\\
Now if $\gamma\in\mathcal{O}_{\quot(\alpha}$ and $\gamma^2\in\zahle[\alpha]$, then $\gamma\cdot f'(\alpha)\in\zahle[\alpha]$. Let $\mathcal{S}=S_i\Delta S_j$ and fix the polynomial
$$P=\left[\frac{\partial f}{\partial x}(x,1)\right]^2\prod_{(a,b)\in\mathcal{S}}(a-bx), \text{ and so } u^2=\left[\frac{\partial f}{\partial x}(x,y)\right](m,1)^2\prod_{a,b\in\mathcal{S}}(a-mb)$$
is a square in $\zahle$. Hence $u$ can be found by taking the product modulo $n$ over all $p<B$ of $p$ raised to half the total order of $p$ in the terms $(a-mb)$ for $(a,b)\in\mathcal{S}$ and multiplying by $f'(m,1)$.\footnote{This technique differs minimally from the technique for $\zahle$ in the GNFS} The reason for computing the square root in this fashion is to ensure we need only $M\log n$ additions and divisions and at most $M\log n$ modular multiplications to compute $u \mod n$ from the exponents. This ensures polynomial running time.\\
\\
Similarly, for at least one of the $\binom{l}{2}$ polynomials considered, there exists $v\in\zahle[\alpha]$ such that
$$v^2=\left[\frac{\partial f}{\partial x}(x,y)\right](\alpha,1)^2\prod_{(a,b\in\mathcal{S}}(a-b\alpha).$$
Using \cite{SQRT} we can compute square roots in the number field to find $v(m,1) \mod n$ in time $\boh(M^2)$. We abuse notation and write $v(m)$ for the element of $\zahle/n\zahle$ obtained by substituting $m$ for $\alpha$. By the definition of $f$ we have $f(m,1)=n$ and so:
\begin{align*}
\begin{split}
v(m)^2 \mod f(m,1) & = \left(\left[\frac{\partial f}{\partial x}(x,y)\right](\alpha,1)^2\prod_{(a,b)\in\mathcal{S}}(a-b\alpha)\mod f(\alpha,1)\right)(m)\\
& =\left[\frac{\partial f}{\partial x}(x,y)\right](m,1)^2\prod_{(a,b)\in\mathcal{S}}(a-mb) \mod(f(m,1))\\
& = u^2 \mod n,
\end{split}
\end{align*}
and so we have constructed a congruence of squares in time
$$L_n\left(\frac{1}{3},\max(2\max(\beta,\beta'),\max(\beta,\beta')+\frac{c}{2},c)\right)^{1+\loh(1)}.$$
As $c\leq\left(\frac{4}{3}+\loh(1)\right)\delta$, we can insist that we have at most $\frac{4}{3}\log\log n$ exceptional values of $h$, and our $f$ lies off a set of probability at most $L_n\left(\frac{2}{3},\frac{\kappa-\delta^{-1}}{3}(1+\loh(1))\right)^{-1}$ the run time bound is as claimed, finishing the proof.
\subsection{Computational efficiency}
The computational efficiency boils down to proving theorem \ref{5.1}. For this we fix $n,\beta,\beta',\sigma,\delta$, and $\kappa$ satisfying the conditions of equations \ref{RNFSeq1} and \ref{RNFSeq2}. Then by theorem \ref{5.3} we can extract a congruence of squares $\mod n$ from $L_n\left(\frac{1}{3},\max(\beta,\beta')+\loh(1)\right)$ pairs $(a,b)\in\mathcal{X}_{f,m,n}$ for a fixed $\left(m,f\right)$ in expected time
$$L_n\left(\frac{1}{3},2\max\left(\frac{2\delta}{3},\beta,\beta'\right)(1+\loh(1))\right).$$
Theorem \ref{5.2} tells us that a fixed $(m,f)$ and this many pairs $(a,b)\in\mathcal{X}_{f,m,n}$ will be found in expected time
$$L_n\left(\frac{1}{3},\max\left(\beta,\beta'\right)+\frac{\delta^{-1}}{3\beta}(1+\loh(1))+\frac{\sigma\delta+\kappa}{3\beta'}(1+\loh(1))\right).$$
Hence we can run the RNFS to obtain a congruence of squares $\mod n$ with the expected time bounded by
$$L_n\left(\frac{1}{3},\lambda(1+\loh(1)\right), \lambda=\max\left(2\max\left(\frac{2\delta}{3},\beta,\beta'\right),\max\left(\beta,\beta'\right)+\left(\frac{\delta^{-1}}{3\beta}+\frac{\kappa+\sigma\delta}{3\beta'}\right)\right).$$
Having chosen $\beta,\beta',\sigma,\delta,\kappa$ satisfying the conditions of \ref{RNFSeq1} and \ref{RNFSeq2}
we can optimize the constants. Note that increasing the lesser of $\beta$ and $\beta'$ cannot increase $\lambda$ or cause the conditions on the constants to be violated, so we can assume $\beta=\beta'$. We then compute:
$$2\sigma\geq\lambda\geq\min_{\beta,\delta}\left(\beta+\frac{2\delta^{-1}+\sigma\delta+\loh(1)}{3\beta}\right)\geq \min_{\beta}\left(\beta+\frac{\sqrt{8\sigma}+\loh(1)}{3\beta}\right)\geq2\sqrt[4]{\frac{8\sigma}{9}}+\loh(1).$$
Fix any $\epsilon>0$ such that $\epsilon=\loh(1)$. If we take $\beta=\beta'=\sigma=\frac{2\delta}{3}=\sqrt[3]{\frac{8}{9}}+\epsilon$ and $\kappa=\sqrt[3]{\frac{1}{3}}+\epsilon$ then the above are all equalities. Furthermore all the conditions are satisfied, giving $\lambda=2\sqrt[3]{\frac{8}{9}}+\loh(1)$. This proves a runtime of
$$L_n\left(\frac{1}{3},\sqrt[3]{\frac{64}{9}}+\loh(1)\right).$$
\subsection{Riemann and the future of factorization}\label{sec: RNFSriem}
As a final note we want to have a short discussion on some of the ways that the Generalized Riemann Hypothesis could impact the analysis performed on the Randomized number field sieve. As discussed: Until now we have not made any assumption regarding the GRH and only limited assumptions on the heuristics, but if we were to accept the GRH then much of our discussion becomes a lot simpler.\\
\\
One of the main reasons it becomes simpler is because the Generalized Riemann Hypothesis allows us to assume that there are no Siegel zeroes. This makes the discussion around algebraic obstructions a lot simpler. For example we obtain the following:
\begin{prop}\label{5.29}
Conditional on GRH, for $\epsilon=\log^{-\frac{1}{4}}(n)=\loh(\delta)$,
$$\mathbb{P}_{f}\left(\abs{E_{K,\epsilon}}>0\right)=0.$$ 
\end{prop}
This follows automatically from our discussion of the size of this set and the fact there are no Siegel zeroes. This makes lemma \ref{5.22} vacuous and allows us to sharpen the $L_n\left(\frac{1}{3}\right)$ bounds in lemma  \ref{5.19} polynomial in $\log n$.\\
 \\
 The GRH also allows us to have far tighter effective bounds on the prime numbers in arithmetic progressions. Without going into details it would allow us to subvert the use of the adapted Bombieri--Vinogradov theorems and use an effective version of Chebotarev's Density Theorem which is completely dependent on the GRH holding.\\
 \\
 It is obvious that in a discussion like this it could easily be assumed that the GRH holds, but that would make the whole duscussion conditional and the beauty of this is that we make no grand assumptions whatsoever. To obtain a complexity that is equivalent to the heuristic complexity without relying on any conditions makes this the best possible outcome.\\
 \\
 Does this mean we are done with the number field sieve after thirty years? No, definitely not. Over the years many different versions have found their way into mathematics and cryptography to solve other problems than factorization. For example the Tower Number Field Sieve which is an adapted algorithm to solve Discrete Logarithm Problems (\cite{TNFS}, \cite{ENFS}) has been worked on for many years and it will be interesting to see if a similar randomization leads to a provable complexity for that as well.

\end{document}